\documentclass[reqno]{amsart}
\usepackage[colorlinks=true, pdfstartview=FitV, linkcolor=blue, citecolor=blue]{hyperref}
\usepackage[letterpaper, margin=1.3in]{geometry} 

\usepackage{amssymb,amsmath,amsfonts}
\usepackage{bm}
\usepackage{enumitem}
\usepackage{array}
\usepackage{graphicx}
\allowdisplaybreaks

  \newcommand{\R}{\mathbb{R}}  
\newcommand{\Z}{\mathbb{Z}}  \newcommand{\Q}{\mathbb{Q}}  
\newcommand{\rar}{\rightarrow}
\newcommand{\mc}{\mathcal}   \newcommand{\mr}{\mathrm}
\newcommand{\balph}{{\bm{\alpha}}}
\newcommand{\bn}{{\bm{n}}}

\newtheorem*{theorem*}{Theorem}
\newtheorem{theorem}{Theorem}
\newtheorem{lemma}[theorem]{Lemma}
\newtheorem{corollary}[theorem]{Corollary}
\newtheorem{proposition}[theorem]{Proposition}

\begin{document}
\title[Cut-and-project sets with interval windows]{Planar cut-and-project sets with interval windows are bi-Lipschitz equivalent to lattices}
\author{Alan Haynes}

\subjclass[2020]{52C23; 51F30; 11J71}

\keywords{Cut-and-project sets, model sets, Delone sets, separated nets, bi-Lipschitz equivalence}

\begin{abstract}
We prove that every cut-and-project set in a plane in $\R^3$ whose window is a line segment is bi-Lipschitz equivalent to a lattice: if $V\subseteq\R^3$ is a plane which is not spanned by vectors in $\Z^3$, $x\in\R^3$, and $\sigma$ is a segment of any length which is transverse to $V$, then the set of points $v\in V$ with $v+x\in\sigma+\Z^3$ is bi-Lipschitz equivalent to $\Z^2$. For planes of Diophantine slope this is a theorem of Burago and Kleiner, and the point of the present paper is that no Diophantine condition on the plane is needed. We also show that the same conclusion holds when the window is a finite union of intervals. Together with the author's construction of a cut-and-project set in a plane in $\R^3$, with a compact window of countable boundary, which is not bi-Lipschitz equivalent to any lattice, this shows that for planes in $\R^3$ the answer to the bi-Lipschitz question depends on the shape of the window rather than on the Diophantine properties of the plane.
\end{abstract}

\maketitle

\section{Introduction}\label{sec.Intro}

A Delone set in $\R^d$ is \textit{rectifiable} if it is bi-Lipschitz (BL) equivalent to $\Z^d$, that is, if there is a bijection from the set onto $\Z^d$ which is Lipschitz and has a Lipschitz inverse. Whether every Delone set in $\R^d$ is rectifiable was asked by Mikhael Gromov \cite{Grom1993}, in connection with quasi-isometries of groups, and by Hillel Furstenberg, in connection with Kakutani equivalence of $\R^d$-actions (see \cite[Section 3]{SmilSolo2026} for the history). For $d=1$ the answer is yes, for elementary reasons, and for $d\ge2$ it is no: Dmitri Burago and Bruce Kleiner \cite{BuraKlei1998} and, independently, Curtis McMullen \cite{McMu1998} constructed Delone sets in the plane which are not rectifiable, by building a set whose density at a sequence of large scales imitates a bounded function which is not the Jacobian of any BL homeomorphism. Sets built in this way are deliberately irregular, and much of the work since then has concerned Delone sets which come from aperiodic order and are described by a finite amount of data: those arising from primitive substitution tilings are rectifiable \cite{Solo2011}, as are all linearly repetitive Delone sets \cite{AlisCoroGamb2013}, while Cortez and Navas \cite{CortNava2016} and Bhat and Dymond \cite{BhatDymo2026} found repetitive subsets of $\Z^2$ which are not.

The most prominent family left is that of \textit{cut-and-project sets} (also called model sets; see \cite[Chapter 7]{BaakGrim2013}), which are obtained by projecting to a subspace $V\subseteq\R^n$ those points of a lattice whose component orthogonal to $V$ lies in a bounded `window.' Burago and Kleiner asked in \cite[Introduction, Problem 3]{BuraKlei2002} whether every cut-and-project set is rectifiable. This is also Problem 2.8.2 of the list \cite{AdicEtal2016} (a list of open problems from a 2015 workshop on Mathematical Quasicrystals at Oberwolfach), and it appears as Open Problem 7.2 in the recent survey of Smilansky and Solomon \cite{SmilSolo2026}. On the positive side, Haynes, Kelly, and Weiss \cite[Theorem 1.1]{HaynKellWeis2014} showed, in all dimensions, that for almost every subspace $V$ the cut-and-project sets in $V$ are rectifiable whenever the boundary of the window has upper Minkowski dimension smaller than the dimension of the internal space. Koivusalo and Lagacé \cite{KoivLaga2025}, using the criteria of Burago and Kleiner \cite{BuraKlei2002} and of Aliste-Prieto, Coronel, and Gambaudo \cite{AlisCoroGamb2013}, found an explicit full-measure class of rectifiable cut-and-project sets, and their sharp discrepancy bounds show that this is the largest class to which those criteria can be applied. In codimension one, the author showed in \cite{Hayn2016} that for every subspace $V$ there are windows, namely segments whose lengths lie in a certain set depending on $V$, for which the cut-and-project set is even within bounded distance of a lattice. 

In the companion paper \cite{Hayn2026} we showed that the answer to the BL question for cut-and-project sets is in general no, by constructing a cut-and-project set in a plane in $\R^3$, with a window having infinitely many components, which is not rectifiable. The present paper proves a complementary positive result: for a plane in $\R^3$ which is not spanned by vectors in $\Z^3$, and a window which is an interval, the cut-and-project set is always rectifiable, independent of the Diophantine properties of the plane.

To state the result, for $\balph=(\alpha_1,\alpha_2)\in\R^2$ let
\[
V_\balph=\{z\in\R^3 : z_3=\alpha_1z_1+\alpha_2z_2\},\qquad \bn_\balph=(\alpha_1,\alpha_2,-1),
\]
and let $\pi_\balph$ denote the orthogonal projection from $\R^3$ onto $V_\balph$. For a bounded set $W\subseteq\R$ we write
\begin{equation}\label{eqn.CPSDef}
Y_\balph(W)=\left\{\pi_\balph(z) : z\in\Z^3,~\langle z,\bn_\balph\rangle\in W\right\}.
\end{equation}
This is the cut-and-project set with physical space $V_\balph$, internal space $V_\balph^\perp$ (identified with $\R$ by the map $w\mapsto\langle w,\bn_\balph\rangle$), lattice $\Z^3$, and window $W$.

\begin{theorem}\label{thm.Main}
Let $\balph=(\alpha_1,\alpha_2)\in\R^2\setminus\Q^2$, and let $W\subseteq\R$ be an interval (open, closed, or half-open) with $0<|W|<1$. Then $Y_\balph(W)$ is a Delone set in $V_\balph$ which is BL equivalent to $\Z^2$.
\end{theorem}

This is the statement cited in \cite{Hayn2026}. Nothing is assumed about $\balph$ beyond $\balph\notin\Q^2$: if $1,\alpha_1,$ and $\alpha_2$ are linearly dependent over $\Q$ (for instance, if $\alpha_2$ is rational), then $V_\balph$ contains a nonzero integer vector $z$ and $Y_\balph(W)+z=Y_\balph(W)$, so that $Y_\balph(W)$ is periodic in one direction, but the conclusion is the same. The restriction $|W|<1$ is a normalization rather than a hypothesis, as Corollary \ref{cor.CPS} below shows. The interval can also be replaced by a finite union of intervals (see Theorem \ref{thm.Union} in Section \ref{sec.Proofs}).

What was known before is the theorem of Burago and Kleiner \cite[Theorem 1.2]{BuraKlei2002}, which gives the conclusion of Theorem \ref{thm.Main} when $\alpha_1$ satisfies their Diophantine condition (1.1), that is, when $|\alpha_1-p/q|>C/q^{\kappa}$ for some $C>0$ and $\kappa>2$ and all integers $p$ and $q\ge1$ (their theorem is stated for a more general window $B\subseteq\R^3$, and concerns $\alpha_1$ only; see the discussion in the introduction of \cite{Hayn2026}).

Theorem \ref{thm.Main} will be deduced in Section \ref{sec.Reduction} from the following statement about the set
\begin{equation}\label{eqn.LambdaDef}
\Lambda_\balph(W)=\left\{m\in\Z^2 : \langle\balph,m\rangle\in W+\Z\right\},
\end{equation}
which is the setting used in the majority of this paper.

\begin{theorem}\label{thm.Plane}
Let $\balph\in\R^2\setminus\Q^2$, and let $W\subseteq\R$ be an interval with $0<|W|<1$. Then $\Lambda_\balph(W)$ is BL equivalent to $\Z^2$.
\end{theorem}

The passage from $Y_\balph(W)$ to this subset of $\Z^2$ is the one used in \cite{Hayn2016} and \cite{Hayn2026}, and it costs nothing, because the two sets differ by a linear change of coordinates and a bounded perturbation (Lemma \ref{lem.NormalForm}); what is gained is that the internal coordinate of each point is the explicit linear form $\langle\balph,m\rangle$, so that the set can be analyzed one row at a time.

The general cut-and-project statement is the following. We say that a $2$-dimensional subspace $V\subseteq\R^3$ is \textit{rational} if it is spanned by vectors in $\Z^3$, which is the same as saying that its orthogonal complement $V^\perp$ contains a nonzero integer vector.

\begin{corollary}\label{cor.CPS}
Let $V\subseteq\R^3$ be a $2$-dimensional subspace which is not rational, that is, with $V^\perp\cap\Z^3=\{0\}$, let $x\in\R^3$, and let $\sigma\subseteq\R^3$ be a line segment (open, closed, or half-open) with endpoints $z_0\not=z_1$ and $z_1-z_0\notin V$. Then the set
\[
Y=\{v\in V : v+x\in\sigma+\Z^3\}
\]
is BL equivalent to $\Z^2$.
\end{corollary}

For an open segment $\sigma$ the set $Y$ is the cut-and-project set $Y_{\mc{S},x}$ of \cite{Hayn2016}, with section $\mc{S}$ the image of $\sigma$ in $\R^3/\Z^3$, so Corollary \ref{cor.CPS} says that for planes in $\R^3$ neither the length of the segment nor its position matters. In the formulation of \cite[Problem 3]{BuraKlei2002}, in which one projects orthogonally to a plane $P$ the points $z\in\Z^n$ for which $(z+B)\cap P\not=\emptyset$, for a bounded set $B$ with nonempty interior, Corollary \ref{cor.CPS} settles the case $n=3$ for every plane $P$ which is not rational, whenever the orthogonal projection of $B$ to the normal line of $P$ is an interval (for instance, whenever $B$ is connected): the set obtained is then the set $Y$ of the corollary, with $x=0$ and $\sigma$ this projection.

Our proof follows the strategy of Burago and Kleiner \cite{BuraKlei2002}, but it departs from their execution at the first step, and the departure is what removes the Diophantine condition. Their criterion \cite[Theorem 1.3]{BuraKlei2002} requires the deviations of the density of the set, on squares of each dyadic scale, from its limiting value to be summable over the scales. For the sets \eqref{eqn.LambdaDef} the natural bound for this deviation is of order $1/q_k$ at every scale between consecutive denominators $q_k$ and $q_{k+1}$ of the convergents of $\alpha_1$, so each $k$ contributes about $\log_2(q_{k+1}/q_k)$ terms of this size, and the sum may diverge when the partial quotients of $\alpha_1$ grow quickly. We therefore replace the dyadic squares by a hierarchy of rectangular cells adapted to the continued fraction of $\alpha_1$, of width $q_j$ or $q_j+q_{j-1}$ and height comparable to $q_j$, and we pass from the cells of level $j$ to those of level $j+1$ in a single step, however large the ratio $q_{j+1}/q_j$ may be. On a cell of level $j+1$ the density of $\Lambda_\balph(W)$, relative to its average over the cell, is $1+O(1/(q_j|W|))$, by an elementary count of $q_j$ consecutive terms of an orbit of the rotation by $\alpha_1$ (Lemma \ref{lem.Block}). The main step of the paper, Proposition \ref{prop.CellLemma}, is that this relative density is the Jacobian of a homeomorphism of the cell which fixes its boundary and is $(1+O(1/(q_j|W|)))$-BL. Composing these maps over all levels produces a BL homeomorphism of the plane whose constant is at most $\exp(O(\sum_j1/(q_j|W|)))$, and $\sum_j1/q_j$ converges for every irrational $\alpha_1$. This homeomorphism maps each cell of a fixed level onto a set whose area is proportional to the number of points of $\Lambda_\balph(W)$ in the cell, and a theorem of Laczkovich \cite{Lacz1992} on uniformly spread sets, which rests on Hall's marriage theorem, converts this into a bijection with $\Z^2$.

The bulk of our effort is focused on proving Proposition \ref{prop.CellLemma}. The corresponding step of \cite[Section 3]{BuraKlei2002} treats a square divided into four equal squares, and uses no information about the density beyond the ratio of its largest and smallest values. In our hierarchy, by contrast, a cell of level $j+1$ is divided into about $a_{j+1}$ columns and $a_{j+1}$ rows of cells of level $j$, and to obtain a BL constant which does not depend on $a_{j+1}$ we use the arithmetic structure of the relative density. Our maps are built from shears. The relative density is constant on each cell of level $j$, so on each column it is a function of the vertical coordinate alone, and such a function is the Jacobian of a vertical shear whose BL constant depends only on the largest and smallest values of the function, and not on how it oscillates. What has to be controlled is that the shears of adjacent columns do not agree on their common edge. The disagreement is removed in a thin rectangle along the edge, by a map with Jacobian exactly $1$ which is the identity on the other three sides of the rectangle (Lemma \ref{lem.Zone}), and our bound for the cost of this map requires the disagreement to be small compared with the width $q_j$ of a column. The same applies to horizontal shears of the rows. Apart from an error which is easy to control, the values of the relative density are determined by sums of a periodic step function along orbits of the rotation by $q_j\alpha_2$, and, depending on the rational approximations to $q_j\alpha_2$, the disagreements can be too large for vertical shears alone, or for horizontal shears alone. Using a rational approximation to $q_j\alpha_2$, taken from Dirichlet's theorem, we split these sums into two parts, one for which the disagreements between vertical shears are small, and one, varying slowly from one row to the next, for which the disagreements between horizontal shears are small, and we compose the two resulting maps (Sections \ref{sec.Density} and \ref{sec.CellLemma}).

This paper is organized as follows: In Section \ref{sec.Reduction} we collect some preliminary results, and we deduce Theorem \ref{thm.Main} from Theorem \ref{thm.Plane} and Corollary \ref{cor.CPS} from Theorem \ref{thm.Main}. In Section \ref{sec.Hierarchy} we construct the hierarchy of cells, state Proposition \ref{prop.CellLemma}, and show how the maps which it provides compose into a single BL homeomorphism of the plane. Section \ref{sec.Tools} contains the tools which we use to construct BL maps with prescribed Jacobians. In Section \ref{sec.Density} we study the relative density on a cell, and in Section \ref{sec.CellLemma} we prove Proposition \ref{prop.CellLemma}. In Section \ref{sec.Bottom} we apply Laczkovich's theorem to pass from this homeomorphism to a bijection with $\Z^2$. Finally, in Section \ref{sec.Proofs} we complete the proof of Theorem \ref{thm.Plane}, extend it to windows which are finite unions of intervals, and make a remark about higher dimensions.

\noindent\textbf{Notation:} For $x,y\in\R^n$ we write $\langle x,y\rangle$ for the standard inner product and $|x|$ for the Euclidean norm, and for $x\in\R$ we write $\|x\|$ for the distance from $x$ to the nearest integer, $\lfloor x\rfloor$ and $\lceil x\rceil$ for the greatest integer less than or equal to $x$ and the least integer greater than or equal to $x$, and $\{x\}=x-\lfloor x\rfloor$. For a $2\times2$ real matrix $A$ we write $\|A\|=\sup\{|Av|:v\in\R^2,~|v|=1\}$ for its operator norm with respect to the Euclidean norm. For a condition $P$ we write $\mathbf{1}[P]$ for the indicator function of $P$, which is $1$ when $P$ holds and $0$ otherwise. For a measurable set $A$ we write $|A|$ for its Lebesgue measure (length or area, as appropriate), and for a finite set $A$ we write $\#A$ for its cardinality. A map $f$ between subsets of Euclidean spaces is \textit{$L$-BL} if
\[L^{-1}|x-y|\le|f(x)-f(y)|\le L|x-y|\]
for all $x$ and $y$ in its domain, and two subsets of Euclidean spaces are \textit{BL equivalent} if there is a bijection between them which is $L$-BL for some $L$. A set $\Lambda\subseteq\R^2$ is a Delone set if it has \textit{separation} $\inf\{|x-y|:x,y\in\Lambda,~x\not=y\}>0$ and finite \textit{covering radius} $\sup_{u\in\R^2}\inf_{x\in\Lambda}|u-x|$; we say that $\Lambda$ is \textit{uniformly discrete} if its separation is positive and \textit{relatively dense} if its covering radius is finite. For $x\in\R$ and $U\subseteq\R/\Z$ we write $x\in U$ to mean that the image of $x$ in $\R/\Z$ lies in $U$. We use $O(\cdot)$ with the understanding that the implied constant is absolute, that is, independent of $\balph$, $W$, and every index, unless a dependence is indicated.

\noindent\textbf{Disclosure of AI use:} AI was used in this research, as well as to assist in the preparation of this manuscript. However, the final details of all proofs, as well as every statement, citation, and claim in the paper, have been checked by the author.

\section{Reduction and preliminaries}\label{sec.Reduction}

In this section we collect what we need about the set $\Lambda_\balph(W)$ before the main construction begins: a normalization of $\balph$ and $W$, the facts about continued fractions which we use, a counting lemma along rows, and the comparison between $Y_\balph(W)$ and $\Lambda_\balph(W)$ which gives Theorem \ref{thm.Main} from Theorem \ref{thm.Plane}. At the end of the section we deduce Corollary \ref{cor.CPS} from Theorem \ref{thm.Main}.

\noindent\emph{Normalizations.}
It is helpful to picture $\Lambda_\balph(W)$ one row at a time. On the horizontal line $m_2=k$ the set consists of those integers $m_1$ for which $m_1\alpha_1+k\alpha_2$ lies in $W$ modulo $1$, which is to say, the return times to $W$ of the orbit of $k\alpha_2$ under rotation by $\alpha_1$. Three observations let us put $\balph$ and $W$ in a convenient form. First, the map $(m_1,m_2)\mapsto(m_2,m_1)$ is an isometry of $\R^2$ which preserves $\Z^2$ and carries $\Lambda_{(\alpha_1,\alpha_2)}(W)$ onto $\Lambda_{(\alpha_2,\alpha_1)}(W)$, so in proving Theorem \ref{thm.Plane} we may assume that $\alpha_1$ is irrational. Second, $\Lambda_\balph(W)$ depends only on $\alpha_1$ and $\alpha_2$ modulo $1$, and only on $W+\Z$, so we may assume that $\alpha_1\in(0,1)$ and that $W$ is contained in $[0,2)$. Third, the sets $\Lambda_\balph(W)$ and $\Lambda_\balph(W^\circ)$, where $W^\circ$ is the interior of $W$, differ only by points $m$ with $\langle\balph,m\rangle$ in $\partial W+\Z$, and since $\alpha_1$ is irrational there are at most two such points on each horizontal line. In Sections \ref{sec.Hierarchy}-\ref{sec.CellLemma} we work with an open window, and we write
\[
W=(a,b),\qquad 0<|W|=b-a<1 ;
\]
the points on the boundary are accounted for in Section \ref{sec.Bottom}.

\noindent\emph{Continued fractions.}
We write $\alpha_1=[0;a_1,a_2,\ldots]$ for the continued fraction expansion of $\alpha_1$, with convergents $p_k/q_k$, so that $q_{-1}=0$, $q_0=1$, and $q_{k+1}=a_{k+1}q_k+q_{k-1}$ for $k\ge 0$, and we put
\[
D_k=q_k\alpha_1-p_k .
\]
We use the following standard facts (see for example \cite{BereHaynVela2020} or \cite[Section 3]{Hayn2016}). The numbers $D_k$ alternate in sign and decrease in absolute value, and
\begin{equation}\label{eqn.CFFacts}
\frac{1}{q_k+q_{k+1}}<|D_k|<\frac{1}{q_{k+1}},\quad\text{so that}\quad q_k|D_k|<\frac{1}{a_{k+1}}\quad\text{and}\quad\frac12<q_{k+1}|D_k|<1 .
\end{equation}
Also $q_kD_m-q_mD_k=q_mp_k-q_kp_m\in\Z$ for all $k$ and $m$, and $q_{k+2}\ge 2q_k$, so that
\begin{equation}\label{eqn.qSum}
\sum_{j\ge k}\frac{1}{q_j}\le\frac{4}{q_k}\qquad\text{for every}\quad k\ge 0 .
\end{equation}

\noindent\emph{Counting along a row.}
For an integer $M\ge 1$ and $\gamma\in\R$ let
\[
R_M(\gamma)=\#\{0\le n<M : n\alpha_1+\gamma\in W+\Z\}-M|W| ,
\]
the discrepancy of $M$ consecutive terms of the orbit of $\gamma$. Then the number of points of $\Lambda_\balph(W)$ in the block $\{s,\ldots,s+M-1\}\times\{k\}$ of a row is $M|W|+R_M(s\alpha_1+k\alpha_2)$. For blocks whose length is a denominator $q_k$ the discrepancy has a simple form. Let
\[
\tau_k=\{q_k|W|\}\qquad\text{and}\qquad\omega_k(x)=\mathbf{1}\big[\{x\}\ge 1-\tau_k\big]-\tau_k\quad(x\in\R),
\]
so that $\omega_k$ is a $1$-periodic function which takes the value $1-\tau_k$ on an interval of length $\tau_k$ and the value $-\tau_k$ on the rest of each period, and has mean zero. We refer to $\omega_k$ as the \textit{wave} of index $k$.

\begin{lemma}\label{lem.Block}
Let $\alpha_1\in(0,1)$ be irrational, let $W=(a,b)$ with $0<|W|<1$, and let $k\ge 0$. Then, for every $\gamma\in\R$,
\[
R_{q_k}(\gamma)=\omega_k\big(q_k(a-\gamma)\big)+\epsilon_k(\gamma),
\]
where $\epsilon_k(\gamma)$ is an integer with $|\epsilon_k(\gamma)|\le 2$, and $\epsilon_k(\gamma)=0$ unless
\[
\|q_k(a-\gamma)\|<q_k|D_k|\qquad\text{or}\qquad\|q_k(b-\gamma)\|<q_k|D_k| .
\]
In particular $|R_{q_k}(\gamma)|<3$ for all $\gamma$, and the set of $\gamma\in\R/\Z$ for which $\epsilon_k(\gamma)\not=0$ has measure at most $4q_k|D_k|<4/a_{k+1}$.
\end{lemma}

The bound $|R_{q_k}|<3$ is the content of \cite[Lemma 2.3]{BuraKlei2002}, and the proof below is the same alternation argument; what we need in addition is the exact form of the main term.

\begin{proof}
Since $p_k$ and $q_k$ are coprime, the residues of $np_k$ modulo $q_k$, for $0\le n<q_k$, run through all of $\Z/q_k\Z$. Writing
\[
n\alpha_1+\gamma=\Big(\gamma+\frac{np_k}{q_k}\Big)+\frac{nD_k}{q_k},
\]
we see that modulo $1$ the points $n\alpha_1+\gamma$ are the $q_k$ points of the grid $\gamma+\frac{1}{q_k}\Z$, each displaced by an amount $nD_k/q_k$ which has the sign of $D_k$ and absolute value less than $|D_k|$. By \eqref{eqn.CFFacts}, $|D_k|$ is less than the grid spacing $1/q_k$.

First we count the grid points in $W+\Z$. The point $\gamma+r/q_k$ lies in $(a,b)+\Z$ if and only if $r$ lies in $(q_k(a-\gamma),q_k(b-\gamma))+q_k\Z$, and since this interval has length $q_k|W|<q_k$, each residue class modulo $q_k$ contributes at most one integer to it. So the number of grid points in $W+\Z$ is the number of integers in the open interval $(x,x+q_k|W|)$, where $x=q_k(a-\gamma)$. If neither endpoint of this interval is an integer, this number is
\[
\lfloor x+q_k|W|\rfloor-\lfloor x\rfloor=\lfloor q_k|W|\rfloor+\big\lfloor\{x\}+\tau_k\big\rfloor=\lfloor q_k|W|\rfloor+\mathbf{1}\big[\{x\}\ge 1-\tau_k\big],
\]
and subtracting $q_k|W|=\lfloor q_k|W|\rfloor+\tau_k$ gives exactly $\omega_k(x)$. If one of the endpoints is an integer, then $\|q_k(a-\gamma)\|=0$ or $\|q_k(b-\gamma)\|=0$, and the number of grid points differs from the expression above by at most $1$.

Next we account for the displacements. Suppose that $D_k>0$, the other case being the mirror image with the roles of $a$ and $b$ exchanged. A grid point can change its membership of $(a,b)+\Z$ under a displacement to the right by less than $|D_k|$ only if it lies in $(a-|D_k|,a]+\Z$, in which case it may enter, or in $(b-|D_k|,b)+\Z$, in which case it may leave. Each of these two sets has length less than the grid spacing, so each contains at most one grid point, and the displacements change the count by $-1$, $0$, or $1$. Moreover a grid point lies in the first set only if $\{q_k(a-\gamma)\}\in[0,q_k|D_k|)$, and in the second only if $\{q_k(b-\gamma)\}\in(0,q_k|D_k|)$, and both of these conditions imply the one in the statement. Adding the two sources of error gives $|\epsilon_k|\le 2$, and since $\omega_k$ takes values in $(-1,1)$ we have that $|R_{q_k}|<3$. Finally, for $0<\eta\le 1/2$ the set of $\gamma\in\R/\Z$ with $\|q_k(a-\gamma)\|<\eta$ has measure $2\eta$, which gives the last assertion (it is trivial when $q_k|D_k|>1/2$).
\end{proof}

\noindent\emph{From $Y_\balph(W)$ to $\Lambda_\balph(W)$.}
Let $\mr{pr}:V_\balph\rar\R^2$ be the linear isomorphism $\mr{pr}(z)=(z_1,z_2)$. The first statement of the following lemma is \cite[Lemma 4]{Hayn2026}, and the second is standard.

\begin{lemma}\label{lem.NormalForm}
Let $\balph\in\R^2$, let $W\subseteq\R$ be an interval with $0<|W|<1$, and let $s_W=\sup\{|s|:s\in W\}$. Then there is a bijection $\Psi:Y_\balph(W)\rar\Lambda_\balph(W)$ with
\[
|\mr{pr}(y)-\Psi(y)|\le s_W\qquad\text{for all}\quad y\in Y_\balph(W).
\]
If moreover $\balph\notin\Q^2$, then $Y_\balph(W)$ is a Delone set in $V_\balph$.
\end{lemma}
\begin{proof}
	The first statement is proved in \cite{Hayn2026}. For the second, note that $Y_\balph(W)$ is the set of times $v\in V_\balph$ at which the orbit of $0$ under the action of $V_\balph$ on $\R^3/\Z^3$ visits the image of the segment $\{-s\bn_\balph/|\bn_\balph|^2:s\in W\}$, which is transverse to $V_\balph$. Since $\balph\notin\Q^2$ the plane $V_\balph$ is not rational, so every orbit of this action is dense (the closure of an orbit is a translate of a subtorus whose tangent space is a rational subspace containing $V_\balph$). It follows from \cite[Proposition 2.2]{HaynKellWeis2014} that $Y_\balph(W)$ is a Delone set.
\end{proof}

\noindent\emph{Proof of Theorem \ref{thm.Main}, assuming Theorem \ref{thm.Plane}:}
Let $\balph$ and $W$ be as in Theorem \ref{thm.Main}. By Lemma \ref{lem.NormalForm}, $Y_\balph(W)$ is a Delone set in $V_\balph$, and the map $\mr{pr}(y)\mapsto\Psi(y)$ is a bijection from $\mr{pr}(Y_\balph(W))$ onto $\Lambda_\balph(W)$ which moves no point by more than $s_W$. Both of these sets are uniformly discrete, the first as the image of a uniformly discrete set under a linear isomorphism and the second as a subset of $\Z^2$, so this map is BL \cite[Proposition 1.8]{SmilSolo2026}. By Theorem \ref{thm.Plane} there is a BL bijection from $\Lambda_\balph(W)$ onto $\Z^2$, and composing it with the previous map and with $\mr{pr}$ gives a BL bijection from $Y_\balph(W)$ onto $\Z^2$. This completes the proof of Theorem \ref{thm.Main}, assuming Theorem \ref{thm.Plane}.

\noindent\emph{Proof of Corollary \ref{cor.CPS}, assuming Theorem \ref{thm.Main}:}
Let $V$, $x$, and $\sigma$ be as in the corollary, and write $\xi=z_1-z_0$ for the direction of $\sigma$, so that $\xi\notin V$ and $\sigma=\{z_0+s\xi:s\in I\}$ for an interval $I$ (open, closed, or half-open) with endpoints $0$ and $1$. The proof has three steps: we first put $V$ in the form $V_\balph$, then dispose of the case in which two distinct points of $\sigma$ differ by an integer vector, and finally use an integer change of coordinates to make the window have length less than $1$.

\noindent\emph{Step 1: coordinates.} An invertible integer matrix $g\in\mr{GL}_3(\Z)$ maps $\Z^3$ onto itself, so
\[
gY=\{v\in gV : v+gx\in g\sigma+\Z^3\},
\]
which is the set of the corollary for the data $gV$, $gx$, $g\sigma$; the subspace $gV$ is again not rational, and $g\sigma$ is a segment of the same type as $\sigma$ with direction $g\xi\notin gV$. If $g$ is a permutation matrix it is an isometry, so $gY$ and $Y$ are isometric. Since $V$ is a plane, some coordinate vector $e_i$ does not lie in $V$, and after exchanging the coordinates $i$ and $3$ we may assume that $e_3\notin V$. Then $V=V_\balph$ for a unique $\balph\in\R^2$, and $\balph\notin\Q^2$. Write $\bn=\bn_\balph$. Since $\xi\notin V$ we have $\langle\xi,\bn\rangle\not=0$.

\noindent\emph{Step 2: a segment containing a period.} Suppose that $\sigma$ contains two points $p\not=p'$ with $p'-p\in\Z^3$. Then $p'-p$ is a nonzero integer multiple of a primitive vector $z^\ast\in\Z^3$ parallel to $\xi$, so $z^\ast\notin V$, and the segment $[p,p']\subseteq\sigma$ contains a translate of $[0,z^\ast]$. It follows that $\sigma+\Z^3$ contains the whole line $\ell=z_0+\R\xi$, and hence equals $\ell+\Z^3$. Each of the lines $z+\ell$, $z\in\Z^3$, meets the plane $V+x$ in exactly one point, so
\[
Y=\{\rho(z_0+z-x):z\in\Z^3\}=\rho(z_0-x)+\rho(\Z^3),
\]
where $\rho:\R^3\rar V$ is the linear projection onto $V$ along $z^\ast$. Completing $z^\ast$ to a basis $z^\ast,z',z''$ of $\Z^3$, we have $\rho(\Z^3)=\Z\rho(z')+\Z\rho(z'')$.
This implies that $Y$ is a translate of a lattice in $V$, so it is BL equivalent to $\Z^2$. Therefore, from now on we assume that no two distinct points of $\sigma$ differ by an integer vector.

\noindent\emph{Step 3: a short window.} If $\sigma+\Z^3$ contains the line $z_0+\R\xi$, then the argument of Step 2 shows that $Y$ is a translate of a lattice in $V$. Otherwise $\sigma$ is contained in an open segment, parallel to $\xi$, whose image in $\R^3/\Z^3$ is a section as in \cite[Proposition 2.2]{HaynKellWeis2014} for the action of $V$, which is minimal since $V$ is not rational. In either case $Y$ is uniformly discrete.

Now choose an integer $N>|\alpha_1|+|\langle\xi,\bn\rangle|$, and let $g\in\mr{GL}_3(\Z)$ be the inverse of the integer matrix with columns $e_2$, $e_3$, and $(1,0,N)^T$, which has determinant $1$. By Step 1, $gY$ is the set of the corollary for the data $gV,gx,g\sigma$, and $gV$ is not rational. We have $e_3\in gV$ if and only if $g^{-1}e_3=(1,0,N)^T\in V$, and $\langle(1,0,N)^T,\bn\rangle=\alpha_1-N\not=0$, so $e_3\notin gV$, and by Step 1 again $gV=V_{\balph'}$ for some $\balph'\notin\Q^2$. Write $\bn'=\bn_{\balph'}$. The linear functionals $z\mapsto\langle g^{-1}z,\bn\rangle$ and $z\mapsto\langle z,\bn'\rangle$ both have kernel $gV$, so they are proportional, and comparing their values $\alpha_1-N$ and $-1$ at $e_3$ gives
\begin{equation}\label{eqn.NewInternal}
\langle gz,\bn'\rangle=\frac{\langle z,\bn\rangle}{N-\alpha_1}\qquad\text{for all}\quad z\in\R^3 .
\end{equation}
Let $W'=\langle gx,\bn'\rangle-\{\langle\zeta,\bn'\rangle:\zeta\in g\sigma\}$. By \eqref{eqn.NewInternal} this is an interval of the same type as $\sigma$ and of length $|\langle\xi,\bn\rangle|/(N-\alpha_1)$, which lies strictly between $0$ and $1$ by the choice of $N$.

Put $K'=\{z\in\Z^3:\langle z,\bn'\rangle\in W'\}$. For $z\in K'$ the line $z-gx+gz_0+\R g\xi$ meets $gV$ in a single point $v_z$, which by the definition of $W'$ has the form $z+\zeta_z-gx$ with $\zeta_z\in g\sigma$. It is straightforward to check that $z\mapsto v_z$ is a bijection from $K'$ onto $gY$ (injectivity is where Step 2 is used), and that $z\mapsto\pi_{\balph'}(z)$ is a bijection from $K'$ onto $Y_{\balph'}(W')$ (here we use $|W'|<1$). The map $\pi_{\balph'}(z)\mapsto v_z$ moves no point by more than $|z-\pi_{\balph'}(z)|+|\zeta_z-gx|$, which is bounded independently of $z$, and both sets are uniformly discrete, by Lemma \ref{lem.NormalForm} and because $Y$ is. Therefore this map is BL \cite[Proposition 1.8]{SmilSolo2026}. Since $Y_{\balph'}(W')$ is BL equivalent to $\Z^2$ by Theorem \ref{thm.Main}, and $g|_V:V\rar gV$ is a linear isomorphism, it follows that $Y$ is BL equivalent to $\Z^2$. This completes the proof of Corollary \ref{cor.CPS}, assuming Theorem \ref{thm.Main}.

\section{The hierarchy}\label{sec.Hierarchy}

From now until the end of Section \ref{sec.CellLemma} we fix $\alpha_1\in(0,1)$ irrational, $\alpha_2\in\R$, and an open interval $W=(a,b)$ with $0<|W|<1$, and we write $\Lambda=\Lambda_\balph(W)$. This section sets up a nested sequence of partitions of the plane into rectangles, one for each sufficiently large denominator $q_j$ of a convergent of $\alpha_1$, such that the number of points of $\Lambda$ in each rectangle of the $j$th partition, divided by its area, differs from $|W|$ by $O(1/q_j)$. The construction follows Burago and Kleiner \cite[Section 3]{BuraKlei2002}, except that their partitions are refined by a fixed factor at every step, whereas ours jump from one denominator $q_j$ to the next.

Throughout we use the small parameter
\begin{equation}\label{eqn.EpsDef}
\varepsilon_j=\frac{1}{q_j|W|} ,
\end{equation}
and an absolute constant $C_0\ge 12$, whose value is fixed in Section \ref{sec.CellLemma}. We let $j_0$ be the least index with $q_{j_0}|W|\ge C_0$, so that $j_0\ge 1$ and $\varepsilon_j\le 1/C_0$ for all $j\ge j_0$.

\noindent\emph{The partitions.}
For each $j\ge j_0$ we will construct a partition $\mc{P}_j$ of $\Z$ into integer intervals of lengths $q_j$ and $q_j+q_{j-1}$, in such a way that $\mc{P}_j$ refines $\mc{P}_{j+1}$, that is, every interval of $\mc{P}_{j+1}$ is a union of intervals of $\mc{P}_j$. The partitions are built from the top down, by cutting each interval of $\mc{P}_{j+1}$ into intervals of $\mc{P}_j$ according to the following rule: an interval of length $q_{j+1}+q_j$ is cut into one interval of length $q_j+q_{j-1}$ followed by $a_{j+1}$ intervals of length $q_j$, and an interval of length $q_{j+1}$ is cut into one interval of length $q_j+q_{j-1}$ followed by $a_{j+1}-1$ intervals of length $q_j$. Both cuts are possible because $q_{j+1}=a_{j+1}q_j+q_{j-1}$, and the result of a cut depends only on $j$ and on the length of the interval being cut.

Since there is no coarsest partition from which to start, we first carry out the construction on finite initial segments of $\Z_{\ge0}$. For $J\ge j_0$ let $I_J=\{0,\ldots,q_J+q_{J-1}-1\}$ and $\mc{P}^{(J)}_J=\{I_J\}$, and for $j=J-1,J-2,\ldots,j_0$ let $\mc{P}^{(J)}_j$ be the partition of $I_J$ obtained by cutting each interval of $\mc{P}^{(J)}_{j+1}$ according to the rule above. The index $J$ is not a level of the construction, as it only determines how much of $\Z_{\ge0}$ is covered. The constructions for different values of $J$ are compatible: cutting $I_{J+1}$, which has length $q_{J+1}+q_J$, produces $I_J$ as its first piece, so $I_J$ is an interval of $\mc{P}^{(J+1)}_J$, and since the rule depends only on lengths, the intervals of $\mc{P}^{(J+1)}_j$ contained in $I_J$ are exactly the intervals of $\mc{P}^{(J)}_j$, for $j_0\le j\le J$. As the sets $I_J$ increase to $\Z_{\ge0}$, it follows that, for each $j\ge j_0$, the union over $J\ge j$ of the sets $\mc{P}^{(J)}_j$ is a partition $\mc{P}_j$ of $\Z_{\ge0}$. Every interval of $\mc{P}_{j+1}$ belongs to $\mc{P}^{(J)}_{j+1}$ for some $J$, and is therefore cut by the rule into intervals of $\mc{P}^{(J)}_j$, which belong to $\mc{P}_j$; so $\mc{P}_j$ refines $\mc{P}_{j+1}$. Finally, we extend each $\mc{P}_j$ to $\Z_{<0}$ by the reflection $n\mapsto -1-n$. Note that an interval of $\mc{P}_{j+1}$ contained in $\Z_{\ge0}$ begins with its interval of length $q_j+q_{j-1}$, while one contained in $\Z_{<0}$ ends with it.

\noindent\emph{Cells.}
For $j\ge j_0$ let $h_j$ be the power of $2$ with $q_j\le h_j<2q_j$. To each interval of $\mc{P}_j$, which we write as $\{s,\ldots,s+\ell-1\}$ with left endpoint $s\in\Z$ and length $\ell\in\{q_j,q_j+q_{j-1}\}$, and to each $t\in\Z$, we associate the half-open rectangle
\[
Q=[s,s+\ell)\times[th_j,(t+1)h_j)\subseteq\R^2 .
\]
We call these rectangles the \textit{cells of level $j$}. Thus a cell of level $j$ has width $\ell$ and height $h_j$, and its sides have ratio between $1/2$ and $2$. The cells of level $j$ tile $\R^2$, and every cell of level $j+1$ is a union of cells of level $j$, because $\mc{P}_{j+1}$ is coarser than $\mc{P}_j$ and $h_j$ divides $h_{j+1}$. We refer to the cells of level $j$ inside a cell $Q$ of level $j+1$ as the \textit{sub-cells} of $Q$; they are arranged in columns of $h_{j+1}/h_j$ sub-cells each, and there are $a_{j+1}$ or $a_{j+1}+1$ columns, of which the first (if $Q\subseteq\{x_1\ge 0\}$) or the last (if $Q\subseteq\{x_1<0\}$) is the \textit{merged} column of width $q_j+q_{j-1}$, and all others have width $q_j$.

\noindent\emph{Densities.} 
For $j\ge j_0$ let $u_j:\R^2\rar(0,\infty)$ be the function which on each cell $Q$ of level $j$ is equal to the constant
\[
u_j(Q)=\frac{\#(\Lambda\cap Q)}{|Q|} .
\]
Since counts are additive, $u_{j+1}$ on a cell $Q$ of level $j+1$ is the area-weighted average of the values of $u_j$ on the sub-cells of $Q$. On a cell $Q=[s,s+\ell)\times[th_j,(t+1)h_j)$ of level $j$, each row $m_2=k$ of $Q$ contributes $q_j|W|+R_{q_j}(s\alpha_1+k\alpha_2)$ points of $\Lambda$ if $\ell=q_j$, and $(q_j+q_{j-1})|W|+R_{q_j}(s\alpha_1+k\alpha_2)+R_{q_{j-1}}((s+q_j)\alpha_1+k\alpha_2)$ points if $\ell=q_j+q_{j-1}$. By Lemma \ref{lem.Block}, in either case
\begin{equation}\label{eqn.DensityBound}
\big|u_j(Q)-|W|\big|<\frac{6}{q_j}=6|W|\varepsilon_j .
\end{equation}
For $j\ge j_0$ we define the \textit{relative density}
\[
v_{j+1}=\frac{u_j}{u_{j+1}} ,
\]
which is constant on cells of level $j$. If $Q$ is a cell of level $j+1$, then
\[
\int_Qv_{j+1}=\frac{1}{u_{j+1}(Q)}\int_Qu_j=\frac{\#(\Lambda\cap Q)}{u_{j+1}(Q)}=|Q| ,
\]
so $v_{j+1}$ has mean $1$ over every cell of level $j+1$, and from \eqref{eqn.DensityBound}, writing $x=6\varepsilon_j\le 1/2$, we have that $(1-x)/(1+x)<v_{j+1}<(1+x)/(1-x)$, hence
\begin{equation}\label{eqn.RelDensityBound}
|v_{j+1}-1|<24\varepsilon_j .
\end{equation}

\noindent\emph{The cell lemma.} 
The main work of the paper is to show that $v_{j+1}$ is the Jacobian of a map of each cell of level $j+1$ which is close to the identity, in the following sense.

\begin{proposition}\label{prop.CellLemma} 
There are absolute constants $C_0\ge 12$ and $C\ge 1$ with the following property. Let $j\ge j_0=j_0(C_0)$, and let $Q$ be a cell of level $j+1$. Then there is a homeomorphism $\psi_Q$ of $\overline{Q}$ onto itself which fixes $\partial Q$ pointwise, which is $(1+C\varepsilon_j)$-BL, and whose Jacobian satisfies $|\det D\psi_Q|=v_{j+1}$ almost everywhere on $Q$.
\end{proposition}

The proof occupies Sections \ref{sec.Tools}-\ref{sec.CellLemma}, and the constant $C_0$ is the threshold which makes every estimate there apply. Here we explain how Proposition \ref{prop.CellLemma} is used. For a bi-Lipschitz map $f$ of a region of the plane, $Df$ exists almost everywhere by Rademacher's theorem, $f$ maps null sets to null sets, and $|f(E)|=\int_E|\det Df|$ for every measurable set $E$ in its domain (the area formula; see \cite[Theorem 3.9]{EvanGari2015}). We use these facts without further comment.

\begin{lemma}\label{lem.Glue} 
For $j\ge j_0$ let $\psi_{j+1}:\R^2\rar\R^2$ be the map which agrees on each cell $Q$ of level $j+1$ with the map $\psi_Q$ of Proposition \ref{prop.CellLemma}. Then $\psi_{j+1}$ is a $(1+C\varepsilon_j)$-BL homeomorphism of $\R^2$ which maps every cell of level $j+1$ onto itself, and $|\det D\psi_{j+1}|=v_{j+1}$ almost everywhere.
\end{lemma}
\begin{proof}
The closures of the cells of level $j+1$ form a locally finite cover of $\R^2$ by rectangles, and the maps $\psi_Q$ agree on common boundaries, where they are the identity. So $\psi_{j+1}$ is well defined and continuous, it maps each cell onto itself, and its inverse, which is $\psi_Q^{-1}$ on each cell, is continuous for the same reason. Let $L=1+C\varepsilon_j$ and let $x,y\in\R^2$. The segment from $x$ to $y$ meets finitely many cells, so it is divided by the cell boundaries into finitely many segments $[x_i,x_{i+1}]$, each contained in the closure of one cell, and therefore
\[
|\psi_{j+1}(x)-\psi_{j+1}(y)|\le\sum_i|\psi_{j+1}(x_i)-\psi_{j+1}(x_{i+1})|\le L\sum_i|x_i-x_{i+1}|=L|x-y| .
\]
The same argument applies to $\psi_{j+1}^{-1}$, so $\psi_{j+1}$ is $L$-BL. The cell boundaries form a null set, so the Jacobian statement follows from the one for the maps $\psi_Q$.
\end{proof}

\noindent\emph{Composition over all levels.} 
Composing the maps of Lemma \ref{lem.Glue} over all levels $j\ge j_0$ produces a single map which spreads the points of $\Lambda$ evenly: it maps each cell of level $j_0$ onto a region whose area is proportional to the number of points of $\Lambda$ in the cell, as in \eqref{eqn.PhiCount} below. An important point is that the product of the BL constants converges, because $\sum_j\varepsilon_j<\infty$ for every irrational $\alpha_1$: this is the reason why no Diophantine condition is needed.

\begin{proposition}\label{prop.Telescope} 
There is an $L_0$-BL homeomorphism $\Phi$ of $\R^2$, where $L_0=e^{4C/C_0}$ is an absolute constant, such that
\begin{equation}\label{eqn.PhiCount}
|\Phi(E)|=\frac{\#(\Lambda\cap E)}{|W|}
\end{equation}
for every finite union $E$ of cells of level $j_0$.
\end{proposition}
\begin{proof}
For $i\ge j_0$ put $\Phi_i=\psi_i\circ\psi_{i-1}\circ\cdots\circ\psi_{j_0+1}$, where the $\psi_{j+1}$ are the maps of Lemma \ref{lem.Glue}, with $\Phi_{j_0}$ the identity. Since $\psi_{i'}$ maps every cell of level $i'$ onto itself, it maps every cell of each coarser level onto itself too, and so $\Phi_i$ maps every cell of level $i$ onto itself. As $v_{i+1}$ is constant on cells of level $i$, this gives $v_{i+1}\circ\Phi_i=v_{i+1}$. By the chain rule, which holds almost everywhere for compositions of bi-Lipschitz maps,
\[
|\det D\Phi_J|=\prod_{i=j_0}^{J-1}v_{i+1}\circ\Phi_i=\prod_{i=j_0}^{J-1}\frac{u_i}{u_{i+1}}=\frac{u_{j_0}}{u_J}\qquad\text{almost everywhere}.
\]
Also $\Phi_J$ is $L_J$-BL with
\[
\log L_J\le\sum_{j=j_0}^{J-1}\log(1+C\varepsilon_j)\le C\sum_{j\ge j_0}\varepsilon_j\le\frac{4C}{q_{j_0}|W|}\le\frac{4C}{C_0},
\]
using \eqref{eqn.qSum} and the definition of $j_0$. So all of the maps $\Phi_J$ are $L_0$-BL.

The origin is a corner of a cell of every level, so it is fixed by each of the maps $\psi_{j+1}$, and $\Phi_J(0)=0$ for every $J$. The maps $\Phi_J$ are therefore uniformly Lipschitz and uniformly bounded on compact sets, and by the Arzel\`{a}--Ascoli theorem there is a subsequence along which they converge, uniformly on compact sets, to a continuous map $\Phi$ of $\R^2$ into itself. The inequalities $L_0^{-1}|x-y|\le|\Phi_J(x)-\Phi_J(y)|\le L_0|x-y|$ pass to the limit, so $\Phi$ is $L_0$-BL onto its image. To see that $\Phi$ is onto, let $y\in\R^2$ and put $x_J=\Phi_J^{-1}(y)$. Since $\Phi_J^{-1}$ is $L_0$-Lipschitz and fixes $0$, we have that $|x_J|\le L_0|y|$, so along a further subsequence, which may depend on $y$, the points $x_J$ converge to some $x\in\R^2$. Then
\[
|\Phi(x)-y|\le|\Phi(x)-\Phi(x_J)|+|\Phi(x_J)-\Phi_J(x_J)|\le L_0|x-x_J|+\sup_{|z|\le L_0|y|}|\Phi(z)-\Phi_J(z)| ,
\]
which tends to $0$ along this subsequence, so $\Phi(x)=y$. Therefore $\Phi$ is a bijection of $\R^2$ onto itself which is $L_0$-BL, and in particular it is a homeomorphism.

Now let $Q$ be a cell of level $j_0$, and let $\nu_J=\sup_Q|\Phi_J-\Phi|$, which tends to $0$ along the subsequence. We claim that the symmetric difference of $\Phi_J(Q)$ and $\Phi(Q)$ lies within distance $2\nu_J$ of the curve $\Phi(\partial Q)$. If $y=\Phi_J(x)$ with $x\in Q$ and $y\notin\Phi(Q)$, then the segment from $y$ to $\Phi(x)\in\Phi(Q)$, which has length at most $\nu_J$, must cross $\partial\Phi(Q)=\Phi(\partial Q)$. If $y=\Phi(x)$ with $x\in Q$ and $y\notin\Phi_J(Q)$, then in the same way the segment from $y$ to $\Phi_J(x)$ crosses $\Phi_J(\partial Q)$, which lies within $\nu_J$ of $\Phi(\partial Q)$. This proves the claim. The curve $\Phi(\partial Q)$ has length at most $L_0|\partial Q|$, and the $2\nu$-neighborhood of a curve of length $\lambda$ has area $O(\nu\lambda+\nu^2)$, so $|\Phi_J(Q)|\rar|\Phi(Q)|$. On the other hand, for $J>j_0$ let $Q^{(J)}$ be the cell of level $J$ which contains $Q$. Since $u_{j_0}$ is constant on $Q$ and $u_J$ is constant on $Q^{(J)}$, the area formula and the formula for $|\det D\Phi_J|$ above give
\[
|\Phi_J(Q)|=\int_Q\frac{u_{j_0}}{u_J}=\frac{|Q|\,u_{j_0}(Q)}{u_J(Q^{(J)})}=\frac{\#(\Lambda\cap Q)}{u_J(Q^{(J)})}\rar\frac{\#(\Lambda\cap Q)}{|W|}\qquad\text{as}\quad J\rar\infty,
\]
where we have used \eqref{eqn.DensityBound}. Therefore $|\Phi(Q)|=\#(\Lambda\cap Q)/|W|$. This is \eqref{eqn.PhiCount} for a single cell, and the general case follows by adding over the cells of a finite union, since $\Phi$ is injective.
\end{proof}

\section{Analytic tools}\label{sec.Tools}

In this section we develop the tools from which the map $\psi_Q$ of Proposition \ref{prop.CellLemma} will be assembled in Section \ref{sec.CellLemma}. Every statement concerns maps of a closed rectangle $Q=[0,w]\times[0,H]$, and each has a mirror image, with the roles of the two coordinates exchanged, which we use without comment.

If a homeomorphism $\varphi$ of $Q$ maps the bottom side of $Q$ onto itself, then $\varphi(x,0)=(g(x),0)$ for a homeomorphism $g$ of $[0,w]$, and we say that $\varphi$ \textit{acts on the bottom side as} $g$; we use the same terminology for the top side, where $\varphi(x,H)=(g(x),H)$, and for the left and right sides, where $\varphi(0,y)=(0,g(y))$ and $\varphi(w,y)=(w,g(y))$ with $g$ a homeomorphism of $[0,H]$. Throughout the section $\varepsilon_\ast\in(0,1/60]$ is an absolute constant, chosen small enough for the proof of Lemma \ref{lem.Zone}; in Section \ref{sec.CellLemma} the threshold $C_0$ is chosen so that every $\varepsilon$ to which these results are applied is at most $\varepsilon_\ast$. We also use two facts about BL maps:

\begin{enumerate}[topsep=2pt,itemsep=4pt,parsep=4pt]
	\item[(F1)] Let $K\subseteq\R^2$ be convex, and let $K=K_1\cup\cdots\cup K_n$ with each $K_i$ closed and convex. If $\varphi:K\rar\R^2$ is $L$-Lipschitz on each $K_i$, then $\varphi$ is $L$-Lipschitz on $K$, because every segment in $K$ is divided into finitely many segments, each contained in one of the $K_i$, and the triangle inequality applies as in the proof of Lemma \ref{lem.Glue}. In particular, if $\varphi$ is a homeomorphism of $K$ onto itself which maps each $K_i$ onto itself and is $L$-BL on each $K_i$, then $\varphi^{-1}$ is $L$-Lipschitz on each $K_i$, and applying the previous statement to $\varphi$ and to $\varphi^{-1}$ shows that $\varphi$ is $L$-BL.
	\item[(F2)] Let $K\subseteq\R^2$ be convex, and let $\varphi:K\rar\R^2$ be Lipschitz with $\|D\varphi-I\|\le\varepsilon\le1/2$ almost everywhere (operator norm). Then $(1-\varepsilon)|z-z'|\le|\varphi(z)-\varphi(z')|\le(1+\varepsilon)|z-z'|$ for all $z,z'\in K$, so $\varphi$ is injective and $(1+2\varepsilon)$-BL onto its image. To see this, note that, by an argument using Rademacher's theorem with Fubini's theorem, for almost every pair $(z,z')$, $\varphi$ is absolutely continuous along the segment $[z,z']$ with derivative $D\varphi$ at almost every point of it, so that \[\big|\varphi(z')-\varphi(z)-(z'-z)\big|=\big|\int_0^1\big(D\varphi(z+t(z'-z))-I\big)(z'-z)\,\mr{d}t\big|\le\varepsilon|z'-z|.\] The general case then follows by continuity.
\end{enumerate}

\noindent\emph{Shears.}
The next two lemmas record what a shear does, first on a rectangle and then on a rectangle cut into vertical strips.

\begin{lemma}\label{lem.Shear}
Let $Q=[0,w]\times[0,H]$, let $L\ge1$, and let $P$ be an increasing $L$-BL homeomorphism of $[0,H]$ onto itself. Then $\varphi(x,y)=(x,P(y))$ is an $L$-BL homeomorphism of $Q$ onto itself with $|\det D\varphi|=P'(y)$ almost everywhere. If $0<\varepsilon\le1/2$ and $f:[0,H]\rar[1-\varepsilon,1+\varepsilon]$ is measurable with $\int_0^Hf=H$, then $P(y)=\int_0^yf$ satisfies these hypotheses with $L=1+2\varepsilon$, with $P'=f$ almost everywhere and $|P(y)-y|\le\varepsilon H$ for all $y$.
\end{lemma}
\begin{proof}
An increasing $L$-BL homeomorphism $P$ of $[0,H]$ onto itself is absolutely continuous with $P'\in[L^{-1},L]$ almost everywhere. Since
\[|\varphi(z)-\varphi(z')|^2=(z_1-z_1')^2+(P(z_2)-P(z_2'))^2\]
and $L\ge1$, the map $\varphi$ is $L$-BL; it is a bijection of $Q$ onto itself because $P$ is a bijection, and $D\varphi=\mr{diag}(1,P'(y))$ wherever $P$ is differentiable. If $f$ is as stated, then $P(y)=\int_0^yf$ is strictly increasing and absolutely continuous with $P(0)=0$, $P(H)=H$, and $P'=f$ almost everywhere, so it is an increasing homeomorphism of $[0,H]$ which is $\max\{1+\varepsilon,(1-\varepsilon)^{-1}\}$-BL, and $(1-\varepsilon)^{-1}\le1+2\varepsilon$. Finally $|P(y)-y|=|\int_0^y(f-1)|\le\varepsilon H$.
\end{proof}

We call a map $\varphi$ as in the statement of Lemma \ref{lem.Shear} a \textit{shear} with \textit{profile} $P$. A shear fixes the bottom and top sides of $Q$ pointwise and acts on each vertical side as its profile.

\begin{lemma}\label{lem.StripShear}
Let $Q=[0,w]\times[0,H]$ be cut into the strips $Q_i=[s_i,s_{i+1}]\times[0,H]$, $0\le i\le N$, where $0=s_0<s_1<\cdots<s_{N+1}=w$, let $0<\varepsilon\le1/2$, and let $u:Q\rar[1-\varepsilon,1+\varepsilon]$ be measurable with $\int_Qu=|Q|$ and $u(x,y)=u^{(i)}(y)$ on each $Q_i$. Let $\bar u^{(i)}=H^{-1}\int_0^Hu^{(i)}$, let $\bar u(x)=\bar u^{(i)}$ for $s_i\le x<s_{i+1}$, and put
\begin{align*}
P_\perp(x)&=\int_0^x\bar u,\qquad&\varphi_\perp(x,y)&=(P_\perp(x),y),\\
P_i(y)&=\frac{1}{\bar u^{(i)}}\int_0^yu^{(i)},\qquad&\varphi_i(x,y)&=(x,P_i(y)).
\end{align*}
Then $\varphi_\perp$ is a $(1+2\varepsilon)$-BL homeomorphism of $Q$ onto itself with $|\det D\varphi_\perp|=\bar u(x)$ almost everywhere, which maps each $Q_i$ onto $Q_i'=[s_i',s_{i+1}']\times[0,H]$, where $s_i'=P_\perp(s_i)$. Each $\varphi_i$ is a $(1+4\varepsilon)$-BL shear of $Q_i'$ onto itself, with profile $P_i$ and $|\det D\varphi_i|=u^{(i)}(y)/\bar u^{(i)}$ almost everywhere. The map $\varphi$ which equals $\varphi_i\circ\varphi_\perp$ on $Q_i$ satisfies $|\det D\varphi|=u$ almost everywhere on $Q$.
\end{lemma}
\begin{proof}
We have
\[\int_0^w\bar u=\sum_i(s_{i+1}-s_i)\bar u^{(i)}=H^{-1}\int_Qu=w,\]
so $\varphi_\perp$ is the shear of Lemma \ref{lem.Shear} with the coordinates exchanged and $f=\bar u$, and $P_\perp$ carries $[s_i,s_{i+1}]$ onto $[s_i',s_{i+1}']$. The function $u^{(i)}/\bar u^{(i)}$ has mean $1$ over $[0,H]$ and takes values in $[(1-\varepsilon)/(1+\varepsilon),(1+\varepsilon)/(1-\varepsilon)]$, so $P_i$ is an increasing homeomorphism of $[0,H]$ onto itself which is $(1+\varepsilon)/(1-\varepsilon)$-BL, and $(1+\varepsilon)/(1-\varepsilon)\le1+4\varepsilon$ for $\varepsilon\le1/2$. Lemma \ref{lem.Shear}, applied on $Q_i'$ with $L=1+4\varepsilon$, gives the statement about $\varphi_i$. By the chain rule, which holds almost everywhere for a composition of bi-Lipschitz maps, $|\det D(\varphi_i\circ\varphi_\perp)|=\big(u^{(i)}(y)/\bar u^{(i)}\big)\cdot\bar u^{(i)}=u$ on $Q_i$.
\end{proof}

The map $\varphi$ of Lemma \ref{lem.StripShear} is in general not continuous: the shears $\varphi_i$ and $\varphi_{i+1}$ act on the common edge $x=s_{i+1}'$ of $Q_i'$ and $Q_{i+1}'$ as $P_i$ and $P_{i+1}$, so the one-sided limits of $\varphi$ at a point $(s_{i+1},y)$ differ by the vertical displacement $P_i(y)-P_{i+1}(y)$, and $\varphi$ is a homeomorphism of $Q$ onto itself if and only if $P_i=P_{i+1}$ for every $0\le i<N$. We call
\[
\mu_i=\sup_{0\le y\le H}|P_i(y)-P_{i+1}(y)|
\]
the \textit{mismatch} at that edge.

\noindent\emph{The zone lemma.}
A mismatch at an edge is removed by composing the shear on one side of the edge with a map of a zone along the edge which fixes three sides of the zone, acts on the fourth side as a given homeomorphism (the correction which makes the maps on the two sides of the edge agree), and has Jacobian exactly $1$. The next lemma shows that if the zone has width $\delta$, and the given homeomorphism is $(1+\varepsilon)$-BL and moves no point by more than $\varepsilon\delta$, then such a map can be taken to be $(1+O(\varepsilon))$-BL.

\begin{lemma}\label{lem.Zone}
Let $\delta>0$, $H\ge5\delta$, $Z=[0,\delta]\times[0,H]$, and $0<\varepsilon\le30\varepsilon_\ast$. Let $\vartheta$ be an increasing homeomorphism of $[0,H]$ onto itself which is $(1+\varepsilon)$-BL and satisfies $|\vartheta(y)-y|\le\varepsilon\delta$ for all $y$. Then there is a $(1+O(\varepsilon))$-BL homeomorphism $\Upsilon$ of $Z$ onto itself such that $|\det D\Upsilon|=1$ almost everywhere, such that $\Upsilon$ is the identity on the three sides $x=0$, $y=0$, and $y=H$ of $Z$, such that
\[
\Upsilon(\delta,y)=(\delta,\vartheta(y))\qquad\text{for all}\quad 0\le y\le H ,
\]
and such that $|\Upsilon(z)-z|=O(\varepsilon\delta)$ for all $z\in Z$.
\end{lemma}

We call $Z$ a \textit{zone} of width $\delta$ and $\vartheta$ the \textit{trace} prescribed on its right side. Composing with the reflection $(x,y)\mapsto(\delta-x,y)$, which preserves area, we may prescribe the trace on the left side instead; and the mirror image of the lemma, with the coordinates exchanged, concerns a zone $[0,w]\times[0,\delta]$ with $w\ge5\delta$ and a trace prescribed on its top or bottom side.

\begin{figure}[h]
\caption{The construction in the proof of Lemma \ref{lem.Zone}, computed for $H=5\delta$ and a piecewise linear trace $\vartheta$ with $|\Delta'|\le 0.4$, which is larger than the lemma allows and is chosen so that the deformation is visible. (a) The function $\chi_x(y)-y$ on $Z$ (positive values in red, negative values in blue), with the graph of the defect $\Delta$ to its right. The lower part of the defect is translated upward and the upper part downward, by $\varsigma_i(x)$ and with the coefficient $\varkappa_i(x)$; the solid and dashed curves follow the peak of each part for $i=1$ and $i=2$. By \eqref{eqn.MeanChi} the integral of $\chi_x(y)-y$ over every horizontal line is $0$. (b)-(d) A square grid in $Z$, its image under $\Upsilon_{\mr{v}}$, and its image under $\Upsilon=\Upsilon_{\mr{h}}\circ\Upsilon_{\mr{v}}$. Every cell has the same area in (d) as in (b). (e) The image of the same grid when the translations are omitted, that is, when $\varsigma_1$ and $\varsigma_2$ are replaced by $0$ and the coefficients are kept. The map which moves points horizontally then shifts the points just above and just below each of the lines drawn in red by different amounts; these lines are at the heights where $\Delta'$ jumps.\\}\label{fig.Zone}
\centering
\includegraphics[width=\textwidth]{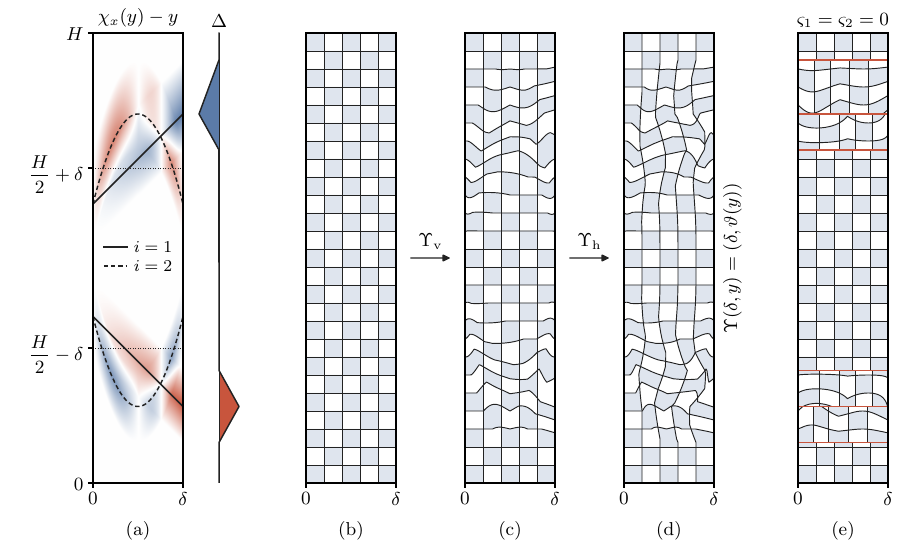}
\end{figure}

Before the proof we describe its structure. The map $\Upsilon$ will be constructed as the composition of a map which moves points vertically, interpolating across the zone between the identity on the left side and $\vartheta$ on the right side, with a map which moves points horizontally and is chosen so that the composite has Jacobian $1$. Write $\Delta=\vartheta^{-1}-\mr{id}$. The obvious choice for the first map is $(x,y)\mapsto(x,y')$, where $y'+\varkappa(x)\Delta(y')=y$, for a function $\varkappa$ with $\varkappa(0)=0$ and $\varkappa(\delta)=1$. For the composite to have Jacobian $1$, the first coordinate of the second map must have $x$-derivative $1+\varkappa(x)\Delta'(y')$, and since the second map has to fix the left side of the zone, it must then be
\[
(x,y')\mapsto\Big(x+\Delta'(y')\int_0^x\varkappa,~y'\Big).
\]
However, $\vartheta$ is only assumed to be BL, so $\Delta'$ need not be continuous (in our applications $\Delta$ is piecewise linear), and then neither is this map (see Figure \ref{fig.Zone}(e)). For this reason we translate the defect vertically as $x$ crosses the zone, replacing $\Delta(y')$ by $\Delta(y'-\varsigma(x))$. The integral above becomes $\int_0^x\varkappa(s)\Delta'(y'-\varsigma(s))\,\mr{d}s$, a weighted average of $\Delta'$ over a range of heights, which is Lipschitz in $y'$ for every bounded $\Delta'$ (see \eqref{eqn.XiLip}). Two adjustments are needed to keep the sides of the zone fixed. To keep the translated defect between the bottom and top sides, we split $\Delta$ into a lower part, which is translated upward, and an upper part, which is translated downward. To make the horizontal map fix the right side of the zone, we add a second translated copy of each part, with a negative coefficient, chosen so that the two copies cancel on average across the zone; this is the identity \eqref{eqn.MeanIdentity} below. Figure \ref{fig.Zone} illustrates the construction.

\begin{proof}
All implied constants in the proof are absolute, and we assume that $\varepsilon_\ast$ is small enough for the finitely many smallness requirements which appear. Let $\Delta=\vartheta^{-1}-\mr{id}$ on $[0,H]$, extended by $0$ to the rest of $\R$, which we call the \textit{defect}. Since $\vartheta^{-1}$ is an increasing $(1+\varepsilon)$-BL homeomorphism of $[0,H]$ it is absolutely continuous with $(\vartheta^{-1})'\in[(1+\varepsilon)^{-1},1+\varepsilon]$ almost everywhere, and $|\vartheta^{-1}(y)-y|=|y'-\vartheta(y')|$ with $y'=\vartheta^{-1}(y)$. So $\Delta$ is a Lipschitz function on $\R$ with $\Delta(0)=\Delta(H)=0$, $|\Delta|\le\varepsilon\delta$, and $|\Delta'|\le\varepsilon$ almost everywhere.

\noindent\emph{Splitting the defect.} Let $\varpi:\R\rar[0,1]$ be the continuous function which is equal to $1$ on $(-\infty,H/2-\delta]$, to $0$ on $[H/2+\delta,\infty)$, and is linear in between, and put $\Delta_-=\varpi\Delta$ and $\Delta_+=(1-\varpi)\Delta$. Then $\Delta=\Delta_-+\Delta_+$, the function $\Delta_-$ vanishes outside $[0,H/2+\delta]$ and $\Delta_+$ vanishes outside $[H/2-\delta,H]$, and since $|\varpi'|\le1/(2\delta)$ we have $|\Delta_\pm|\le\varepsilon\delta$ and $|\Delta_\pm'|\le2\varepsilon$ almost everywhere.

\noindent\emph{A sliding interpolation.} For $0\le s\le\delta$ put $r=r(s)=1-2s/\delta\in[-1,1]$, and let
\[
\varkappa_1(s)=\frac{s}{\delta},\qquad\varsigma_1(s)=\delta-s,\qquad\varkappa_2(s)=-2|r|(1-r^2),\qquad\varsigma_2(s)=\delta r^2 .
\]
Then $0\le\varsigma_1,\varsigma_2\le\delta$, we have $\varkappa_1(0)=\varkappa_2(0)=\varkappa_2(\delta)=0$, $\varkappa_1(\delta)=1$, $\varsigma_1(\delta)=0$, and $|\varkappa_2|<1$, and all four functions are Lipschitz, with $|\varkappa_i'|=O(1/\delta)$ and $|\varsigma_i'|=O(1)$. For $0\le s\le\delta$ and $y\in\R$ let
\[
\chi_s(y)=y+\sum_{i=1}^{2}\varkappa_i(s)\Big[\Delta_-\big(y-\varsigma_i(s)\big)+\Delta_+\big(y+\varsigma_i(s)\big)\Big] .
\]
For each $s$, the lower half of the defect is translated upward by $\varsigma_i(s)$ and the upper half downward, with the coefficient $\varkappa_i(s)$; both coefficients vanish at $s=0$, and at $s=\delta$ the first translation is $0$ and its coefficient is $1$. We note the following properties of the family.
\begin{enumerate}[topsep=2pt,itemsep=4pt,parsep=4pt]
\item[(P1)] $\chi_0=\mr{id}$, and $\chi_\delta=\mr{id}+\Delta$, which agrees with $\vartheta^{-1}$ on $[0,H]$.
\item[(P2)] $\chi_s(y)=y$ for $y\le0$ and for $y\ge H$. If $y\le0$ then $y-\varsigma_i(s)\le0$ and $y+\varsigma_i(s)\le\delta<H/2-\delta$, so all four terms vanish, and if $y\ge H$ then $y-\varsigma_i(s)\ge H-\delta>H/2+\delta$ and $y+\varsigma_i(s)\ge H$; here we use $H\ge5\delta$.
\item[(P3)] $|\chi_s'-1|=O(\varepsilon)$ almost everywhere, so for $\varepsilon_\ast$ small enough each $\chi_s$ is a BL homeomorphism of $\R$ with $\frac12\le\chi_s'\le\frac32$, and by (P2) it maps $[0,H]$ onto itself.
\item[(P4)] $|\chi_s(y)-y|=O(\varepsilon\delta)$ for all $s$ and $y$.
\item[(P5)] The map $(s,y)\mapsto\chi_s(y)$ is Lipschitz on $[0,\delta]\times\R$, and $|\partial_s\chi_s(y)|=O(\varepsilon)$ almost everywhere, because $\partial_s$ of each term is bounded by $|\varkappa_i'||\Delta_\pm|+|\varkappa_i||\varsigma_i'||\Delta_\pm'|$.
\end{enumerate}

\noindent\emph{The mean identity.} The functions $\varkappa_2$ and $\varsigma_2$ are chosen so that, for every bounded measurable function $g$ on $[0,\delta]$,
\begin{equation}\label{eqn.MeanIdentity}
\int_0^\delta\Big[\varkappa_1(s)\,g\big(\varsigma_1(s)\big)+\varkappa_2(s)\,g\big(\varsigma_2(s)\big)\Big]\,\mr{d}s=0 .
\end{equation}
To see this, substitute $c=\delta-s$ in the first term, which becomes $\int_0^\delta(1-c/\delta)g(c)\,\mr{d}c$. On $[0,\delta/2]$ the variable $r$ decreases from $1$ to $0$, with $\mr{d}s=-\frac{\delta}{2}\mr{d}r$, so the second term integrated over $[0,\delta/2]$ is
\[
-\delta\int_0^1r(1-r^2)\,g(\delta r^2)\,\mr{d}r=-\frac12\int_0^\delta\Big(1-\frac{c}{\delta}\Big)g(c)\,\mr{d}c
\]
after the substitution $c=\delta r^2$; the integral over $[\delta/2,\delta]$ has the same value, since $\varkappa_2$ and $\varsigma_2$ are even functions of $r$. The three contributions cancel. Applying \eqref{eqn.MeanIdentity} with $g(c)=\Delta_-(y-c)$ and with $g(c)=\Delta_+(y+c)$ we find that
\begin{equation}\label{eqn.MeanChi}
\int_0^\delta\chi_s(y)\,\mr{d}s=\delta y\qquad\text{for every}\quad y\in\R .
\end{equation}

\noindent\emph{The map.} For $(x,y)\in Z$ let
\[
\Xi(x,y)=x+\sum_{i=1}^{2}\int_0^x\varkappa_i(s)\Big[\Delta_-'\big(y-\varsigma_i(s)\big)+\Delta_+'\big(y+\varsigma_i(s)\big)\Big]\,\mr{d}s .
\]
This is defined for every $y$: the derivatives $\Delta_\pm'$ exist almost everywhere, and the preimage of a null set under $\varsigma_1$ or $\varsigma_2$ is null, since these are polynomials in $s$ whose derivatives vanish at most at one point; so for fixed $y$ the integrand is defined for almost every $s$. We set
\[
\Upsilon=\Upsilon_{\mr{h}}\circ\Upsilon_{\mr{v}},\qquad\Upsilon_{\mr{v}}(x,y)=\big(x,\chi_x^{-1}(y)\big),\qquad\Upsilon_{\mr{h}}(x,y')=\big(\Xi(x,y'),y'\big).
\]
The map $\Upsilon_{\mr{v}}$ is the vertical interpolation, and $\Upsilon_{\mr{h}}$ is the horizontal correction: for fixed $y$, the function $x\mapsto\Xi(x,y)$ is absolutely continuous with
\begin{equation}\label{eqn.XiDeriv}
\partial_x\Xi(x,y)=\chi_x'(y)\qquad\text{for almost every}\quad x\in[0,\delta],
\end{equation}
which is what will make the Jacobian equal to $1$.

\noindent\emph{$\Xi$ is Lipschitz in $y$.} Fix $x$, and consider one of the integrals $\int_0^x\varkappa_i(s)\phi(y\mp\varsigma_i(s))\,\mr{d}s$ with $\phi=\Delta_\mp'$. On each of the subintervals of $[0,x]$ on which $\varsigma_i$ is strictly monotone we substitute $c=\varsigma_i(s)$; the piece becomes $\int_Ig(c)\phi(y\mp c)\,\mr{d}c$, where $I\subseteq[0,\delta]$ is an interval and $g(c)=\varkappa_i(s)/|\varsigma_i'(s)|$, which is $1-c/\delta$ for $i=1$ and $-\frac12(1-c/\delta)$ on each of the two monotone pieces of $\varsigma_2$ for $i=2$, because $|\varsigma_2'|=4|r|$ cancels the factor $|r|$ in $\varkappa_2$. In every case $|g|\le1$ and $g$ is $(1/\delta)$-Lipschitz on $I$. Let $\hat g$ be equal to $g$ on $I$ and to $0$ elsewhere, so that the piece is $\int_\R\phi(u)\hat g(\pm(y-u))\,\mr{d}u$. Since $|g|\le1$, $g$ is $(1/\delta)$-Lipschitz on $I$, and $|I|\le\delta$, we have $\int_\R|\hat g(c+h)-\hat g(c)|\,\mr{d}c\le3|h|$ for every $h\in\R$, and since $|\phi|\le2\varepsilon$ the piece changes by at most $6\varepsilon|y_2-y_1|$ when $y_1$ is replaced by $y_2$. There are at most six pieces in all, so
\begin{equation}\label{eqn.XiLip}
|\Xi(x,y_1)-\Xi(x,y_2)|=O(\varepsilon)|y_1-y_2|\qquad\text{uniformly in}\quad x .
\end{equation}
This is the estimate for which the translations were introduced, and it uses nothing about $\Delta_\pm'$ beyond boundedness.

\noindent\emph{Boundary values of $\Xi$.} We have $\Xi(0,y)=0$. By \eqref{eqn.XiDeriv} and (P3), $x\mapsto\Xi(x,y)$ is an increasing BL map of $[0,\delta]$ onto $[0,\Xi(\delta,y)]$. For $y_1<y_2$, Fubini's theorem and the absolute continuity of $\Delta_\pm$ give
\begin{align*}
\int_{y_1}^{y_2}\big(\Xi(\delta,y)-\delta\big)\,\mr{d}y&=\sum_{i=1}^2\int_0^\delta\varkappa_i(s)\Big[\Delta_-(y-\varsigma_i(s))+\Delta_+(y+\varsigma_i(s))\Big]_{y=y_1}^{y=y_2}\,\mr{d}s\\
&=\int_0^\delta\Big[\chi_s(y)-y\Big]_{y=y_1}^{y=y_2}\,\mr{d}s=0
\end{align*}
by \eqref{eqn.MeanChi}. Since $\Xi(\delta,\cdot)$ is continuous by \eqref{eqn.XiLip}, it follows that $\Xi(\delta,y)=\delta$ for all $y$. Finally $\Xi(x,0)=x$ and $\Xi(x,H)=x$, because for $0<c\le\delta$ we have $\Delta_-'(-c)=0$ and $\Delta_+'(c)=0$, as $c<H/2-\delta$, and likewise $\Delta_-'(H-c)=0$ and $\Delta_+'(H+c)=0$.

\noindent\emph{$\Upsilon$ is a homeomorphism with the stated boundary values.} For each $x$ the map $\chi_x^{-1}$ is a homeomorphism of $[0,H]$ onto itself by (P2) and (P3), and $(x,y)\mapsto\chi_x^{-1}(y)$ is Lipschitz on $Z$ by (P3) and (P5). So $\Upsilon_{\mr{v}}$ is a Lipschitz homeomorphism of $Z$ onto itself, with the Lipschitz inverse $(x,y')\mapsto(x,\chi_x(y'))$. For each $y'$ the map $x\mapsto\Xi(x,y')$ is an increasing homeomorphism of $[0,\delta]$ onto itself, and $\Xi$ is Lipschitz on $Z$, by \eqref{eqn.XiDeriv}, (P3), and \eqref{eqn.XiLip}; so $\Upsilon_{\mr{h}}$ is a Lipschitz bijection of the compact set $Z$ onto itself, and hence a homeomorphism. Therefore $\Upsilon$, as a composition of Lipschitz homeomorphisms, is a Lipschitz homeomorphism of $Z$ onto itself. On the side $x=0$ we have $\chi_0=\mr{id}$ and $\Xi(0,\cdot)=0$, so $\Upsilon$ is the identity there. On the side $x=\delta$, by (P1) and the boundary values of $\Xi$, $\Upsilon(\delta,y)=\big(\Xi(\delta,\vartheta(y)),\vartheta(y)\big)=(\delta,\vartheta(y))$. On the sides $y=0$ and $y=H$ the map $\chi_x^{-1}$ fixes $0$ and $H$ by (P2), and $\Xi(x,0)=\Xi(x,H)=x$.

\noindent\emph{The Jacobian and the constant.} The map $\Upsilon$ is Lipschitz, hence differentiable almost everywhere, and since $\Upsilon_{\mr{v}}$ is bi-Lipschitz the chain rule $D\Upsilon=(D\Upsilon_{\mr{h}}\circ\Upsilon_{\mr{v}})\,D\Upsilon_{\mr{v}}$ is valid almost everywhere, all partial derivatives below existing almost everywhere (the exceptional sets are null, and so are their preimages under $\Upsilon_{\mr{v}}$). Write $y'=\chi_x^{-1}(y)$, and at a point where the derivatives exist put $\chi'=\chi_x'(y')$, $\dot\chi=\partial_s\chi_s(y')|_{s=x}$, and $\Xi_y=\partial_{y'}\Xi(x,y')$. Since $\Upsilon_{\mr{v}}(x,y)=(x,y')$, and differentiating the identity $\chi_x(y')=y$ with respect to $x$ and with respect to $y$ gives $\dot\chi+\chi'\,\partial_xy'=0$ and $\chi'\,\partial_yy'=1$, we have that
\[
D\Upsilon_{\mr{v}}(x,y)=\begin{pmatrix}1 & 0\\[2pt] -\dot\chi/\chi' & 1/\chi'\end{pmatrix}.
\]
Since $\Upsilon_{\mr{h}}(x,y')=(\Xi(x,y'),y')$, and $\partial_x\Xi(x,y')=\chi'$ by \eqref{eqn.XiDeriv}, we also have that
\[
D\Upsilon_{\mr{h}}(x,y')=\begin{pmatrix}\chi' & \Xi_y\\[2pt] 0 & 1\end{pmatrix}.
\]
Therefore
\[
D\Upsilon(x,y)=D\Upsilon_{\mr{h}}(x,y')\,D\Upsilon_{\mr{v}}(x,y)=\begin{pmatrix}\chi'-\Xi_y\dot\chi/\chi' & \Xi_y/\chi'\\[2pt] -\dot\chi/\chi' & 1/\chi'\end{pmatrix},
\]
and $\det D\Upsilon=\det D\Upsilon_{\mr{h}}\cdot\det D\Upsilon_{\mr{v}}=\chi'\cdot(1/\chi')=1$.
By (P3), (P5), and \eqref{eqn.XiLip} we have $|\chi'-1|$, $|1/\chi'-1|$, $|\dot\chi|$, $|\Xi_y|=O(\varepsilon)$, so every entry of $D\Upsilon-I$ is $O(\varepsilon)$, and (F2), applied on the convex set $Z$, shows that $\Upsilon$ is $(1+O(\varepsilon))$-BL once $\varepsilon_\ast$ is small enough. Finally $|\Upsilon(x,y)-(x,y)|\le|\Xi(x,y')-x|+|y'-y|$; the first term is $\big|\int_0^x(\chi_s'(y')-1)\,\mr{d}s\big|=O(\varepsilon\delta)$ by \eqref{eqn.XiDeriv} and (P3), and the second is $O(\varepsilon\delta)$ by (P4).
\end{proof}

\noindent\emph{The weighted zone lemma.}
We will also need a version of Lemma \ref{lem.Zone} in which the map preserves, instead of area, the measure $u(y)\,\mr{d}x\,\mr{d}y$, for a function $u$ of $y$ alone. This will be provided by the following corollary.

\begin{corollary}\label{cor.WeightedZone}
Let $\delta>0$, $H\ge10\delta$, $Z=[0,\delta]\times[0,H]$, and $0<\varepsilon\le3\varepsilon_\ast$. Let $u:[0,H]\rar[1-\varepsilon,1+\varepsilon]$ be measurable, and let $\vartheta$ be an increasing homeomorphism of $[0,H]$ onto itself which is $(1+\varepsilon)$-BL and satisfies $|\vartheta(y)-y|\le\varepsilon\delta$ for all $y$. Then there is a $(1+O(\varepsilon))$-BL homeomorphism $\Upsilon_u$ of $Z$ onto itself which is the identity on the sides $x=0$, $y=0$, and $y=H$, satisfies $\Upsilon_u(\delta,y)=(\delta,\vartheta(y))$ for $0\le y\le H$ and $|\Upsilon_u(z)-z|=O(\varepsilon\delta)$ for all $z\in Z$, and whose Jacobian satisfies
\begin{equation}\label{eqn.WeightedJac}
|\det D\Upsilon_u(z)|=\frac{u(z_2)}{u\big((\Upsilon_u(z))_2\big)}\qquad\text{for almost every}\quad z\in Z .
\end{equation}
\end{corollary}
\begin{proof}
Let $P_u(y)=\int_0^yu$ and $\varphi_u(x,y)=(x,P_u(y))$. Then $P_u$ is an increasing $(1+2\varepsilon)$-BL homeomorphism of $[0,H]$ onto $[0,H_u]$, where $H_u=P_u(H)\ge(1-\varepsilon)H\ge5\delta$, with $P_u'=u$ almost everywhere, so $\varphi_u$ is a $(1+2\varepsilon)$-BL homeomorphism of $Z$ onto $Z_u=[0,\delta]\times[0,H_u]$ with $|\det D\varphi_u(z)|=u(z_2)$. Put $\vartheta_u=P_u\circ\vartheta\circ P_u^{-1}$, an increasing homeomorphism of $[0,H_u]$ onto itself. It is $(1+2\varepsilon)^2(1+\varepsilon)$-BL, hence $(1+10\varepsilon)$-BL, and for $Y=P_u(y)$ we have \[|\vartheta_u(Y)-Y|=|P_u(\vartheta(y))-P_u(y)|\le(1+2\varepsilon)|\vartheta(y)-y|\le10\varepsilon\delta.\] So Lemma \ref{lem.Zone}, applied on $Z_u$ with $10\varepsilon\le30\varepsilon_\ast$ in place of $\varepsilon$, gives a $(1+O(\varepsilon))$-BL homeomorphism $\Upsilon^\circ$ of $Z_u$ onto itself with $|\det D\Upsilon^\circ|=1$ almost everywhere, equal to the identity on $x=0$, $Y=0$, and $Y=H_u$, with $\Upsilon^\circ(\delta,Y)=(\delta,\vartheta_u(Y))$, and with $|\Upsilon^\circ(z)-z|=O(\varepsilon\delta)$.

Let $\Upsilon_u=\varphi_u^{-1}\circ\Upsilon^\circ\circ\varphi_u$. This is a homeomorphism of $Z$ onto itself, and it is $(1+O(\varepsilon))$-BL, being a composition of three such maps. Since $\varphi_u$ maps the sides $x=0$, $y=0$, and $y=H$ of $Z$ onto the corresponding sides of $Z_u$, the map $\Upsilon_u$ is the identity on them, and $\Upsilon_u(\delta,y)=\varphi_u^{-1}(\delta,\vartheta_u(P_u(y)))=(\delta,\vartheta(y))$. The displacement of $\Upsilon_u$ is at most $(1+2\varepsilon)$ times that of $\Upsilon^\circ$, because $\varphi_u^{-1}$ is $(1+2\varepsilon)$-Lipschitz. For the Jacobian, the chain rule for bi-Lipschitz maps gives, for almost every $z$,
\begin{align*}
|\det D\Upsilon_u(z)|&=|\det D\varphi_u^{-1}(\Upsilon^\circ(\varphi_u(z)))|\cdot|\det D\Upsilon^\circ(\varphi_u(z))|\cdot|\det D\varphi_u(z)|\\
&=\frac{1}{u\big((\Upsilon_u(z))_2\big)}\cdot1\cdot u(z_2),
\end{align*}
using $|\det D\varphi_u^{-1}(x,Y)|=1/u(P_u^{-1}(Y))$ and $P_u^{-1}\big((\Upsilon^\circ(\varphi_u(z)))_2\big)=(\Upsilon_u(z))_2$.
\end{proof}

\noindent\emph{The gluing lemma.}
We now glue the maps of a strip decomposition across their edges. The zone at an edge may be placed in either of the two strips meeting there, and the zone maps in each strip are those of Corollary \ref{cor.WeightedZone}, for a function which is prescribed on that strip.

\begin{lemma}\label{lem.StripGlue} 
Let $Q=[0,w]\times[0,H]$ be cut into the strips $Q_i=[s_i,s_{i+1}]\times[0,H]$, $0\le i\le N$, where $0=s_0<s_1<\cdots<s_{N+1}=w$, let $0<\varepsilon\le\varepsilon_\ast$, let $L\ge1$, and suppose that we are given the following data.
\begin{enumerate}[topsep=2pt,itemsep=4pt,parsep=4pt]
\item[(i)] For each $0\le i\le N$, an $L$-BL homeomorphism $\varphi_i$ of $Q_i$ onto itself which fixes the bottom and top sides of $Q_i$ pointwise and maps each of its vertical sides onto itself. We write $\varphi_i^{\mr{L}}$ and $\varphi_i^{\mr{R}}$ for the increasing homeomorphisms of $[0,H]$ by which $\varphi_i$ acts on the left and on the right side of $Q_i$, and we assume that these are $(1+\varepsilon)$-BL.
\item[(ii)] For each $0\le i<N$, a \textit{zone} $Z_i$ for the edge $x=s_{i+1}$: a rectangle of height $H$ and width $\delta_i\le H/10$ with one vertical side on that edge, contained in $Q_i$ or in $Q_{i+1}$ (so that $\delta_i$ is at most the width of the strip containing $Z_i$). The zones $Z_0,\ldots,Z_{N-1}$ are required to have pairwise disjoint interiors.
\item[(iii)] For each $0\le i\le N$, a measurable \textit{weight} $u^{(i)}:[0,H]\rar[1-\varepsilon,1+\varepsilon]$.
\end{enumerate}
Suppose that the mismatches $\mu_i=\sup_y|\varphi_i^{\mr{R}}(y)-\varphi_{i+1}^{\mr{L}}(y)|$ satisfy $\mu_i\le\varepsilon\delta_i$ for every $0\le i<N$. Then there is a homeomorphism $\tilde\varphi$ of $Q$ onto itself with the following properties.
\begin{enumerate}[topsep=2pt,itemsep=4pt,parsep=4pt]
\item[(a)] On each $Q_i$, $\tilde\varphi=\Upsilon_i\circ\varphi_i$, where $\Upsilon_i$ is a $(1+O(\varepsilon))$-BL homeomorphism of $Q_i$ onto itself which is the identity outside the zones contained in $Q_i$, moves no point by more than $O(\varepsilon\max_k\delta_k)$, and satisfies $|\det D\Upsilon_i(z)|=u^{(i)}(z_2)/u^{(i)}\big((\Upsilon_i(z))_2\big)$ almost everywhere.
\item[(b)] For almost every $z\in Q_i$,
\[
|\det D\tilde\varphi(z)|=|\det D\varphi_i(z)|\cdot\frac{u^{(i)}\big((\varphi_i(z))_2\big)}{u^{(i)}\big((\tilde\varphi(z))_2\big)} .
\]
In particular $|\det D\tilde\varphi|=|\det D\varphi_i|$ almost everywhere on $Q_i$ if $u^{(i)}\equiv1$; and if $\varphi^\circ$ is a BL map on $Q$ with $|\det D\varphi^\circ(z)|=u^{(i)}(z_2)$ for almost every $z$ in the zones contained in $Q_i$, then $|\det D(\varphi^\circ\circ\tilde\varphi)|=|\det D(\varphi^\circ\circ\varphi_i)|$ almost everywhere on $Q_i$.
\item[(c)] $\tilde\varphi$ fixes the bottom and top sides of $Q$ pointwise, acts on the left side of $Q$ as $\varphi_0^{\mr{L}}$ and on the right side as $\varphi_N^{\mr{R}}$, and is $(1+O(\varepsilon))L$-BL.
\end{enumerate}
\end{lemma}

The model case is that of Lemma \ref{lem.StripShear}: if the strips of the present lemma are taken to be the strips $Q_i'$ of that lemma, the strip maps to be its shears $\varphi_i$, and all weights to be $1$, then $L=1+O(\varepsilon)$, the side maps are the profiles $P_i$, the mismatches are those defined after Lemma \ref{lem.StripShear}, and the lemma produces a $(1+O(\varepsilon))$-BL homeomorphism $\tilde\varphi$ of $Q$ with $|\det D\tilde\varphi|=|\det D\varphi_i|$ almost everywhere on $Q_i'$, so that $\tilde\varphi\circ\varphi_\perp$ realizes the density $u$ of Lemma \ref{lem.StripShear} exactly. Thus, if the side maps are $(1+\varepsilon)$-BL, an edge whose mismatch is at most $\varepsilon\delta_i$ is repaired by a map which is $(1+O(\varepsilon))$-BL, and the zone may be placed on whichever side of the edge has room. The role of the weights is explained by (b): if the glued map is followed by a map $\varphi^\circ$ whose Jacobian on the zones is the weight, the zone maps are invisible to the Jacobian of the composite.

\begin{proof}
\noindent\emph{The zone maps.} Fix $0\le i<N$ and suppose first that $Z_i\subseteq Q_i$, so that $Z_i=[s_{i+1}-\delta_i,s_{i+1}]\times[0,H]$. Put $\vartheta_i=\varphi_{i+1}^{\mr{L}}\circ(\varphi_i^{\mr{R}})^{-1}$, an increasing homeomorphism of $[0,H]$ onto itself which is $(1+\varepsilon)^2$-BL, hence $(1+3\varepsilon)$-BL, and satisfies $|\vartheta_i(y)-y|=|\varphi_{i+1}^{\mr{L}}(y')-\varphi_i^{\mr{R}}(y')|\le\mu_i\le\varepsilon\delta_i$ with $y'=(\varphi_i^{\mr{R}})^{-1}(y)$. Since $H\ge10\delta_i$, Corollary \ref{cor.WeightedZone}, translated in $x$ and applied with the weight $u^{(i)}$ and with $3\varepsilon\le3\varepsilon_\ast$ in place of $\varepsilon$, gives a $(1+O(\varepsilon))$-BL homeomorphism of $Z_i$ onto itself with trace $\vartheta_i$ on the edge $x=s_{i+1}$, equal to the identity on the other three sides of $Z_i$, moving no point by more than $O(\varepsilon\delta_i)$, and with the Jacobian \eqref{eqn.WeightedJac} for the weight $u^{(i)}$; we extend it by the identity to the rest of $Q_i$. If instead $Z_i\subseteq Q_{i+1}$, we put $\vartheta_i=\varphi_i^{\mr{R}}\circ(\varphi_{i+1}^{\mr{L}})^{-1}$ and build the map on $Q_{i+1}$ in the same way, with the weight $u^{(i+1)}$ and with the trace prescribed on the left side of $Z_i$ (the reflected form of the corollary).

\noindent\emph{The strip maps.} For each $0\le i\le N$ let $\Upsilon_i$ be the composition of the zone maps supported in $Q_i$. There are none, one, or two of them, and in the last case their zones have disjoint interiors. By (F1), the map $\Upsilon_i$ is a $(1+O(\varepsilon))$-BL homeomorphism of $Q_i$ onto itself. It has the Jacobian stated in (a), by \eqref{eqn.WeightedJac} on the zones contained in $Q_i$ and because it is the identity on the open complement of the zones; it fixes the bottom and top sides of $Q_i$; and it acts on each vertical side of $Q_i$ as the trace prescribed there, or as the identity if the zone for that edge lies in the neighboring strip or the side is a side of $Q$ (a zone whose width equals that of the strip touches the far side of the strip, where the zone map is the identity). Now define
\[
\hat\varphi_i=\Upsilon_i\circ\varphi_i\qquad\text{on}\quad Q_i .
\]
The map $\hat\varphi_i$ is a $(1+O(\varepsilon))L$-BL homeomorphism of $Q_i$ onto itself which fixes the bottom and top sides of $Q_i$ pointwise, and by the chain rule, which holds almost everywhere because the preimage under $\varphi_i$ of a null set is null, $|\det D\hat\varphi_i(z)|=|\det D\Upsilon_i(\varphi_i(z))|\cdot|\det D\varphi_i(z)|$ almost everywhere, which is the identity in (b).

\noindent\emph{Agreement on the edges.} Consider the edge $x=s_{i+1}$, and suppose that $Z_i\subseteq Q_i$. On this edge $\hat\varphi_i$ acts as $\vartheta_i\circ\varphi_i^{\mr{R}}=\varphi_{i+1}^{\mr{L}}$, while $\hat\varphi_{i+1}$ acts as $\varphi_{i+1}^{\mr{L}}$, since $Q_{i+1}$ contains no zone for this edge and the zone at its other edge, if there is one, is the identity here. So the two maps agree on the edge, and the case $Z_i\subseteq Q_{i+1}$ is symmetric. The left side of $Q$ is a side of $Q_0$, which contains no zone there, so $\hat\varphi_0$ acts on it as $\varphi_0^{\mr{L}}$; likewise $\hat\varphi_N$ acts on the right side of $Q$ as $\varphi_N^{\mr{R}}$.

\noindent\emph{Assembly.} Let $\tilde\varphi$ be the map which equals $\hat\varphi_i$ on $Q_i$. By the agreement on the edges it is well defined and continuous, it is a bijection of $Q$ onto itself, and so it is a homeomorphism; (a) and the boundary behavior in (c) hold by construction, and the BL bound in (c) follows from (F1). The identity in (b) was proved above, and gives $|\det D\tilde\varphi|=|\det D\varphi_i|$ when $u^{(i)}\equiv1$. For the last assertion of (b), let $\varphi^\circ$ be as stated. Off a null set, $\tilde\varphi(z)$ lies in the interior of a zone $Z$ contained in $Q_i$ if and only if $\varphi_i(z)$ does. If it does, then $|\det D\varphi^\circ(\tilde\varphi(z))|=u^{(i)}\big((\tilde\varphi(z))_2\big)$ and $|\det D\varphi^\circ(\varphi_i(z))|=u^{(i)}\big((\varphi_i(z))_2\big)$ for almost every such $z$, and multiplying the identity in (b) by the first of these gives $|\det D(\varphi^\circ\circ\tilde\varphi)(z)|=|\det D(\varphi^\circ\circ\varphi_i)(z)|$. If it does not, then $\tilde\varphi(z)=\varphi_i(z)$ and $D\Upsilon_i(\varphi_i(z))=I$, so the two Jacobians agree trivially. The exceptional sets are null, with null preimages under $\tilde\varphi$ and $\varphi_i$.
\end{proof}

\noindent\emph{Boundary fixing.}
The final result in this section, which is another application of the zone lemma above, turns a map which preserves each side of a rectangle into one which fixes the boundary pointwise, at no cost to the Jacobian.

\begin{lemma}\label{lem.BoundaryFix} 
Let $Q=[0,w]\times[0,H]$, let $0<\varepsilon\le\varepsilon_\ast$, and let $\varphi$ be an $L$-BL homeomorphism of $Q$ onto itself which maps each side of $Q$ onto itself, acting on each side as an increasing homeomorphism which is $(1+\varepsilon)$-BL and moves no point by more than $\varepsilon\min\{w,H\}$. Then there is a $(1+O(\varepsilon))L$-BL homeomorphism $\hat\varphi$ of $Q$ onto itself which fixes $\partial Q$ pointwise, such that $|\det D\hat\varphi|=|\det D\varphi|$ almost everywhere and $|\hat\varphi(z)-\varphi(z)|=O(\varepsilon\min\{w,H\})$ for all $z\in Q$.
\end{lemma}
\begin{proof}
Put $\delta=\min\{w,H\}/5$, and write $\varphi^{\mr{L}},\varphi^{\mr{R}}$ (homeomorphisms of $[0,H]$) and $\varphi^{\mr{B}},\varphi^{\mr{T}}$ (homeomorphisms of $[0,w]$) for the maps by which $\varphi$ acts on the left, right, bottom, and top sides of $Q$. Each of these, and each of their inverses, is $(1+\varepsilon)$-BL and moves no point by more than $\varepsilon\min\{w,H\}=5\varepsilon\delta$. Let $\varrho^{\mr{L}}$ be the map of Lemma \ref{lem.Zone}, applied with $5\varepsilon\le30\varepsilon_\ast$ in place of $\varepsilon$, on the zone $[0,\delta]\times[0,H]$, with the trace $(\varphi^{\mr{L}})^{-1}$ prescribed on the side $x=0$ and the identity on the side $x=\delta$, extended by the identity to the rest of $Q$; the zone qualifies because $H\ge5\delta$. Let $\varrho^{\mr{R}}$ be defined likewise on $[w-\delta,w]\times[0,H]$, with the trace $(\varphi^{\mr{R}})^{-1}$ on the side $x=w$, and let $\varrho^{\mr{B}}$ and $\varrho^{\mr{T}}$ be defined on the horizontal collars $[0,w]\times[0,\delta]$ and $[0,w]\times[H-\delta,H]$, by the mirror image of the lemma, with the traces $(\varphi^{\mr{B}})^{-1}$ and $(\varphi^{\mr{T}})^{-1}$ on the bottom and top sides of $Q$; here $w\ge5\delta$. Each of the four maps is a homeomorphism of $Q$ onto itself, it is $(1+O(\varepsilon))$-BL by (F1), it has $|\det D\varrho^\bullet|=1$ almost everywhere, and it moves no point by more than $O(\varepsilon\delta)$. Moreover,  each of them acts on its own side of $Q$ as the inverse of the corresponding side map of $\varphi$ and is the identity on the other three sides. Hence $\varrho=\varrho^{\mr{T}}\circ\varrho^{\mr{B}}\circ\varrho^{\mr{R}}\circ\varrho^{\mr{L}}$ is a $(1+O(\varepsilon))$-BL homeomorphism of $Q$ onto itself with $|\det D\varrho|=1$ almost everywhere, which acts on $\partial Q$ as the inverse of $\varphi|_{\partial Q}$, and $\hat\varphi=\varrho\circ\varphi$ fixes $\partial Q$ pointwise. It is $(1+O(\varepsilon))L$-BL, and $|\det D\hat\varphi(z)|=|\det D\varrho(\varphi(z))|\cdot|\det D\varphi(z)|=|\det D\varphi(z)|$ almost everywhere. Finally $|\hat\varphi(z)-\varphi(z)|=|\varrho(\varphi(z))-\varphi(z)|=O(\varepsilon\delta)$.
\end{proof}

\section{The density on a cell}\label{sec.Density}

In this section we fix $j\ge j_0$ and a cell $Q$ of level $j+1$, and we describe the relative density $v_{j+1}$ on $Q$ in the detail which the construction of Section \ref{sec.CellLemma} requires.

\noindent\emph{The setting.}
We assume that $Q\subseteq\{x_1\ge0\}$, so that the merged column of $Q$ (i.e. the sub-cells having width $q_j+q_{j-1}$) is its first column. For a cell in $\{x_1<0\}$ the merged column is the last one; such cells are reduced to the present case by a reflection at the end of Section \ref{sec.CellLemma}. Write
\[
Q=[s_0,s_0+w_Q)\times[y_0,y_0+h_{j+1}),\qquad w_Q\in\{q_{j+1},q_{j+1}+q_j\},\quad y_0\in h_{j+1}\Z,
\]
and $T_j=h_{j+1}/h_j$, a power of $2$. The sub-cells of $Q$ are the rectangles
\[
Q_{i,t}=[s_i,s_i+w_i)\times[y_0+th_j,y_0+(t+1)h_j),\qquad 0\le i\le N_Q,\quad 0\le t<T_j,
\]
where $w_0=q_j+q_{j-1}$, and $w_i=q_j$ and $s_i=s_0+q_{j-1}+iq_j$ for $i\ge1$; here $N_Q=a_{j+1}-1$ if $w_Q=q_{j+1}$ and $N_Q=a_{j+1}$ otherwise. We call $i$ the \textit{column} and $t$ the \textit{row-block} of $Q_{i,t}$, and we also use this terminology for the vertical strips $[s_i,s_i+w_i)\times[y_0,y_0+h_{j+1})$ and the horizontal strips $[s_0,s_0+w_Q)\times[y_0+th_j,y_0+(t+1)h_j)$ of $Q$ which they index. The \textit{rows} of $Q$ are the horizontal lines $x_2=y_0+m$, $0\le m<h_{j+1}$, and we refer to this line as row $m$; it lies in the row-block $\lfloor m/h_j\rfloor$. If $a_{j+1}=1$ and $w_Q=q_{j+1}$ then $N_Q=0$ and $Q$ consists of its merged column alone; such cells are treated separately in Section \ref{sec.CellLemma}, and the statements below which concern columns $i\ge1$ are empty for them.

For a column $i$ and a row $m$ put $\gamma_i(m)=s_i\alpha_1+(y_0+m)\alpha_2$. By the definition of $R_M$, the number of points of $\Lambda$ on row $m$ of column $i\ge1$ is $q_j|W|+R_{q_j}(\gamma_i(m))$, and the number on row $m$ of column $0$ is $w_0|W|+R_{q_j}(\gamma_0(m))+R_{q_{j-1}}(\gamma_0(m)+q_j\alpha_1)$. We write
\begin{align*}
R_i[t]&=\sum_{th_j\le m<(t+1)h_j}R_{q_j}(\gamma_i(m))\qquad(0\le i\le N_Q),\\
R_0'[t]&=\sum_{th_j\le m<(t+1)h_j}R_{q_{j-1}}(\gamma_0(m)+q_j\alpha_1),
\end{align*}
for the \textit{block counts} of $Q$ (the square brackets distinguish them from the quantities $R_M(\gamma)$), so that
\[
\#(\Lambda\cap Q_{i,t})=w_ih_j|W|+R_i[t]\quad(i\ge1),\qquad\#(\Lambda\cap Q_{0,t})=w_0h_j|W|+R_0[t]+R_0'[t] .
\]
By Lemma \ref{lem.Block}, $|R_i[t]|<3h_j$ for every $i$ and $|R_0'[t]|<3h_j$. By the definition of $v_{j+1}$ in Section \ref{sec.Hierarchy},
\begin{equation}\label{eqn.vFormula}
v_{j+1}\big|_{Q_{i,t}}=\frac{1}{u_{j+1}(Q)}\Big(|W|+\frac{R_i[t]+\mathbf 1[i=0]\,R_0'[t]}{w_ih_j}\Big),\qquad u_{j+1}(Q)=|W|\big(1+O(\varepsilon_j)\big),
\end{equation}
the second statement being \eqref{eqn.DensityBound} at level $j+1$, together with $\varepsilon_{j+1}\le\varepsilon_j$.

With $a$ the left endpoint of the window $W=(a,b)$, for $0\le i\le N_Q$ let
\[
x_i=q_j\big(a-s_i\alpha_1-y_0\alpha_2\big)\in\R/\Z,\qquad\beta_j=q_j\alpha_2,\qquad\eta_j=q_jD_j .
\]
Then $q_j(a-\gamma_i(m))\equiv x_i-m\beta_j$ modulo $1$, and we call $x_i-m\beta_j\in\R/\Z$ the \textit{phase} of row $m$ in column $i$. Since $q_j\alpha_1=p_j+D_j$ and $q_jq_{j-1}\alpha_1\equiv q_{j-1}D_j$ modulo $1$,
\begin{equation}\label{eqn.PhaseStep}
x_{i+1}=x_i-\eta_j\quad(1\le i<N_Q),\qquad x_1=x_0-(q_j+q_{j-1})D_j .
\end{equation}
By \eqref{eqn.CFFacts} we have $|\eta_j|=q_j|D_j|<1/a_{j+1}$ and, since $(N_Q-1)q_j<q_{j+1}$, also $(N_Q-1)|\eta_j|<q_{j+1}|D_j|<1$: the phases $x_1,\ldots,x_{N_Q}$ form an arithmetic progression in $\R/\Z$ of step $-\eta_j$ which sweeps out less than one full turn of the circle. The vertical step $\beta_j$, on the other hand, is an arbitrary real number. The wave $\omega_j$ of Section \ref{sec.Reduction} is constant except at the two points $0$ and $-\tau_j$ of $\R/\Z$, where it jumps by $-1$ and by $+1$ respectively, and we call these its \textit{jump points}; a point $x\in\R/\Z$ is \textit{within $\delta$ of a jump point} if its distance to $\{0,-\tau_j\}+\Z$ is less than $\delta$. (If $\tau_j=0$ then $\omega_j\equiv0$ and the two jump points coincide; the statements below remain valid, with this point counted twice.)

\noindent\emph{Two ways of realizing the density.} 
Lemma \ref{lem.StripShear} suggests two ways of realizing $v_{j+1}$ on $Q$, by a vertical shear in each column, or by a horizontal shear in each row-block, with the mismatches at the common edges removed in either case by Lemma \ref{lem.StripGlue}. Neither of them works by itself:
\begin{enumerate}[topsep=2pt,itemsep=4pt,parsep=4pt]
	\item[(i)] For vertical shears, if $\beta_j$ is very close to a rational number whose denominator $q$ is much smaller than $a_{j+1}$, then the block counts of two adjacent columns can differ by about $h_j/q$ on every row-block in a large part of the cell, and the mismatch at the edge between them is too large to be removed by a map which is $(1+O(\varepsilon_j))$-BL.
	\item[(ii)] For horizontal shears, if $\beta_j$ is close to a rational number whose denominator does not divide $h_j$, then the block counts of consecutive row-blocks can differ by an amount of order $1$ in a positive proportion of the columns, with no decay in $a_{j+1}$, and the mismatches between consecutive row-blocks can again be too large.
\end{enumerate}
The decomposition of Lemma \ref{lem.Split} below separates the two effects: it isolates the part of each block count which is responsible for the second effect, and this part is realized by vertical shears, while the rest, which changes little from one row-block to the next, is realized by horizontal shears.

\noindent\emph{The sub-cell counts.}
Using the notation of the sub-cells, Lemma \ref{lem.Block} gives us the following result.

\begin{lemma}\label{lem.SubCell}
For $0\le i\le N_Q$ and $0\le t<T_j$,
\[
R_i[t]=\sum_{th_j\le m<(t+1)h_j}\Big(\omega_j\big(x_i-m\beta_j\big)+\epsilon_j\big(\gamma_i(m)\big)\Big),
\]
where $\epsilon_j$ is the integer-valued error of Lemma \ref{lem.Block}; so $|\epsilon_j(\gamma_i(m))|\le2$, and $\epsilon_j(\gamma_i(m))=0$ unless the phase $x_i-m\beta_j$ is within $|\eta_j|$ of a jump point.
\end{lemma}
\begin{proof}
Apply Lemma \ref{lem.Block} with $k=j$ and $\gamma=\gamma_i(m)$, and sum over the rows of the row-block. The main term is $\omega_j(q_j(a-\gamma_i(m)))=\omega_j(x_i-m\beta_j)$, because $\omega_j$ is $1$-periodic. The error vanishes unless $\|q_j(a-\gamma_i(m))\|<q_j|D_j|$ or $\|q_j(b-\gamma_i(m))\|<q_j|D_j|$; the first condition says that $x_i-m\beta_j$ is within $|\eta_j|$ of $0$, and since $q_j(b-\gamma)=q_j(a-\gamma)+q_j|W|\equiv x_i-m\beta_j+\tau_j$ modulo $1$, the second says that it is within $|\eta_j|$ of $-\tau_j$.
\end{proof}

\noindent\emph{Adjacent sub-cells.} 
The next lemma says that adjacent columns have the same row counts except on the rows whose phase is near a jump point, and it records how many columns can have their phase in a given set.

\begin{lemma}\label{lem.Adjacent} 
Let $0\le m<h_{j+1}$.
\begin{enumerate}[topsep=2pt,itemsep=4pt,parsep=4pt]
\item[(i)] For $1\le i<N_Q$ the difference $R_{q_j}(\gamma_i(m))-R_{q_j}(\gamma_{i+1}(m))$ has absolute value at most $2$, and it is $0$ unless the phase $x_i-m\beta_j$ is within $2|\eta_j|$ of a jump point.
\item[(ii)] If $A\subseteq\R/\Z$ is a union of $n$ arcs, then for every $y\in\R/\Z$,
\[
\#\{1\le i\le N_Q : x_i-y\in A\}\le\frac{|A|}{|\eta_j|}+2n .
\]
\end{enumerate}
\end{lemma}
\begin{proof}
(i) We have $\gamma_{i+1}(m)-\gamma_i(m)=q_j\alpha_1\equiv D_j$ modulo $1$. Write $\gamma=\gamma_i(m)$ and $\gamma'=\gamma+D_j$. The count $\#\{0\le n<q_j:n\alpha_1+\gamma\in(a,b)+\Z\}$ can change when $\gamma$ is replaced by $\gamma'$ only if $n\alpha_1+\gamma$ and $n\alpha_1+\gamma'$ lie on different sides of $a$ or of $b$ modulo $1$. This happens only if $a-n\alpha_1$ or $b-n\alpha_1$ lies in the closed interval with endpoints $\gamma$ and $\gamma'$, modulo $1$. As in the proof of Lemma \ref{lem.Block}, the points $a-n\alpha_1$, $0\le n<q_j$, are the points of the grid $a+\frac{1}{q_j}\Z$, each displaced by less than $|D_j|$, so any two of them are at distance greater than $1/q_j-|D_j|>|D_j|$ (by \eqref{eqn.CFFacts}, $|D_j|<1/q_{j+1}\le1/(2q_j)$, because a column $i+1\le N_Q$ with $i\ge1$ exists only if $a_{j+1}\ge2$). Hence at most one of them lies in an interval of length $|D_j|$, and the same holds for the points $b-n\alpha_1$. This establishes the upper bound of $2$ in part (i). Moreover, if $a-n\alpha_1$ lies in the interval, then multiplying by $q_j$ and using $q_jn\alpha_1=np_j+nD_j$ we find that $q_j(a-n\alpha_1-\gamma)\equiv x_i-m\beta_j-nD_j$ lies in the closed interval with endpoints $0$ and $\eta_j$, modulo $1$; since $|nD_j|<|\eta_j|$, the phase $x_i-m\beta_j$ is then within $2|\eta_j|$ of $0$. In the same way, since $q_j(b-\gamma)\equiv x_i-m\beta_j+\tau_j$, if $b-n\alpha_1$ lies in the interval then $x_i-m\beta_j$ is within $2|\eta_j|$ of $-\tau_j$.

(ii) By \eqref{eqn.PhaseStep} the points $x_i-y$, $1\le i\le N_Q$, lift to an arithmetic progression in $\R$ of step $|\eta_j|$ in absolute value, which is contained in an interval $J$ of length less than $1$. The preimage in $\R$ of an arc $B\subseteq\R/\Z$ meets $J$ in at most two intervals of total length at most $|B|$. An arithmetic progression of step $\sigma$ has at most $\ell/\sigma+1$ terms in an interval of length $\ell$, so the progression has at most $|B|/|\eta_j|+2$ terms in the preimage of $B$. Summing over the $n$ arcs of $A$ gives the bound.
\end{proof}

\noindent\emph{A rational approximation and the drift.}
We will use a decomposition in Lemma \ref{lem.Split} below which rests on a rational approximation to $\beta_j$. By Dirichlet's theorem, there are coprime integers $p$ and $q$ with $1\le q\le q_{j+1}$ and $|q\beta_j-p|<1/q_{j+1}$. For the rest of this section and in Section \ref{sec.CellLemma} we fix such a pair, so that
\begin{equation}\label{eqn.Dirichlet}
1\le q\le q_{j+1}\qquad\text{and}\qquad|\theta|<\frac{1}{qq_{j+1}},\qquad\text{where}\quad\theta=\beta_j-\frac pq .
\end{equation}
We call $\theta$ the \textit{drift}. We stress that nothing is assumed about the size of $q$: it may be as small as $1$ or as large as $q_{j+1}$, and $\beta_j$ itself may be rational. What \eqref{eqn.Dirichlet} gives, through \eqref{eqn.CFFacts}, $h_{j+1}<2q_{j+1}$, and $h_j<2q_j$, are the three inequalities
\begin{equation}\label{eqn.DriftBounds}
h_{j+1}|\theta|<\frac2q,\qquad h_j|\theta|<\frac{4|\eta_j|}{q},\qquad q|\eta_j|<q_j .
\end{equation}
In the language of cells, the phase of row $m$ in column $i$ is the point $x_i-mp/q$ of the grid $x_i+\frac1q\Z$, displaced by $-m\theta$. The phases themselves are not close to one another, since $\beta_j$ is arbitrary, but \eqref{eqn.DriftBounds} says that their displacements from the grid change by less than $2/q$, that is, by less than two spacings of the grid, over the whole height of the cell, and by less than $4/q$ times the horizontal phase step $|\eta_j|$ over one row-block.

For $0\le t<T_j$ and $0\le i\le N_Q$ we define the \textit{drifted phases}
\[
x_i^{(t)}=x_i-th_j\theta\in\R/\Z ,
\]
so that $x_{i+1}^{(t)}=x_i^{(t)}-\eta_j$ for $1\le i<N_Q$ and $x_i^{(t+1)}=x_i^{(t)}-h_j\theta$. For a row $m$ in the row-block $t$, the phase of row $m$ in column $i$ is
\begin{equation}\label{eqn.PhaseDecomp}
x_i-m\beta_j=x_i^{(t)}-\frac{mp}{q}-(m-th_j)\,\theta ,
\end{equation}
and since $0\le m-th_j<h_j$ the last term has absolute value less than $h_j|\theta|$. So within one row-block the phases of a column lie within $h_j|\theta|$ of the $q$ points $x_i^{(t)}-rp/q$, $0\le r<q$, which the rows of the block visit in the order of their residues modulo $q$.

\noindent\emph{Splitting the block counts.} 
For a residue $0\le r<q$ and a row-block $t$ put
\[
\iota_{r,t}=\#\{th_j\le m<(t+1)h_j : m\equiv r~\mathrm{mod}~q\}-\frac{h_j}{q},
\]
and let
\[
\bar\omega(x)=\frac1q\sum_{r=0}^{q-1}\omega_j\Big(x-\frac rq\Big)\qquad(x\in\R/\Z)
\]
be the average of $\omega_j$ over the grid $x+\frac1q\Z$. For $0\le i\le N_Q$ and $0\le t<T_j$ we define the \textit{averaged}, \textit{remainder}, \textit{crossing}, and \textit{correction} parts of the block count $R_i[t]$ by
\begin{align*}
R^{\mr{av}}_i[t]&=h_j\,\bar\omega\big(x_i^{(t)}\big),\\
R^{\mr{rm}}_i[t]&=\sum_{r=0}^{q-1}\iota_{r,t}\,\omega_j\Big(x_i^{(t)}-\frac{rp}{q}\Big),\\
R^{\mr{cr}}_i[t]&=\sum_{th_j\le m<(t+1)h_j}\Big[\omega_j\big(x_i-m\beta_j\big)-\omega_j\Big(x_i^{(t)}-\frac{mp}{q}\Big)\Big],\\
R^{\epsilon}_i[t]&=\sum_{th_j\le m<(t+1)h_j}\epsilon_j\big(\gamma_i(m)\big).
\end{align*}
The averaged part is what the block count would be if the rows of the block visited each residue class exactly $h_j/q$ times and did not drift within the block; the remainder part accounts for the unequal numbers of visits; the crossing part accounts for the rows whose phase drifts across a jump point within the block; and the correction part is the error of Lemma \ref{lem.Block}.

\begin{lemma}\label{lem.Split}
	Let $0\le i\le N_Q$ and $0\le t<T_j$.
	\begin{enumerate}[topsep=2pt,itemsep=4pt,parsep=4pt]
		\item[(i)] $R_i[t]=R^{\mr{av}}_i[t]+R^{\mr{rm}}_i[t]+R^{\mr{cr}}_i[t]+R^{\epsilon}_i[t]$.
		\item[(ii)] $|\bar\omega|\le1/q$, and $\bar\omega$ is constant except at the points of $\frac1q\Z$, where it jumps by $-1/q$, and at the points of $-\tau_j+\frac1q\Z$, where it jumps by $+1/q$ (the two jumps canceling where these sets meet).
		\item[(iii)] $\sum_{r=0}^{q-1}\iota_{r,t}=0$, and for every residue $r$ and every set $I$ of consecutive row-blocks,
		\[
		\Big|\sum_{t\in I}\iota_{r,t}\Big|<1 .
		\]
		In particular $|\iota_{r,t}|<1$.
		\item[(iv)] $|R^{\mr{rm}}_i[t]|\le h_j$.
	\end{enumerate}
\end{lemma}
\begin{proof}
	(i) By Lemma \ref{lem.SubCell} and the definition of the crossing part, $R_i[t]-R^{\epsilon}_i[t]-R^{\mr{cr}}_i[t]$ is the sum of $\omega_j(x_i^{(t)}-mp/q)$ over the rows $m$ of the block. Each term depends only on the residue $r$ of $m$ modulo $q$, and the residue $r$ occurs $h_j/q+\iota_{r,t}$ times, so the sum is
	\[
	\frac{h_j}{q}\sum_{r=0}^{q-1}\omega_j\Big(x_i^{(t)}-\frac{rp}{q}\Big)+R^{\mr{rm}}_i[t]=h_j\,\bar\omega\big(x_i^{(t)}\big)+R^{\mr{rm}}_i[t],
	\]
	because $r\mapsto rp$ permutes $\Z/q\Z$.
	
	(ii) By the definition of $\omega_j$, $q\bar\omega(x)$ is the number of points of the grid $x+\frac1q\Z$ which lie in the arc $[1-\tau_j,1)$ of $\R/\Z$, minus $q\tau_j$. An arc of length $\tau_j$ contains $\lfloor q\tau_j\rfloor$ or $\lfloor q\tau_j\rfloor+1$ points of a grid of spacing $1/q$, so $|q\bar\omega|\le1$, and the jumps of $\bar\omega$ are those of its $q$ terms divided by $q$.
	
	(iii) The sum over $r$ is $h_j-q(h_j/q)=0$. If $I=\{t_1,\ldots,t_2-1\}$ is a set of consecutive row-blocks, then the rows of these row-blocks are the integers $t_1h_j\le m<t_2h_j$, so
	\[
	\sum_{t'\in I}\iota_{r,t'}=\#\{t_1h_j\le m<t_2h_j : m\equiv r~\mathrm{mod}~q\}-\frac{(t_2-t_1)h_j}{q},
	\]
	and the right hand side has absolute value less than $1$.
	
	(iv) Let $S$ be the set of residues $r$ with $\omega_j(x_i^{(t)}-rp/q)=1-\tau_j$. Since $\omega_j$ takes only the values $1-\tau_j$ and $-\tau_j$, and $\sum_r\iota_{r,t}=0$, we have that $R^{\mr{rm}}_i[t]=\sum_{r\in S}\iota_{r,t}$, which is the number of rows of the block whose residue lies in $S$, minus $\#S\cdot h_j/q$. Both of these numbers lie between $0$ and $h_j$.
\end{proof}

\noindent\emph{Three estimates.}
We now come to the estimates which are needed for the vertical and the horizontal shears. For $1\le i<N_Q$ and $0\le t<T_j$ let $A_{i,t}$ be the half-open arc of $\R/\Z$ of length $|\eta_j|$ whose endpoints are the drifted phases $x_i^{(t)}$ and $x_{i+1}^{(t)}=x_i^{(t)}-\eta_j$, namely $(x_{i+1}^{(t)},x_i^{(t)}]$ if $\eta_j>0$ and $(x_i^{(t)},x_{i+1}^{(t)}]$ if $\eta_j<0$. As $t$ increases, $A_{i,t}$ moves rigidly along the circle, by $-h_j\theta$ per step.

For a residue $0\le r<q$ and a jump point $\mathfrak j\in\{0,-\tau_j\}$, consider the row-blocks $t$ for which the point $rp/q+\mathfrak j$ of $\R/\Z$ lies in $A_{i,t}$. These row-blocks are divided into \textit{runs} of the pair $(r,\mathfrak j)$: a run is a set $\{t_1,t_1+1,\ldots,t_2-1\}$ of consecutive row-blocks, all with this property, such that neither $t_1-1$ nor $t_2$ has it. We let $N_i$ be the total number of runs, summed over the $2q$ pairs $(r,\mathfrak j)$.

\begin{lemma}\label{lem.Estimates} 
Using the notation of \eqref{eqn.Dirichlet}, we have the following estimates.
\begin{enumerate}[topsep=2pt,itemsep=4pt,parsep=4pt]
\item[(E1)] For $1\le i<N_Q$ and every set $I$ of consecutive row-blocks,
\[
\Big|\sum_{t\in I}\big(R^{\mr{rm}}_i[t]-R^{\mr{rm}}_{i+1}[t]\big)\Big|\le N_i,\qquad\text{and}\qquad N_i<2q|\eta_j|+6<2q_j+6 .
\]
\item[(E2)] For every $0\le t<T_j$,
\[
\sum_{i=0}^{N_Q}\big|R^{\mr{cr}}_i[t]\big|=O(h_j)\qquad\text{and}\qquad\sum_{i=0}^{N_Q}\big|R^{\epsilon}_i[t]\big|=O(h_j),
\]
and for every $0\le t<T_j-1$,
\[
\sum_{i=0}^{N_Q}\big|R^{\mr{av}}_i[t]-R^{\mr{av}}_i[t+1]\big|=O(h_j).
\]
\item[(E3)] For $0\le i\le N_Q$ and $0\le t_0\le T_j$,
\[
\Big|\sum_{t<t_0}R^{\mr{rm}}_i[t]\Big|<q+7 .
\]
\end{enumerate}
\end{lemma}

\begin{proof}
\noindent\emph{(E1).} Fix $1\le i<N_Q$. We first express the difference of the remainder parts of two adjacent columns through the jump points. Since $x_{i+1}^{(t)}=x_i^{(t)}-\eta_j$ and $|\eta_j|<1$, the translate $A_{i,t}-rp/q$ is the half-open arc of length $|\eta_j|$ between the points $x_i^{(t)}-rp/q$ and $x_{i+1}^{(t)}-rp/q$ of $\R/\Z$, closed at the first of these if $\eta_j>0$ and at the second if $\eta_j<0$. A jump point $\mathfrak j$ lies in this arc if and only if $rp/q+\mathfrak j\in A_{i,t}$. Using this, we have that
\[
\omega_j\Big(x_i^{(t)}-\frac{rp}{q}\Big)-\omega_j\Big(x_{i+1}^{(t)}-\frac{rp}{q}\Big)=\sum_{\mathfrak j\in\{0,-\tau_j\}:\ rp/q+\mathfrak j\in A_{i,t}}\pm1 ,
\]
where the sign is the jump of $\omega_j$ at $\mathfrak j$ if $\eta_j>0$ and its negative if $\eta_j<0$. Multiplying by $\iota_{r,t}$ and summing over $r$ gives
\begin{equation}\label{eqn.RmDiff}
R^{\mr{rm}}_i[t]-R^{\mr{rm}}_{i+1}[t]=\sum_{(r,\mathfrak j):\ rp/q+\mathfrak j\in A_{i,t}}\pm\,\iota_{r,t},
\end{equation}
with signs which depend only on $\mathfrak j$. Now sum over $t\in I$ and exchange the two finite sums: the result is $\sum_{(r,\mathfrak j)}\pm\sum_t\iota_{r,t}$, the inner sum running over the $t\in I$ with $rp/q+\mathfrak j\in A_{i,t}$. For a fixed pair $(r,\mathfrak j)$ this set of $t$ is the intersection of $I$ with the runs of the pair. Hence by Lemma \ref{lem.Split}(iii), the whole sum has absolute value at most $N_i$.

For the bound on $N_i$ we lift to $\R$. Let $\tilde A_0\subseteq\R$ be an interval of length $|\eta_j|$ which is a lift of the arc $A_{i,0}$, and put $\tilde A_t=\tilde A_0-th_j\theta$, so that $\tilde A_t$ is a lift of $A_{i,t}$. As $t$ increases these intervals move monotonically, through a total distance $(T_j-1)h_j|\theta|$, so all of them are contained in an interval $J$ of length at most $|\eta_j|+h_{j+1}|\theta|$.

Now fix a pair $(r,\mathfrak j)$, and let $\tilde y\in\R$ be a lift of the point $rp/q+\mathfrak j$. Since the intervals $\tilde A_t$ move monotonically, the set of $t$ with $\tilde y\in\tilde A_t$ is a set of consecutive row-blocks, so if it is nonempty then it is contained in a single run of the pair. Every run arises in this way from at least one lift $\tilde y\in J$, since $rp/q+\mathfrak j\in A_{i,t}$ if and only if $\tilde A_t$ contains a lift of this point. Therefore the number of runs of the pair is at most the number of its lifts which lie in $J$.

Finally, for each $\mathfrak j$, as $r$ runs through the residues the lifts of the points $rp/q+\mathfrak j$ make up the set $\mathfrak j+\frac1q\Z$, each point of which is a lift for exactly one $r$, and this set has at most $q\ell+1$ points in an interval of length $\ell$. Adding the contributions of the two jump points, we have that
\[
N_i\le2\big(q(|\eta_j|+h_{j+1}|\theta|)+1\big)<2q|\eta_j|+6<2q_j+6 ,
\]
by the first and third inequalities of \eqref{eqn.DriftBounds}.

\noindent\emph{(E2).} All three bounds come from Lemma \ref{lem.Adjacent}(ii), which we apply to the drifted progression $x_1^{(t)},\ldots,x_{N_Q}^{(t)}$ (a translate of $x_1,\ldots,x_{N_Q}$, with the same steps); the terms with $i=0$ are bounded trivially.

\emph{The crossing part.} Let $m$ be a row of the row-block $t$, and consider the term
\[
\omega_j\big(x_i-m\beta_j\big)-\omega_j\Big(x_i^{(t)}-\frac{mp}{q}\Big)
\]
of the sum which defines $R^{\mr{cr}}_i[t]$. Since $\omega_j$ takes only two values, which differ by $1$, this term has absolute value at most $1$. By \eqref{eqn.PhaseDecomp} the first argument is obtained from the second by subtracting $(m-th_j)\theta$, and $0\le m-th_j<h_j$. Since $\omega_j$ is constant away from its jump points, the term is nonzero only if some jump point $\mathfrak j$ has the form $x_i^{(t)}-mp/q-s\theta$ with $0\le s\le h_j$, that is, only if
\[
x_i^{(t)}-\Big(\frac{mp}{q}+\mathfrak j\Big)\in\{s\theta:0\le s\le h_j\} ,
\]
which is a closed arc of $\R/\Z$ of length at most $h_j|\theta|$. For given $m$ and $\mathfrak j$, Lemma \ref{lem.Adjacent}(ii) with $y=mp/q+\mathfrak j$ shows that this happens for at most $h_j|\theta|/|\eta_j|+2<4/q+2\le6$ columns $i\ge1$, by \eqref{eqn.DriftBounds}. Two jump points and $h_j$ rows give $\sum_{i\ge1}|R^{\mr{cr}}_i[t]|\le12h_j$, and $|R^{\mr{cr}}_0[t]|\le h_j$.

\emph{The correction part.} By Lemma \ref{lem.SubCell}, $\epsilon_j(\gamma_i(m))\not=0$ only if $x_i-m\beta_j$ lies in the union of the two open arcs of length $2|\eta_j|$ centered at the jump points, and by Lemma \ref{lem.Adjacent}(ii) with $y=m\beta_j$ this happens for at most $4+4=8$ columns $i\ge1$. Since $|\epsilon_j|\le2$, summing over the $h_j$ rows gives $\sum_{i\ge1}|R^\epsilon_i[t]|\le16h_j$, and $|R^\epsilon_0[t]|\le2h_j$.

\emph{The averaged part.} By Lemma \ref{lem.Split}(ii), $|\bar\omega(x)-\bar\omega(x')|$ is at most $1/q$ times the number of points of $\frac1q\Z\cup(-\tau_j+\frac1q\Z)$, counted with multiplicity, in the arc between $x$ and $x'$. Here $x=x_i^{(t)}$ and $x'=x_i^{(t+1)}=x-h_j\theta$, so the arc has length $h_j|\theta|$. Lift the progression $x_1^{(t)},\ldots,x_{N_Q}^{(t)}$ to $\R$, where it occupies an interval of length less than $1$, and lift the arcs with it. A fixed lift of a point of $\frac1q\Z\cup(-\tau_j+\frac1q\Z)$ lies in the lifted arcs of at most $h_j|\theta|/|\eta_j|+1<5$ columns, since the lifted arcs are translates of one interval of length $h_j|\theta|$ by the terms of an arithmetic progression of step $|\eta_j|$, and only the lifts in an interval of length less than $1+h_j|\theta|$ lie in any of them; there are at most $2(q(1+h_j|\theta|)+1)=O(q)$ such lifts, by the first inequality of \eqref{eqn.DriftBounds}. Hence \[\sum_{i\ge1}|R^{\mr{av}}_i[t]-R^{\mr{av}}_i[t+1]|\le(h_j/q)\cdot5\cdot O(q)=O(h_j),\] and the term with $i=0$ is at most $(h_j/q)\cdot2(qh_j|\theta|+1)\le6h_j/q$, since $qh_j|\theta|\le qh_{j+1}|\theta|<2$.

\noindent\emph{(E3).} Fix $i$ and $t_0$, and put $M=t_0h_j$. Consider first the sum obtained by replacing the drifted phases $x_i^{(t)}$ by $x_i$ in the definition of the remainder parts. Since $\sum_{t<t_0}\iota_{r,t}=\#\{0\le m<M:m\equiv r~\mathrm{mod}~q\}-M/q$, since $\omega_j(x_i-mp/q)$ depends only on the residue of $m$, and since $\sum_r\omega_j(x_i-rp/q)=q\bar\omega(x_i)$ as in the proof of Lemma \ref{lem.Split}(i),
\[
\sum_{t<t_0}\sum_{r=0}^{q-1}\iota_{r,t}\,\omega_j\Big(x_i-\frac{rp}{q}\Big)=\sum_{m<M}\omega_j\Big(x_i-\frac{mp}{q}\Big)-M\bar\omega(x_i) .
\]
Write $M=nq+M'$ with $0\le M'<q$. The $n$ full periods contribute $nq\bar\omega(x_i)$ to the first term on the right, which cancels against the same amount in the second, leaving $\sum_{m<M'}\omega_j(x_i-mp/q)-M'\bar\omega(x_i)$, of absolute value less than $M'+M'/q<q+1$ because $|\omega_j|<1$ and $|\bar\omega|\le1/q$.

Now we put the drift back. For a fixed residue $r$, the value $\omega_j(x_i^{(t)}-rp/q)$ changes with $t$ exactly when the moving point $x_i^{(t)}$ crosses one of the points $rp/q+\mathfrak j$. Lifted to $\R$, the point $x_i^{(t)}$ moves monotonically through a total distance of at most $h_{j+1}|\theta|$, so it crosses each lift of such a point at most once, and a crossing which takes effect at the row-block $t^\ast$ changes $\omega_j(x_i^{(t)}-rp/q)$ by $\pm1$ for all $t\ge t^\ast$, hence changes the sum by $\pm\sum_{t^\ast\le t<t_0}\iota_{r,t}$, which has absolute value less than $1$ by Lemma \ref{lem.Split}(iii). The lifts which are crossed, over all $r$ and both $\mathfrak j$, are points of $\frac1q\Z\cup(-\tau_j+\frac1q\Z)$ in an interval of length at most $h_{j+1}|\theta|<2/q$, so there are fewer than $2(2+1)=6$ of them. Therefore the drifted sum differs from the undrifted one by less than $6$, and its absolute value is less than $q+7$.
\end{proof}

\noindent\emph{Adjacent full-cycle counts.}
For $0\le i\le N_Q$ and $0\le t<T_j$ we call $R_i[t]-R^{\mr{rm}}_i[t]$ the \textit{full-cycle count}, and we let $G_i$ be the step function on $[0,h_{j+1})$ which is equal to $(R_i[t]-R^{\mr{rm}}_i[t])/h_j$ on $[th_j,(t+1)h_j)$. In Section \ref{sec.CellLemma} the remainder parts are removed by vertical shears of the columns, and the density determined by the $G_i$ is then realized by horizontal shears. Since the vertical shears do not preserve the row-blocks, this requires a comparison of $G_i$ and $G_{i+1}$ over intervals which are not unions of row-blocks.

\begin{lemma}\label{lem.FullCycleDiff} 
Assume \eqref{eqn.Dirichlet}, and let $1\le i<N_Q$.
\begin{enumerate}[topsep=2pt,itemsep=4pt,parsep=4pt]
\item[(i)] For integers $0\le m_1<m_1+M\le h_{j+1}$, the number of rows $m_1\le m<m_1+M$ whose phase $x_i-m\beta_j$ is within $2|\eta_j|$ of a jump point is at most $2(M/q+1)(3+4q|\eta_j|)$.
\item[(ii)] For $0\le y\le z\le h_{j+1}$,
\[
\Big|\int_y^z\big(G_i-G_{i+1}\big)\Big|=O\Big(\Big(\frac{z-y}{q}+\min\Big\{1,\frac{z-y}{h_j}\Big\}\Big)\big(1+q|\eta_j|\big)\Big).
\]
\end{enumerate}
\end{lemma}
\begin{proof}
(i) For $m_1\le m<m_1+M$ we have, as in \eqref{eqn.PhaseDecomp} but with $m_1$ as the reference row, $x_i-m\beta_j=x_i-m_1\theta-(m-m_1)\theta-mp/q$. If the phase of row $m$ is within $2|\eta_j|$ of the jump point $\mathfrak j$, then $mp/q+\mathfrak j$ lies within $2|\eta_j|$ of one of the points $x_i-m_1\theta-s\theta$, $0\le s<M$, which fill an arc of length at most $M|\theta|$; so $mp/q+\mathfrak j$ lies in an arc of length at most $M|\theta|+4|\eta_j|$. If this length is at least $1$ there is nothing to prove: by \eqref{eqn.DriftBounds} we have $M|\theta|<2/q$, so either $q\le6$ or $4q|\eta_j|>2q/3$, and in both cases the claimed bound exceeds $M$. Otherwise, since $mp/q$ modulo $1$ depends only on the residue of $m$ and runs through $\frac1q\Z/\Z$ as the residue varies, the residue of $m$ belongs to a set of at most $q(M|\theta|+4|\eta_j|)+1<3+4q|\eta_j|$ residues, each of which occurs at most $\lceil M/q\rceil\le M/q+1$ times among $M$ consecutive integers. Finally, we multiply by two to account for both jump points.

(ii) The row-blocks contained in $[y,z]$ form a set of consecutive row-blocks, which is nonempty only if $z-y\ge h_j$. Besides these, $[y,z]$ meets at most two row-blocks, one at each end, in intervals whose lengths $\ell$ are less than $h_j$ and add up to at most $z-y$; we refer to these as partial row-blocks. On row-block $t$,
\[
G_i-G_{i+1}=\frac{1}{h_j}\big(R_i[t]-R_{i+1}[t]\big)-\frac{1}{h_j}\big(R^{\mr{rm}}_i[t]-R^{\mr{rm}}_{i+1}[t]\big),
\]
so $\int_y^z(G_i-G_{i+1})$ is the difference of the integrals over $[y,z]$ of the two step functions on the right, which we call the row term and the remainder term. Each of these integrals is a sum of contributions from the row-blocks which meet $[y,z]$: a row-block contained in $[y,z]$ contributes $R_i[t]-R_{i+1}[t]$ to the row term and $R^{\mr{rm}}_i[t]-R^{\mr{rm}}_{i+1}[t]$ to the remainder term, and a partial row-block contributes $\ell/h_j$ times these quantities.

\emph{The row term.} By Lemma \ref{lem.Adjacent}(i), $R_i[t]-R_{i+1}[t]$ is a sum over the rows of the row-block $t$ of terms of absolute value at most $2$ which vanish unless the phase of the row is within $2|\eta_j|$ of a jump point. The row-blocks contained in $[y,z]$, say those with $t_1\le t<t_2$, consist of the consecutive rows $t_1h_j\le m<t_2h_j$, and there are $(t_2-t_1)h_j\le z-y$ of these, so by (i) the sum of $R_i[t]-R_{i+1}[t]$ over these row-blocks has absolute value at most $4((z-y)/q+1)(3+4q|\eta_j|)$; and there are such row-blocks only when $z-y\ge h_j$. For a partial row-block, (i) with $M=h_j$ gives $|R_i[t]-R_{i+1}[t]|\le4(h_j/q+1)(3+4q|\eta_j|)$, so its contribution to the row term has absolute value at most
\[
\frac{\ell}{h_j}\cdot4\Big(\frac{h_j}{q}+1\Big)(3+4q|\eta_j|)=4\Big(\frac{\ell}{q}+\frac{\ell}{h_j}\Big)(3+4q|\eta_j|).
\]
Since the lengths $\ell$ add up to at most $z-y$ and to at most $2h_j$, the row term is $O\big(((z-y)/q+\min\{1,(z-y)/h_j\})(1+q|\eta_j|)\big)$.

\emph{The remainder term.} The row-blocks contained in $[y,z]$ are consecutive, so by (E1) the sum of $R^{\mr{rm}}_i[t]-R^{\mr{rm}}_{i+1}[t]$ over them has absolute value at most $N_i=O(1+q|\eta_j|)$, and again there are such row-blocks only when $z-y\ge h_j$. For a partial row-block, \eqref{eqn.RmDiff} and $|\iota_{r,t}|<1$ give $|R^{\mr{rm}}_i[t]-R^{\mr{rm}}_{i+1}[t]|<\#\{(r,\mathfrak j):rp/q+\mathfrak j\in A_{i,t}\}\le2(q|\eta_j|+1)$, since an arc of length $|\eta_j|$ contains at most $q|\eta_j|+1$ points of $\frac1q\Z/\Z$; so its contribution to the remainder term has absolute value at most $(\ell/h_j)\cdot2(q|\eta_j|+1)$. The remainder term is therefore $O\big(\min\{1,(z-y)/h_j\}(1+q|\eta_j|)\big)$, and the lemma follows.
\end{proof}

\section{Proof of Proposition \ref{prop.CellLemma}}\label{sec.CellLemma}

In this section we prove Proposition \ref{prop.CellLemma}. We fix $j\ge j_0$ and a cell $Q$ of level $j+1$, and we use the notation of Section \ref{sec.Density}. Since every quantity of Section \ref{sec.Density} depends on $Q$ only through its phases and block counts, we may translate $Q$ so that
\[
Q=[0,w_Q]\times[0,h_{j+1}],
\]
replacing the half-open cell by its closure; column $i$ then occupies $[s_i,s_i+w_i]\times[0,h_{j+1}]$ with $s_0=0$ and $s_{i+1}=s_i+w_i$, and the row-block $t$ occupies $[0,w_Q]\times[th_j,(t+1)h_j]$. We write $v(i,y)$ for the value of $v_{j+1}$ on column $i$ at height $0\le y<h_{j+1}$, which is given by \eqref{eqn.vFormula}, so that $v(i,\cdot)$ is a step function with jumps at the \textit{block points} $h_j,2h_j,\ldots,(T_j-1)h_j$, and by \eqref{eqn.RelDensityBound} we have $|v(i,y)-1|<24\varepsilon_j$.

We construct $\psi_Q$ as a composition of three maps. The first is obtained by gluing together vertical shears of the columns of $Q$. For each column $i\ge1$ we define a density $u^{(i)}$, the \textit{full-cycle density} of the column, which is obtained from $v(i,\cdot)$ by removing the remainder parts from the block counts of the column and renormalizing, and we apply to the column the vertical shear whose profile carries the measure $v(i,y)\,\mr{d}y$ to $u^{(i)}(y)\,\mr{d}y$ (Section \ref{subsec.ColumnMaps}). The merged column is sheared with the profile of column $1$, and its full-cycle density $u^{(0)}_\sharp$ is defined to be the density to which this shear carries $v(0,\cdot)$ (Section \ref{subsec.MergedColumn}). The mismatch between the shears of adjacent columns is $O(\varepsilon_jq_j)$ (Section \ref{subsec.EdgeMismatch}), and Lemma \ref{lem.StripGlue} glues them into a homeomorphism of $Q$ (Section \ref{subsec.RemainderGlue}). The second map is obtained in the same way from horizontal shears of horizontal strips of $Q$, and its Jacobian is the function $u^{\mr{fc}}$ which is equal to $u^{(i)}$ on column $i\ge1$ and to $u^{(0)}_\sharp$ on column $0$ (Section \ref{subsec.FullCycleMap}). The composition of these two maps is $(1+O(\varepsilon_j))$-BL and has Jacobian $v_{j+1}$, and composing it with a third map, provided by Lemma \ref{lem.BoundaryFix}, gives a map with the same Jacobian which fixes $\partial Q$ pointwise (Section \ref{subsec.Assembly}). Cells which consist of a single column, and cells in $\{x_1<0\}$, are treated in Section \ref{subsec.Remaining}. Several of the steps below require $\varepsilon_j$ to be smaller than an absolute constant; since $\varepsilon_j\le1/C_0$, the finitely many conditions of this kind all hold once $C_0$ is large enough, and $C_0$ is fixed accordingly in Section \ref{subsec.Assembly}.

\subsection{The full-cycle density and the column maps}\label{subsec.ColumnMaps}
For a column $1\le i\le N_Q$ and $0\le y<h_{j+1}$ put
\begin{equation}\label{eqn.uFC}
u^{(i)}(y)=\frac{c_i}{u_{j+1}(Q)}\Big(|W|+\frac{G_i(y)}{q_j}\Big),
\end{equation}
where the constant $c_i>0$ is chosen so that
\[
\int_0^{h_{j+1}}u^{(i)}=\int_0^{h_{j+1}}v(i,\cdot).
\]
By \eqref{eqn.vFormula}, $v(i,\cdot)$ is the same expression with $c_i$ replaced by $1$ and $G_i$ replaced by the full block count $R_i[t]/h_j$; so $u^{(i)}$ is the density which column $i$ would have if its remainder parts were removed from its block counts, rescaled so as to keep the mass of the column. Both $v(i,\cdot)$ and $u^{(i)}$ are step functions with jumps at the block points. Let
\[
F_i(y)=\int_0^yv(i,\cdot),\qquad F^{\mr{fc}}_i(y)=\int_0^yu^{(i)},\qquad P_i=(F^{\mr{fc}}_i)^{-1}\circ F_i\qquad(0\le y\le h_{j+1}).
\]
We call $P_i$ the \textit{column map} of column $i$: it is the increasing map of $[0,h_{j+1}]$ which carries the measure $v(i,y)\,\mr{d}y$ to the measure $u^{(i)}(y)\,\mr{d}y$.

\begin{lemma}\label{lem.ColumnMaps} 
Let $1\le i\le N_Q$.
\begin{enumerate}[topsep=2pt,itemsep=4pt,parsep=4pt]
\item[(i)] $|G_i|<4$, $c_i=1+O(\varepsilon_j)$, and $u^{(i)}=1+O(\varepsilon_j)$ on $[0,h_{j+1})$.
\item[(ii)] $P_i$ is an increasing piecewise linear homeomorphism of $[0,h_{j+1}]$ onto itself with $|P_i'-1|=O(\varepsilon_j)$ wherever it is differentiable; in particular it is $(1+O(\varepsilon_j))$-BL, and
\[
u^{(i)}\big(P_i(y)\big)\,P_i'(y)=v(i,y)\qquad\text{for almost every}\quad y .
\]
\item[(iii)] $|P_i(y)-y|=O(\varepsilon_jq)$ for all $0\le y\le h_{j+1}$.
\end{enumerate}
\end{lemma}
\begin{proof}
(i) We have $|R_i[t]|<3h_j$ by Lemma \ref{lem.Block} and $|R^{\mr{rm}}_i[t]|\le h_j$ by Lemma \ref{lem.Split}(iv), so $|G_i|<4$ and $|W|+G_i/q_j=|W|(1+O(\varepsilon_j))$. Let
\[
M_i=\int_0^{h_{j+1}}\Big(|W|+\frac{G_i}{q_j}\Big)=|W|h_{j+1}+\frac{1}{q_j}\sum_{t<T_j}\big(R_i[t]-R^{\mr{rm}}_i[t]\big)
\]
be the mass of the unnormalized full-cycle density of column $i$, so that $M_i=|W|h_{j+1}(1+O(\varepsilon_j))$. Integrating \eqref{eqn.uFC} and \eqref{eqn.vFormula} over $[0,h_{j+1}]$, the defining property of $c_i$ reads
\[
c_iM_i=|W|h_{j+1}+\frac{1}{q_j}\sum_{t<T_j}R_i[t]=M_i+\frac{1}{q_j}\sum_{t<T_j}R^{\mr{rm}}_i[t],
\]
and therefore
\begin{equation}\label{eqn.cMinusOne}
	c_i-1=\frac{\sum_{t<T_j}R^{\mr{rm}}_i[t]}{q_j\,M_i}.
\end{equation}
Again by Lemma \ref{lem.Split}(iv), we have that $|\sum_tR^{\mr{rm}}_i[t]|\le h_{j+1}$. Therefore \[|c_i-1|=O(h_{j+1}/(q_j|W|h_{j+1}))=O(\varepsilon_j),\] and then $u^{(i)}=1+O(\varepsilon_j)$ follows from \eqref{eqn.uFC}, since $u_{j+1}(Q)=|W|(1+O(\varepsilon_j))$ by \eqref{eqn.vFormula}.

(ii) $F_i$ and $F^{\mr{fc}}_i$ are increasing and continuous, they are linear between consecutive block points, with slopes $1+O(\varepsilon_j)$, and they have the same values at $0$ and at $h_{j+1}$. Hence $P_i$ is an increasing piecewise linear homeomorphism of $[0,h_{j+1}]$ onto itself; differentiating $F^{\mr{fc}}_i\circ P_i=F_i$ at a point $y$ at which $P_i$ is differentiable and $P_i(y)$ is not a block point, which excludes finitely many points, gives $u^{(i)}(P_i(y))P_i'(y)=v(i,y)$, and so $P_i'=1+O(\varepsilon_j)$. The BL constant of an increasing homeomorphism of an interval is the larger of the supremum of its derivative and the supremum of the derivative of its inverse, so $P_i$ is $(1+O(\varepsilon_j))$-BL.

(iii) By the definitions, we have that
\[
\int_y^{P_i(y)}u^{(i)}=F^{\mr{fc}}_i(P_i(y))-F^{\mr{fc}}_i(y)=F_i(y)-F^{\mr{fc}}_i(y).
\]
By (i) we have that $u^{(i)}\ge1/2$ once $\varepsilon_j$ is small enough, so the left hand side has absolute value at least $\frac12|P_i(y)-y|$, and therefore $|P_i(y)-y|\le2|F_i(y)-F^{\mr{fc}}_i(y)|$.

For a block point $y=t_0h_j$, we claim that
\begin{equation}\label{eqn.FirstOrder}
u_{j+1}(Q)\,q_j\,\big(F_i(y)-F^{\mr{fc}}_i(y)\big)=\sum_{t<t_0}R^{\mr{rm}}_i[t]-\Big(\sum_{t<T_j}R^{\mr{rm}}_i[t]\Big)\,\frac{F^{\mr{fc}}_i(y)}{F^{\mr{fc}}_i(h_{j+1})}\ .
\end{equation}
To see this, note that \[u_{j+1}(Q)F_i(y)=|W|y+q_j^{-1}\sum_{t<t_0}R_i[t]\] and \[u_{j+1}(Q)F^{\mr{fc}}_i(y)=c_i(|W|y+q_j^{-1}\sum_{t<t_0}(R_i[t]-R^{\mr{rm}}_i[t])).\] The expression in parentheses in the second formula is a constant multiple of $F^{\mr{fc}}_i(y)$, and its value at $y=h_{j+1}$ is $M_i$, so it is equal to $M_iF^{\mr{fc}}_i(y)/F^{\mr{fc}}_i(h_{j+1})$. Subtracting the two formulas, we find that
\[
u_{j+1}(Q)\big(F_i(y)-F^{\mr{fc}}_i(y)\big)=\frac{1}{q_j}\sum_{t<t_0}R^{\mr{rm}}_i[t]-(c_i-1)M_i\,\frac{F^{\mr{fc}}_i(y)}{F^{\mr{fc}}_i(h_{j+1})},
\]
and \eqref{eqn.FirstOrder} follows from \eqref{eqn.cMinusOne}.

Both sums on the right of \eqref{eqn.FirstOrder} have absolute value less than $q+7$ by (E3), and \[0\le F^{\mr{fc}}_i(y)/F^{\mr{fc}}_i(h_{j+1})\le1\] because $F^{\mr{fc}}_i$ is increasing with $F^{\mr{fc}}_i(0)=0$, so $|F_i(y)-F^{\mr{fc}}_i(y)|<2(q+7)/(u_{j+1}(Q)q_j)=O(\varepsilon_jq)$ at the block points. Between consecutive block points both $F_i$ and $F^{\mr{fc}}_i$ are linear, so the bound holds everywhere.
\end{proof}

We will need a version of \eqref{eqn.FirstOrder} which holds at every height, and not only at the block points. For $1\le i\le N_Q$ let $\Gamma_i$ be the continuous function on $[0,h_{j+1}]$ which is linear on each of the intervals $[th_j,(t+1)h_j]$, $0\le t<T_j$, and which satisfies
\[
\Gamma_i(t_0h_j)=\sum_{t<t_0}R^{\mr{rm}}_i[t]\qquad\text{for}\quad 0\le t_0\le T_j .
\]
Then, for every $0\le y\le h_{j+1}$,
\begin{equation}\label{eqn.FirstOrderExt}
	u_{j+1}(Q)\,q_j\,\big(F_i(y)-F^{\mr{fc}}_i(y)\big)=\Gamma_i(y)-\Gamma_i(h_{j+1})\,\frac{F^{\mr{fc}}_i(y)}{F^{\mr{fc}}_i(h_{j+1})} .
\end{equation}
To see this, note that both sides are linear on each of these intervals, because $F_i$ and $F^{\mr{fc}}_i$ are, and that they agree at the endpoints of each interval, by \eqref{eqn.FirstOrder} at the block points and because both sides vanish at $y=0$ and at $y=h_{j+1}$.

The displacement $O(\varepsilon_jq)$ in Lemma \ref{lem.ColumnMaps}(iii) is not small when $q$ is large: it can be much larger than the bound $O(\varepsilon_jq_j)$ which we need for the mismatch between adjacent columns. For this reason the proof of the next lemma has to compare the full-cycle densities of adjacent columns over intervals of heights of length $O(\varepsilon_jq)$.

\subsection{The mismatch at an interior edge}\label{subsec.EdgeMismatch}
For $1\le i<N_Q$ let
\[
\mu_i=\sup_{0\le y\le h_{j+1}}\big|P_i(y)-P_{i+1}(y)\big|
\]
be the mismatch, in the sense defined after Lemma \ref{lem.StripShear}, between the column maps of the columns $i$ and $i+1$.

\begin{lemma}\label{lem.EdgeMismatch} 
Let $1\le i<N_Q$.
\begin{enumerate}[topsep=2pt,itemsep=4pt,parsep=4pt]
\item[(i)] $|c_i-c_{i+1}|=O\big(\varepsilon_jq_j(1/h_{j+1}+\varepsilon_j/q)\big)$.
\item[(ii)] For $0\le y\le z\le h_{j+1}$,
\[
\Big|\int_y^z\big(u^{(i)}-u^{(i+1)}\big)\Big|=O\Big(\varepsilon_jq_j\Big(\frac{z-y}{q}+\min\Big\{1,\frac{z-y}{h_j}\Big\}\Big)\Big).
\]
\item[(iii)] $\mu_i=O(\varepsilon_jq_j)$.
\end{enumerate}
\end{lemma}
\begin{proof}
Throughout we use that $\Gamma_i(h_{j+1})=\sum_{t<T_j}R^{\mr{rm}}_i[t]$ is the total remainder count of column $i$, so that $|\Gamma_i(h_{j+1})|\le\min\{h_{j+1},q+7\}$ by Lemma \ref{lem.Split}(iv) and (E3), and \[|\Gamma_i(h_{j+1})-\Gamma_{i+1}(h_{j+1})|\le N_i=O(q_j)\] by (E1). We also use repeatedly that $1+q|\eta_j|=O(q_j)$, by \eqref{eqn.DriftBounds}, so that Lemma \ref{lem.FullCycleDiff}(ii) reads
\begin{equation}\label{eqn.GDiff}
\Big|\int_y^z\big(G_i-G_{i+1}\big)\Big|=O\Big(q_j\Big(\frac{z-y}{q}+\min\Big\{1,\frac{z-y}{h_j}\Big\}\Big)\Big)\qquad(0\le y\le z\le h_{j+1}).
\end{equation}

(i) By \eqref{eqn.cMinusOne},
\[
c_i-c_{i+1}=\frac{\Gamma_i(h_{j+1})-\Gamma_{i+1}(h_{j+1})}{q_jM_i}+\frac{\Gamma_{i+1}(h_{j+1})}{q_j}\cdot\frac{M_{i+1}-M_i}{M_iM_{i+1}} .
\]
The first term is $O(q_j/(q_j|W|h_{j+1}))=O(\varepsilon_jq_j/h_{j+1})$. For the second, we have that \[q_j(M_{i+1}-M_i)=\sum_t(R_{i+1}[t]-R_i[t])-(\Gamma_{i+1}(h_{j+1})-\Gamma_i(h_{j+1})),\] and by Lemma \ref{lem.Adjacent}(i) the first sum is a sum over all $h_{j+1}$ rows of column $i$ of terms of absolute value at most $2$ which vanish unless the phase of the row is within $2|\eta_j|$ of a jump point; by Lemma \ref{lem.FullCycleDiff}(i) with $M=h_{j+1}$ there are $O((h_{j+1}/q+1)q_j)$ such rows. So $|M_{i+1}-M_i|=O(h_{j+1}/q+1)$, and the second term is at most \[|\Gamma_{i+1}(h_{j+1})|q_j^{-1}\cdot O(h_{j+1}/q+1)\cdot(|W|h_{j+1})^{-2}=O\big(\varepsilon_j^2q_j(1/q+1/h_{j+1})\big),\] using $|\Gamma_{i+1}(h_{j+1})|\le h_{j+1}$. Since $\varepsilon_j\le1$ the two terms together give (i).

(ii) By \eqref{eqn.uFC},
\begin{equation}\label{eqn.uDiffSplit}
u^{(i)}-u^{(i+1)}=\frac{1}{u_{j+1}(Q)}\Big[(c_i-c_{i+1})\Big(|W|+\frac{G_i}{q_j}\Big)+\frac{c_{i+1}}{q_j}\big(G_i-G_{i+1}\big)\Big].
\end{equation}
The first term in the bracket is $O(|W||c_i-c_{i+1}|)$ by Lemma \ref{lem.ColumnMaps}(i), so over $[y,z]$ it contributes $O\big(\varepsilon_jq_j(z-y)(1/h_{j+1}+\varepsilon_j/q)\big)$ by (i), and $(z-y)/h_{j+1}\le\min\{1,(z-y)/h_j\}$. The second contributes \[O(1/(|W|q_j))\cdot O(q_j((z-y)/q+\min\{1,(z-y)/h_j\}))\] by \eqref{eqn.GDiff}, and $1/(|W|q_j)=\varepsilon_j$.

(iii) Fix $y$, and put $z_i=P_i(y)$, $z_{i+1}=P_{i+1}(y)$, and
\[
X=\big(F_i-F^{\mr{fc}}_i\big)(y)-\big(F_{i+1}-F^{\mr{fc}}_{i+1}\big)(y).
\]
Since $F^{\mr{fc}}_{i'}(z_{i'})=F_{i'}(y)$ for $i'\in\{i,i+1\}$ and $F^{\mr{fc}}_i(z_{i+1})=F^{\mr{fc}}_{i+1}(z_{i+1})+\int_0^{z_{i+1}}(u^{(i)}-u^{(i+1)})$, we have
\begin{equation}\label{eqn.MeanValue}
	\int_{z_{i+1}}^{z_i}u^{(i)}=F_i(y)-F^{\mr{fc}}_i(z_{i+1})=X-\int_y^{z_{i+1}}\big(u^{(i)}-u^{(i+1)}\big),
\end{equation}
and the left side has absolute value at least $\frac12|z_i-z_{i+1}|$. The integral on the right is taken over an interval of length $|z_{i+1}-y|=O(\varepsilon_jq)$, by Lemma \ref{lem.ColumnMaps}(iii), so by (ii) it is $O(\varepsilon_jq_j(\varepsilon_j+1))=O(\varepsilon_jq_j)$. It remains to bound $X$. Applying \eqref{eqn.FirstOrderExt} to column $i$ and to column $i+1$, subtracting the second identity from the first, and writing $\rho_{i'}(y)=F^{\mr{fc}}_{i'}(y)/F^{\mr{fc}}_{i'}(h_{j+1})\in[0,1]$, we have that
\begin{align*}
	u_{j+1}(Q)\,q_j\,X&=\Gamma_i(y)-\Gamma_i(h_{j+1})\rho_i(y)-\Gamma_{i+1}(y)+\Gamma_{i+1}(h_{j+1})\rho_{i+1}(y)\\
	&=\big(\Gamma_i-\Gamma_{i+1}\big)(y)-\big(\Gamma_i(h_{j+1})-\Gamma_{i+1}(h_{j+1})\big)\rho_i(y)\\
	&\hspace*{20bp}-\Gamma_{i+1}(h_{j+1})\big(\rho_i-\rho_{i+1}\big)(y).
\end{align*}
The first term on the right hand side is $(\Gamma_i-\Gamma_{i+1})(y)$. This function is linear on each of the intervals $[t_0h_j,(t_0+1)h_j]$, and its values at the endpoints of such an interval are the sums of $R^{\mr{rm}}_i[t]-R^{\mr{rm}}_{i+1}[t]$ over $t<t_0$ and over $t\le t_0$, each of which has absolute value at most $N_i$ by (E1). So the first term has absolute value at most $N_i=O(q_j)$. The same is true of the second term, by (E1) applied to the set of all row-blocks, and because $0\le\rho_i(y)\le1$. For the third, $\rho_i-\rho_{i+1}$ is the sum of $(F^{\mr{fc}}_i(y)-F^{\mr{fc}}_{i+1}(y))/F^{\mr{fc}}_i(h_{j+1})$ and of $\rho_{i+1}(y)(F^{\mr{fc}}_{i+1}(h_{j+1})-F^{\mr{fc}}_i(h_{j+1}))/F^{\mr{fc}}_i(h_{j+1})$; both numerators are integrals of $u^{(i)}-u^{(i+1)}$ over intervals of length at most $h_{j+1}$, hence $O(\varepsilon_jq_j(h_{j+1}/q+1))$ by (ii), and $F^{\mr{fc}}_i(h_{j+1})\ge h_{j+1}/2$. So $|\rho_i-\rho_{i+1}|=O(\varepsilon_jq_j(1/q+1/h_{j+1}))$, and since $|\Gamma_{i+1}(h_{j+1})|\le\min\{h_{j+1},q+7\}$, which is at most $8\min\{h_{j+1},q\}$, the third term is $O(\varepsilon_jq_j)$. Altogether the right hand side is $O(q_j)$, so by \eqref{eqn.vFormula} we have that $X=O(1/|W|)=O(\varepsilon_jq_j)$, and \eqref{eqn.MeanValue} gives $|z_i-z_{i+1}|=O(\varepsilon_jq_j)$.
\end{proof}

\subsection{The merged column}\label{subsec.MergedColumn}
Column $0$ is wider than the other columns, and its sub-cell counts contain the terms $R_0'[t]$ coming from the strip of width $q_{j-1}$, to which the decomposition of Lemma \ref{lem.Split} does not apply. For both of these reasons the comparison of adjacent columns in Lemmas \ref{lem.FullCycleDiff} and \ref{lem.EdgeMismatch} is not available for columns $0$ and $1$, and we do not try to define the full-cycle density of column $0$ from its own counts. Instead we shear column $0$ by the column map $P_1$ of its neighbor. Some shear of column $0$ is needed, because $P_1$ may move points by much more than the mismatch which Lemma \ref{lem.StripGlue} can remove at the edge between the two columns, and with this choice there is no mismatch at that edge. We then define the full-cycle density of column $0$ to be whatever this shear produces. Since column $0$ is a single column, of width less than $2q_j$, the trivial bounds for its counts are all that the horizontal shears of Section \ref{subsec.FullCycleMap} will require of this density. Put
\begin{equation}\label{eqn.Transported}
u^{(0)}_\sharp(y)=u^{(1)}(y)\,\frac{v\big(0,P_1^{-1}(y)\big)}{v\big(1,P_1^{-1}(y)\big)}\qquad(0\le y<h_{j+1}) .
\end{equation}
The images $P_1(h_j),P_1(2h_j),\ldots,P_1((T_j-1)h_j)$ of the block points under $P_1$ are called the \textit{shifted points}.

\begin{lemma}\label{lem.ColumnZero} 
The function $u^{(0)}_\sharp$ has the following properties.
\begin{enumerate}[topsep=2pt,itemsep=4pt,parsep=4pt]
\item[(i)] $u^{(0)}_\sharp\big(P_1(y)\big)\,P_1'(y)=v(0,y)$ for almost every $y$, and consequently $\int_0^{h_{j+1}}u^{(0)}_\sharp=\int_0^{h_{j+1}}v(0,\cdot)$.
\item[(ii)] $u^{(0)}_\sharp=1+O(\varepsilon_j)$ on $[0,h_{j+1})$.
\item[(iii)] $u^{(0)}_\sharp$ is a step function whose jumps lie at the block points and at the shifted points.
\end{enumerate}
\end{lemma}
\begin{proof}
(i) By Lemma \ref{lem.ColumnMaps}(ii) with $i=1$, $u^{(1)}(P_1(y))P_1'(y)=v(1,y)$ for almost every $y$, and multiplying by $v(0,y)/v(1,y)$ gives the identity. Integrating it over $[0,h_{j+1}]$, with the substitution $z=P_1(y)$, which is legitimate because $P_1$ is piecewise linear and increasing, gives the equality of masses.

(ii) By \eqref{eqn.RelDensityBound}, $v(0,y)/v(1,y)=1+O(\varepsilon_j)$, and $u^{(1)}=1+O(\varepsilon_j)$ by Lemma \ref{lem.ColumnMaps}(i).

(iii) The factor $u^{(1)}$ is constant between consecutive block points, and the ratio $v(0,\cdot)/v(1,\cdot)$, which is constant between consecutive block points, is composed with $P_1^{-1}$, so it is constant between consecutive shifted points.
\end{proof}

By (i), the shear of column $0$ by $P_1$ carries the measure $v(0,y)\,\mr{d}y$ to the measure $u^{(0)}_\sharp(y)\,\mr{d}y$, and $u^{(0)}_\sharp$, unlike the $u^{(i)}$, is not constant on the row-blocks. For this reason the horizontal shears of Section \ref{subsec.FullCycleMap} are applied on a decomposition of $Q$ which is finer than the row-blocks: the strips on which the full-cycle density is constant are cut by the block points and by the shifted points together. Let
\[
0=\zeta_0<\zeta_1<\cdots<\zeta_K=h_{j+1}
\]
be the distinct elements of $\{0,h_{j+1}\}\cup\{th_j:0<t<T_j\}\cup\{P_1(th_j):0<t<T_j\}$; we call $\zeta_1,\ldots,\zeta_{K-1}$ the \textit{refinement points}, and we put
\[
S_k=[0,w_Q]\times[\zeta_{k-1},\zeta_k]\qquad(1\le k\le K) .
\]
Since $P_1'\ge1-O(\varepsilon_j)$, consecutive shifted points are at least $\frac34h_j$ apart once $\varepsilon_j$ is small enough, and every shifted point is at least $\frac34h_j$ from $0$ and from $h_{j+1}$; block points are $h_j$ apart. So a shifted point and a block point can be close to each other, but two points of the same kind cannot, and we say that two distinct refinement points, one a block point and the other a shifted point, are \textit{partners} if their distance is less than $h_j/3$ (a refinement point which is both a block point and a shifted point therefore has no partner). A refinement point has at most one partner. Put
\[
\delta=\min\Big\{\frac{h_j}{6},\frac{w_Q}{20}\Big\},
\]
so that $q_j/20\le\delta\le h_j/6$; the first inequality holds because $h_j\ge q_j$ and $w_Q\ge q_{j+1}\ge q_j$.

\begin{lemma}\label{lem.Refinement}
To each refinement point $\zeta_k$, $0<k<K$, one can assign a closed interval $Z_k\subseteq[0,h_{j+1}]$ of length $\delta$ having $\zeta_k$ as an endpoint, in such a way that no refinement point lies in the interior of any $Z_k$, and the intervals $Z_1,\ldots,Z_{K-1}$ have pairwise disjoint interiors.
\end{lemma}
\begin{proof}
If $\zeta_k$ has a partner, let $Z_k$ be the interval of length $\delta$ with endpoint $\zeta_k$ which lies on the side of $\zeta_k$ away from its partner; otherwise let $Z_k=[\zeta_k-\delta,\zeta_k]$. Every refinement point other than $\zeta_k$ and its partner is at distance at least $h_j/3>\delta$ from $\zeta_k$, and the partner lies on the other side, so no refinement point lies in the interior of $Z_k$; and $Z_k\subseteq[0,h_{j+1}]$ because $\zeta_k$ is at distance at least $\frac34h_j>\delta$ from both ends of $[0,h_{j+1}]$. If $\zeta_k$ and $\zeta_{k'}$ are partners, their intervals lie on opposite outer sides of the segment between them; otherwise $|\zeta_k-\zeta_{k'}|\ge h_j/3\ge2\delta$, and both intervals lie within $\delta$ of their own points. In either case the interiors are disjoint.
\end{proof}

\subsection{Gluing the vertical shears}\label{subsec.RemainderGlue}
We now glue the vertical shears of the columns. Let $\varphi^{\mr{rm}}_i(x,y)=(x,P_i(y))$ on column $i$ for $1\le i\le N_Q$, and let $\varphi^{\mr{rm}}_0(x,y)=(x,P_1(y))$ on column $0$; so columns $0$ and $1$ carry the same shear, and we treat their union $[0,s_1+q_j]\times[0,h_{j+1}]$ as a single strip. Put
\[
u^{\mr{fc}}(x,y)=u^{(i)}(y)\quad\text{on column }i\ge1,\qquad u^{\mr{fc}}(x,y)=u^{(0)}_\sharp(y)\quad\text{on column }0,
\]
the \textit{full-cycle density} of $Q$.

\begin{lemma}\label{lem.RemainderGlue} 
Suppose that $\varepsilon_j$ is small enough. There is a $(1+O(\varepsilon_j))$-BL homeomorphism $\tilde\varphi^{\mr{rm}}$ of $Q$ onto itself with the following properties.
\begin{enumerate}[topsep=2pt,itemsep=4pt,parsep=4pt]
\item[(i)] $\tilde\varphi^{\mr{rm}}$ fixes the bottom and top sides of $Q$ pointwise, acts on the left side of $Q$ as $P_1$ and on the right side as $P_{N_Q}$.
\item[(ii)] If $\varphi^\circ$ is a BL map on $Q$ with $|\det D\varphi^\circ(z)|=u^{\mr{fc}}(z)$ for almost every $z\in Q$, then for every $0\le i\le N_Q$,
\[
|\det D(\varphi^\circ\circ\tilde\varphi^{\mr{rm}})|=|\det D(\varphi^\circ\circ\varphi^{\mr{rm}}_i)|\qquad\text{almost everywhere on column }i .
\]
\end{enumerate}
\end{lemma}
\begin{proof}
We apply Lemma \ref{lem.StripGlue} to $Q$, cut into the strips $[0,s_1+q_j]\times[0,h_{j+1}]$ (columns $0$ and $1$ together) and, for $2\le i\le N_Q$, the columns $[s_i,s_i+q_j]\times[0,h_{j+1}]$, with $H=h_{j+1}$. The strip maps are $\varphi^{\mr{rm}}_1$ on the first strip and $\varphi^{\mr{rm}}_i$ on column $i\ge2$. By Lemmas \ref{lem.ColumnMaps}(ii) and \ref{lem.Shear}, each of them is a $(1+O(\varepsilon_j))$-BL shear of its strip onto itself, with profile $P_i$ ($P_1$ for the first strip); so it fixes the bottom and top sides pointwise and acts on both vertical sides as its profile, an increasing $(1+O(\varepsilon_j))$-BL homeomorphism of $[0,h_{j+1}]$. The mismatch at the edge $x=s_{i+1}$ between the strips with maps $\varphi^{\mr{rm}}_i$ and $\varphi^{\mr{rm}}_{i+1}$ is $\mu_i=O(\varepsilon_jq_j)$ by Lemma \ref{lem.EdgeMismatch}(iii), for $1\le i<N_Q$; there are no other edges. At that edge we place the zone
\[
Z^{\mr{rm}}_i=\Big[s_{i+1}-\frac{q_j}{5},s_{i+1}\Big]\times[0,h_{j+1}]
\]
inside column $i$, of width $q_j/5$. The zones lie in distinct columns, so their interiors are disjoint; and when an interior edge exists we have $N_Q\ge2$, hence $a_{j+1}\ge2$ and $h_{j+1}\ge q_{j+1}\ge2q_j$, so the width $q_j/5$ is at most $h_{j+1}/10$, as Lemma \ref{lem.StripGlue}(ii) requires. The weight of the first strip is $u^{(1)}$, and that of column $i\ge2$ is $u^{(i)}$; all are $1+O(\varepsilon_j)$. The hypothesis $\mu_i\le\varepsilon\delta_i$ of Lemma \ref{lem.StripGlue} holds with $\delta_i=q_j/5$ and $\varepsilon=C'\varepsilon_j$ for an absolute constant $C'$, and $\varepsilon\le\varepsilon_\ast$ once $\varepsilon_j$ is small enough. Lemma \ref{lem.StripGlue} now gives a $(1+O(\varepsilon_j))$-BL homeomorphism $\tilde\varphi^{\mr{rm}}$ of $Q$ with the boundary behavior (i), by its part (c).

For (ii), let $\varphi^\circ$ be as stated. The zone contained in the first strip lies in column $1$, where $u^{\mr{fc}}(z)=u^{(1)}(z_2)$ is the weight of that strip, and on column $i\ge2$ we have $u^{\mr{fc}}(z)=u^{(i)}(z_2)$, the weight of that column; column $N_Q$ contains no zone. So the last assertion of Lemma \ref{lem.StripGlue}(b) applies on every strip, and it gives $|\det D(\varphi^\circ\circ\tilde\varphi^{\mr{rm}})|=|\det D(\varphi^\circ\circ\varphi^{\mr{rm}}_i)|$ almost everywhere on column $i$ for $i\ge1$, and the same with $\varphi^{\mr{rm}}_1=\varphi^{\mr{rm}}_0$ on column $0$.
\end{proof}

Nothing has been said about the Jacobian of $\tilde\varphi^{\mr{rm}}$ itself, and nothing needs to be: it is composed with the map of Section \ref{subsec.FullCycleMap}, whose Jacobian on the zones is the weight, and (ii) says that the composite does not see the zone maps.

\subsection{The full-cycle map}\label{subsec.FullCycleMap}
We now realize $u^{\mr{fc}}$ by horizontal shears. By Lemma \ref{lem.ColumnZero}(iii) and the definition of the $u^{(i)}$, the function $u^{\mr{fc}}$ is constant on the interior of each product of a column with one of the strips $S_k$; in particular on $S_k$ it depends on $x$ only, and we write $u^{\mr{fc}}_k(x)$ for its value there, with the mean $\bar u^{\mr{fc}}_k=w_Q^{-1}\int_0^{w_Q}u^{\mr{fc}}_k$. By Lemmas \ref{lem.ColumnMaps}(i) and \ref{lem.ColumnZero}(ii), $u^{\mr{fc}}=1+O(\varepsilon_j)$. Also, $\int_Qu^{\mr{fc}}=|Q|$, because on each column the mass of $u^{\mr{fc}}$ equals the mass of $v_{j+1}$ (by the choice of $c_i$ for $i\ge1$ and by Lemma \ref{lem.ColumnZero}(i) for $i=0$) and $v_{j+1}$ has mean $1$ over $Q$. So Lemma \ref{lem.StripShear}, with the coordinates exchanged, applies to $u^{\mr{fc}}$ and the horizontal strips $S_k$. In its notation, put
\[
P^{\mr{fc}}_\perp(y)=\int_0^y\bar u^{\mr{fc}},\qquad\varphi^{\mr{fc}}_\perp(x,y)=\big(x,P^{\mr{fc}}_\perp(y)\big),\qquad P^{\mr{fc}}_k(x)=\frac{1}{\bar u^{\mr{fc}}_k}\int_0^xu^{\mr{fc}}_k,\qquad\varphi^{\mr{fc}}_k(x,y)=\big(P^{\mr{fc}}_k(x),y\big),
\]
where $\bar u^{\mr{fc}}(y)=\bar u^{\mr{fc}}_k$ for $\zeta_{k-1}\le y<\zeta_k$.

Then $\varphi^{\mr{fc}}_\perp$ is a $(1+O(\varepsilon_j))$-BL homeomorphism of $Q$ onto itself which fixes the bottom and top sides of $Q$ pointwise, acts on the vertical sides as $P^{\mr{fc}}_\perp$, and maps $S_k$ onto the strip $S_k'=[0,w_Q]\times[P^{\mr{fc}}_\perp(\zeta_{k-1}),P^{\mr{fc}}_\perp(\zeta_k)]$, whose height is $\bar u^{\mr{fc}}_k(\zeta_k-\zeta_{k-1})$. The map $\varphi^{\mr{fc}}_k$ is a $(1+O(\varepsilon_j))$-BL homeomorphism of $S_k'$ onto itself which fixes its vertical sides pointwise, acts on its bottom and top sides as $P^{\mr{fc}}_k$, and satisfies $|\det D(\varphi^{\mr{fc}}_k\circ\varphi^{\mr{fc}}_\perp)|=u^{\mr{fc}}$ almost everywhere on $S_k$. For $0<k<K$ the edge between $S_k'$ and $S_{k+1}'$ is the image under $\varphi^{\mr{fc}}_\perp$ of the line at height $\zeta_k$, and the mismatch at this edge is
\[\mu^{\mr{fc}}_k=\sup_x|P^{\mr{fc}}_{k+1}(x)-P^{\mr{fc}}_k(x)|.\]
The displacements are $|P^{\mr{fc}}_\perp(y)-y|\le\int_0^{h_{j+1}}|\bar u^{\mr{fc}}-1|=O(\varepsilon_jh_{j+1})$ and $|P^{\mr{fc}}_k(x)-x|\le\int_0^{w_Q}|u^{\mr{fc}}_k/\bar u^{\mr{fc}}_k-1|=O(\varepsilon_jw_Q)$.

\begin{lemma}\label{lem.FullCycleMap}
Suppose that $\varepsilon_j$ is small enough.
\begin{enumerate}[topsep=2pt,itemsep=4pt,parsep=4pt]
\item[(i)] For $0<k<K$, $\mu^{\mr{fc}}_k=O(\varepsilon_jq_j)$.
\item[(ii)] There is a $(1+O(\varepsilon_j))$-BL homeomorphism $\varphi^{\mr{fc}}$ of $Q$ onto itself with $|\det D\varphi^{\mr{fc}}|=u^{\mr{fc}}$ almost everywhere, which maps each side of $Q$ onto itself, acting on the vertical sides as $P^{\mr{fc}}_\perp$, on the bottom side as $P^{\mr{fc}}_1$, and on the top side as $P^{\mr{fc}}_K$.
\end{enumerate}
\end{lemma}
\begin{proof}
(i) We first show that for two measurable functions $f_-,f_+:[0,w_Q]\rar[1/2,2]$ with means $\bar f_\pm$ and profiles $P_\pm(x)=\bar f_\pm^{-1}\int_0^xf_\pm$,
\begin{equation}\label{eqn.ProfileDiff}
\sup_{0\le x\le w_Q}|P_+(x)-P_-(x)|\le4\int_0^{w_Q}|f_+-f_-| .
\end{equation}
This holds because
\[P_+(x)-P_-(x)=\bar f_+^{-1}\int_0^x(f_+-f_-)-P_-(x)(\bar f_+-\bar f_-)/\bar f_+,\]
where $0\le P_-\le w_Q$ and $w_Q|\bar f_+-\bar f_-|\le\int_0^{w_Q}|f_+-f_-|$.

We apply \eqref{eqn.ProfileDiff} with $f_-=u^{\mr{fc}}_k$ and $f_+=u^{\mr{fc}}_{k+1}$, the values of $u^{\mr{fc}}$ on the two sides of the line at height $\zeta_k$, for which $P_-=P^{\mr{fc}}_k$ and $P_+=P^{\mr{fc}}_{k+1}$. This gives
\[
\mu^{\mr{fc}}_k\le4\int_0^{w_Q}\big|u^{\mr{fc}}_{k+1}-u^{\mr{fc}}_k\big| ,
\]
so it suffices to show that the integral on the right is $O(\varepsilon_jq_j)$. On column $0$ the integrand is the difference of two values of $u^{(0)}_\sharp$, which is $O(\varepsilon_j)$ by Lemma \ref{lem.ColumnZero}(ii), and column $0$ has width less than $2q_j$, so it contributes $O(\varepsilon_jq_j)$ to the integral. On a column $i\ge1$ the integrand is the absolute value of the difference of the values of $u^{(i)}$ on the two sides of $\zeta_k$. If $\zeta_k$ is not a block point then $u^{(i)}$ is constant near $\zeta_k$, and these columns contribute nothing. Suppose then that $\zeta_k=t_0h_j$ is a block point. By \eqref{eqn.uFC}, on column $i\ge1$ we have that
\[
u^{\mr{fc}}_{k+1}-u^{\mr{fc}}_k=\frac{c_i}{q_j\,u_{j+1}(Q)}\Big(G_i(t_0h_j)-G_i\big((t_0-1)h_j\big)\Big),
\]
where $G_i(t_0h_j)$ and $G_i((t_0-1)h_j)$ are the values of $G_i$ on the row-blocks $t_0$ and $t_0-1$. By the definition of $G_i$ and Lemma \ref{lem.Split}(i), $h_jG_i(th_j)=R^{\mr{av}}_i[t]+R^{\mr{cr}}_i[t]+R^{\epsilon}_i[t]$ for every row-block $t$, so the three bounds of (E2) give
\[
\sum_{i=1}^{N_Q}\Big|G_i(t_0h_j)-G_i\big((t_0-1)h_j\big)\Big|=O(1).
\]
Since the columns $i\ge1$ have width $q_j$, and $c_i=O(1)$ by Lemma \ref{lem.ColumnMaps}(i), it follows that
\[
\int_{s_1}^{w_Q}\big|u^{\mr{fc}}_{k+1}-u^{\mr{fc}}_k\big|=\frac{1}{u_{j+1}(Q)}\sum_{i=1}^{N_Q}c_i\,\Big|G_i(t_0h_j)-G_i\big((t_0-1)h_j\big)\Big|=O\Big(\frac{1}{|W|}\Big)=O(\varepsilon_jq_j),
\]
by \eqref{eqn.vFormula}. This completes the proof of (i).

(ii) We apply Lemma \ref{lem.StripGlue}, with the coordinates exchanged, to the strips $S_1',\ldots,S_K'$, which are stacked vertically in $Q$ and have length $w_Q$ (so that $w_Q$ plays the part of $H$), with the strip maps $\varphi^{\mr{fc}}_k$, the weights $1$, and the following zones. For $0<k<K$ let $Z^{\mr{fc}}_k=\varphi^{\mr{fc}}_\perp([0,w_Q]\times Z_k)$, where $Z_k$ is the interval of Lemma \ref{lem.Refinement}. Since $Z_k$ has $\zeta_k$ as an endpoint and no refinement point in its interior, $[0,w_Q]\times Z_k$ lies in $S_k$ or in $S_{k+1}$, so $Z^{\mr{fc}}_k$ is a rectangle of full length $w_Q$ with one horizontal side on the edge between $S_k'$ and $S_{k+1}'$, lying in one of them; its height is $\bar u^{\mr{fc}}\delta$ for the mean of the strip containing it, which lies between $\delta/2$ and $2\delta$; and the zones have pairwise disjoint interiors, because the intervals $Z_k$ do and $\varphi^{\mr{fc}}_\perp$ is a homeomorphism which maps horizontal lines to horizontal lines. The heights are at most $2\delta\le w_Q/10$, as required, and at least $\delta/2\ge q_j/40$.

By (i) the mismatch at each edge is $O(\varepsilon_jq_j)$, which is at most $\varepsilon$ times the height of the zone at that edge for $\varepsilon=C''\varepsilon_j$ with an absolute constant $C''$, and $\varepsilon\le\varepsilon_\ast$ once $\varepsilon_j$ is small enough. Lemma \ref{lem.StripGlue} gives a $(1+O(\varepsilon_j))$-BL homeomorphism $\tilde\varphi^{\mr{fc}}$ of $Q$ onto itself which, with the coordinates exchanged in its part (c), fixes the left and right sides of $Q$ pointwise, acts on the bottom side as $P^{\mr{fc}}_1$ and on the top side as $P^{\mr{fc}}_K$, and has, by its part (b) with the weights $1$, Jacobian $|\det D\varphi^{\mr{fc}}_k|$ almost everywhere on $S_k'$. Put $\varphi^{\mr{fc}}=\tilde\varphi^{\mr{fc}}\circ\varphi^{\mr{fc}}_\perp$. By the chain rule, $|\det D\varphi^{\mr{fc}}|=|\det D(\varphi^{\mr{fc}}_k\circ\varphi^{\mr{fc}}_\perp)|=u^{\mr{fc}}$ almost everywhere on $S_k$, and the side maps are as stated because $\varphi^{\mr{fc}}_\perp$ fixes the bottom and top sides pointwise and acts on the vertical sides as $P^{\mr{fc}}_\perp$.
\end{proof}

\subsection{Assembly, the Jacobian, and the boundary}\label{subsec.Assembly}
We can now prove Proposition \ref{prop.CellLemma} for a cell with $N_Q\ge1$. Let
\[
\varphi=\varphi^{\mr{fc}}\circ\tilde\varphi^{\mr{rm}},
\]
a $(1+O(\varepsilon_j))$-BL homeomorphism of $Q$ onto itself. We first check that $|\det D\varphi|=v_{j+1}$ almost everywhere on $Q$. By Lemma \ref{lem.FullCycleMap}(ii), $\varphi^{\mr{fc}}$ satisfies the hypothesis on $\varphi^\circ$ in Lemma \ref{lem.RemainderGlue}(ii), so on column $i\ge1$ we have $|\det D\varphi|=|\det D(\varphi^{\mr{fc}}\circ\varphi^{\mr{rm}}_i)|$ almost everywhere, and by the chain rule and Lemma \ref{lem.ColumnMaps}(ii),
\[
|\det D(\varphi^{\mr{fc}}\circ\varphi^{\mr{rm}}_i)(x,y)|=u^{\mr{fc}}\big(x,P_i(y)\big)\,P_i'(y)=u^{(i)}\big(P_i(y)\big)\,P_i'(y)=v(i,y)
\]
for almost every $(x,y)$ in column $i$ (the point $(x,P_i(y))$ lies in column $i$, where $u^{\mr{fc}}=u^{(i)}$). On column $0$ the same computation with $\varphi^{\mr{rm}}_0$, whose profile is $P_1$, and Lemma \ref{lem.ColumnZero}(i) gives $|\det D\varphi(x,y)|=u^{(0)}_\sharp(P_1(y))P_1'(y)=v(0,y)$. The column edges form a null set, so $|\det D\varphi|=v_{j+1}$ almost everywhere.

Next we record the behavior of $\varphi$ on $\partial Q$. By Lemma \ref{lem.RemainderGlue}(i) and Lemma \ref{lem.FullCycleMap}(ii), the map $\varphi$ maps each side of $Q$ onto itself, acting on the left side as $P^{\mr{fc}}_\perp\circ P_1$, on the right side as $P^{\mr{fc}}_\perp\circ P_{N_Q}$, on the bottom side as $P^{\mr{fc}}_1$, and on the top side as $P^{\mr{fc}}_K$. Each of these is an increasing $(1+O(\varepsilon_j))$-BL homeomorphism of its side, being a composition of such maps. For the displacements, Lemma \ref{lem.ColumnMaps}(iii) and the displacement bounds recorded before Lemma \ref{lem.FullCycleMap} give, for $i\in\{1,N_Q\}$,
\[
\big|P^{\mr{fc}}_\perp(P_i(y))-y\big|\le\big|P^{\mr{fc}}_\perp(P_i(y))-P_i(y)\big|+|P_i(y)-y|=O(\varepsilon_jh_{j+1})+O(\varepsilon_jq),
\]
while the bottom and top side maps satisfy $|P^{\mr{fc}}_k(x)-x|=O(\varepsilon_jw_Q)$. Since $q\le q_{j+1}$, $h_{j+1}<2q_{j+1}$, and $w_Q<2q_{j+1}$, while $\min\{w_Q,h_{j+1}\}\ge q_{j+1}$, every side map of $\varphi$ moves points by at most an absolute multiple of $\varepsilon_j\min\{w_Q,h_{j+1}\}$. So Lemma \ref{lem.BoundaryFix} applies to $\varphi$, with $\varepsilon$ an absolute multiple of $\varepsilon_j$, once $\varepsilon_j$ is small enough, and it produces a $(1+O(\varepsilon_j))$-BL homeomorphism $\psi_Q=\hat\varphi$ of $Q$ onto itself which fixes $\partial Q$ pointwise and has $|\det D\psi_Q|=|\det D\varphi|=v_{j+1}$ almost everywhere.

Finally we fix the constants. The construction of this section used $\varepsilon_j\le c$ for finitely many absolute constants $c$ (in Lemmas \ref{lem.ColumnMaps}-\ref{lem.FullCycleMap}, in the separation of the refinement points, and in the choice of $\varepsilon$ when Lemmas \ref{lem.StripGlue} and \ref{lem.BoundaryFix} were applied, each of which requires an absolute multiple of $\varepsilon_j$ to be at most $\varepsilon_\ast$). We choose $C_0\ge12$ so large that all of these conditions hold whenever $\varepsilon_j\le1/C_0$. With this $C_0$, and for $j\ge j_0$, every cell of level $j+1$ in $\{x_1\ge0\}$ with $N_Q\ge1$ carries a map $\psi_Q$ as required, which is $(1+C\varepsilon_j)$-BL for an absolute constant $C$.

\subsection{Single-column cells, and the cells in \texorpdfstring{$\{x_1<0\}$}{the left half-plane}}\label{subsec.Remaining}
For a cell with $N_Q=0$, we have $a_{j+1}=1$ and $w_Q=q_{j+1}=q_j+q_{j-1}$: the cell is its merged column, cut into $T_j$ row-blocks, and $v_{j+1}$ depends on the height $y$ only. Put $P(y)=\int_0^yv(0,\cdot)$ and $\varphi(x,y)=(x,P(y))$. Since $v(0,\cdot)$ takes values in $[1-24\varepsilon_j,1+24\varepsilon_j]$ and has mean $1$ over $[0,h_{j+1}]$, Lemma \ref{lem.Shear} shows that $\varphi$ is a $(1+48\varepsilon_j)$-BL shear of $Q$ with profile $P$ and $|\det D\varphi|=v_{j+1}$ almost everywhere; so it fixes the bottom and top sides of $Q$ pointwise and acts on the vertical sides as $P$, which moves points by at most $24\varepsilon_jh_{j+1}<48\varepsilon_j\min\{w_Q,h_{j+1}\}$, since $h_{j+1}<2q_{j+1}=2w_Q$. Lemma \ref{lem.BoundaryFix}, applied with $\varepsilon=48\varepsilon_j$, then gives $\psi_Q$ as required, with a constant $1+O(\varepsilon_j)$, after enlarging $C_0$ and $C$ by absolute factors if necessary.

Finally let $Q$ be a cell in $\{x_1<0\}$, translated as above to $[0,w_Q]\times[0,h_{j+1}]$, so that its merged column is its last column, and let $\varrho(x,y)=(w_Q-x,y)$ be the reflection of $Q$ in its vertical midline. If the columns of $Q$ are numbered from right to left, and their block counts and phases are defined from their left edges as before, then $v_{j+1}\circ\varrho$ is given by \eqref{eqn.vFormula}, and $x_{i+1}=x_i+\eta_j$ for $0\le i<N_Q$, in place of \eqref{eqn.PhaseStep}. The arguments of Sections \ref{sec.Density} and \ref{sec.CellLemma} apply with the sign of this step reversed, so Sections \ref{subsec.ColumnMaps}-\ref{subsec.Assembly} give a $(1+C\varepsilon_j)$-BL homeomorphism $\psi'$ of $Q$ onto itself which fixes $\partial Q$ pointwise and has $|\det D\psi'|=v_{j+1}\circ\varrho$ almost everywhere. Since $\varrho$ is an isometry of $Q$, the map $\psi_Q=\varrho\circ\psi'\circ\varrho$ has the required properties. This completes the proof of Proposition \ref{prop.CellLemma}.

\section{From the Jacobian to the bijection}\label{sec.Bottom}

In this section we complete the passage from Proposition \ref{prop.Telescope} to a bijection with $\Z^2$. The window is now allowed to be closed or half-open: let $\alpha_1\in(0,1)$ be irrational, let $\alpha_2\in\R$, let $W\subseteq\R$ be an interval with $0<|W|<1$, and let $W^\circ$ be its interior. We write
\[
\Lambda^\circ=\Lambda_\balph(W^\circ)\subseteq\Lambda=\Lambda_\balph(W),
\]
and we apply Sections \ref{sec.Hierarchy}-\ref{sec.CellLemma} to the open window $W^\circ$, which has the same length as $W$: the parameters $\varepsilon_j$ and $j_0$, the cells, and the map $\Phi$ of Proposition \ref{prop.Telescope} are those of $\Lambda^\circ$. As noted in Section \ref{sec.Reduction}, the set $\Lambda\setminus\Lambda^\circ$ has at most two points on each horizontal line $m_2=k$. It is convenient to normalize the density to $1$: put
\[
\Phi_W=|W|^{1/2}\,\Phi\qquad\text{and}\qquad L_W=L_0|W|^{-1/2},
\]
so that $\Phi_W$ is an $L_W$-BL homeomorphism of $\R^2$ (we use $|W|<1$), and by \eqref{eqn.PhiCount}
\begin{equation}\label{eqn.PhiWCount}
|\Phi_W(E)|=\#(\Lambda^\circ\cap E)\qquad\text{for every finite union $E$ of cells of level $j_0$}.
\end{equation}
Since $q_{j_0}|W|\ge C_0\ge1$ we have $L_W\le L_0q_{j_0}^{1/2}$, so everything below depends on $\balph$ and $W$ only through $q_{j_0}$. A \textit{unit square} is a set $[m_1,m_1+1)\times[m_2,m_2+1)$ with $(m_1,m_2)\in\Z^2$, and for a finite union $U$ of unit squares $|\partial U|$ is the length of its boundary.

\noindent\emph{The count estimate.} 
The next lemma converts \eqref{eqn.PhiWCount} into a count estimate on unions of unit squares.

\begin{lemma}\label{lem.Count} 
For every finite union $U$ of unit squares,
\[
\Big|\#\big(\Phi_W(\Lambda^\circ)\cap U\big)-|U|\Big|\le4(L_W+2q_{j_0})^2\,|\partial U|
\]
and
\[
\Big|\#\big(\Phi_W(\Lambda)\cap U\big)-|U|\Big|\le5(L_W+2q_{j_0})^2\,|\partial U| .
\]
\end{lemma}
\begin{proof}
\noindent\emph{Reduction to boundary cells.} Call a closed edge of the unit grid \textit{separating} if it is a side of a unit square contained in $U$ and of a unit square not contained in $U$. It is easy to check that $\partial U$ is the union of the separating edges, so that $|\partial U|$ is their number. Put $\Omega=\Phi_W^{-1}(U)$, a bounded set. The cells $Q$ of level $j_0$ tile $\R^2$, so the sets $\Phi_W(Q\cap\Omega)$ partition $U$, and
\[
\#\big(\Phi_W(\Lambda^\circ)\cap U\big)-|U|=\sum_Q\Big(\#(\Lambda^\circ\cap Q\cap\Omega)-|\Phi_W(Q\cap\Omega)|\Big),
\]
a finite sum. By \eqref{eqn.PhiWCount} with $E=Q$, both terms of the summand lie between $0$ and $\#(\Lambda^\circ\cap Q)$, and they are equal if $Q\subseteq\Omega$, while both vanish if $Q\cap\Omega=\emptyset$. We call the remaining cells, those which meet both $\Omega$ and its complement, the \textit{boundary cells}, and we write $N_U$ for their number. A boundary cell is connected, so it meets $\partial\Omega=\Phi_W^{-1}(\partial U)$, because $\Phi_W$ is a homeomorphism. A cell of level $j_0$ has width at most $2q_{j_0}$ and height $h_{j_0}<2q_{j_0}$, so it contains fewer than $4q_{j_0}^2$ points of $\Z^2$, and the display above gives
\begin{equation}\label{eqn.BoundaryCells}
\Big|\#\big(\Phi_W(\Lambda^\circ)\cap U\big)-|U|\Big|\le4q_{j_0}^2N_U .
\end{equation}

\noindent\emph{Counting the boundary cells.} Each boundary cell meets $\Phi_W^{-1}(e)$ for some separating edge $e$, a set of diameter at most $L_W$. The cells of level $j_0$ are the products of the intervals of $\mc{P}_{j_0}$, which have length at least $q_{j_0}$, with the intervals $[th_{j_0},(t+1)h_{j_0})$, of length $h_{j_0}\ge q_{j_0}$. An interval of length at most $D$ meets at most $D/q_{j_0}+2$ intervals of a partition of $\R$ into intervals of length at least $q_{j_0}$, since all but the first and the last of those it meets are contained in it; so a set of diameter at most $D$, whose two coordinate projections have length at most $D$, meets at most $(D/q_{j_0}+2)^2$ cells. Therefore $N_U\le|\partial U|(L_W/q_{j_0}+2)^2$, and \eqref{eqn.BoundaryCells} gives the first bound of the lemma.

\noindent\emph{The boundary of the window.} Since $\Lambda\setminus\Lambda^\circ$ has at most two points on each horizontal line $m_2=k$,
\[
0\le\#\big(\Phi_W(\Lambda)\cap U\big)-\#\big(\Phi_W(\Lambda^\circ)\cap U\big)\le2\,\#\{k\in\Z:(\R\times\{k\})\cap\Omega\not=\emptyset\},
\]
and it remains to count the lines which meet $\Omega$. A horizontal line which meets the bounded set $\Omega$ also meets its complement, so, being connected, it meets $\partial\Omega$, and hence it meets $\Phi_W^{-1}(e)$ for some separating edge $e$. This set has diameter at most $L_W$, so it meets at most $L_W+1$ of the lines $\R\times\{k\}$, $k\in\Z$. Therefore at most $(L_W+1)|\partial U|$ of these lines meet $\Omega$, the two counts differ by at most $2(L_W+1)|\partial U|$, and $4(L_W+2q_{j_0})^2+2(L_W+1)\le5(L_W+2q_{j_0})^2$, because $L_W\ge1$.
\end{proof}

\noindent\emph{Laczkovich's criterion.}
Lemma \ref{lem.Count} is the hypothesis of a theorem of Laczkovich \cite[Theorem 1.1]{Lacz1992}, which characterizes, by count estimates of this kind, the discrete sets in $\R^d$ which admit a bijection of bounded displacement with a multiple of $\Z^d$. We use the following lemma from his proof.

\begin{lemma}\label{lem.Laczkovich}
Let $Y\subseteq\R^2$ be a set whose intersection with every bounded set is finite, and suppose that there is a constant $c$ such that
\[
\Big|\#(Y\cap U)-|U|\Big|\le c\,|\partial U|
\]
for every finite union $U$ of unit squares. Then there is a bijection $f:Y\rar\Z^2$ with $|f(y)-y|\le2\sqrt2\,(4c+1)$ for all $y\in Y$.
\end{lemma}

This is \cite[Lemma 2.3]{Lacz1992}, applied with $S_1=Y$, $S_2=\Z^2$, $\alpha=1$, $C_1=c$, and $C_2=0$: every half-open unit square contains exactly one point of $\Z^2$, so $\#(\Z^2\cap U)=|U|$ for every finite union $U$ of unit squares. 

We can now prove the following proposition, which is Theorem \ref{thm.Plane} in the case when $\alpha_1\in(0,1)$ is irrational, with a BL constant which depends only on $q_{j_0}$.

\begin{proposition}\label{prop.Bottom}
Let $\alpha_1\in(0,1)$ be irrational, let $\alpha_2\in\R$, and let $W\subseteq\R$ be an interval, open, closed, or half-open, with $0<|W|<1$. Then there is an $L$-BL bijection from $\Lambda_\balph(W)$ onto $\Z^2$ with $L=O(q_{j_0}^3)$, where $q_{j_0}$ is the least denominator of a convergent of $\alpha_1$ with $q_{j_0}|W|\ge C_0$.
\end{proposition}
\begin{proof}
Since $\Lambda\subseteq\Z^2$ and $\Phi_W$ is $L_W$-BL, the set $Y=\Phi_W(\Lambda)$ has separation at least $1/L_W$, and by Lemma \ref{lem.Count} it satisfies the hypotheses of Lemma \ref{lem.Laczkovich} with $c=5(L_W+2q_{j_0})^2$, so there is a bijection $f:Y\rar\Z^2$ which moves no point by more than $K=2\sqrt2(4c+1)$. Since $\Z^2$ has separation $1$, for $y\not=y'$ in $Y$ the triangle inequality gives
\begin{align*}
|f(y)-f(y')|&\le|y-y'|+2K\le(1+2KL_W)|y-y'| ,\\
|y-y'|&\le|f(y)-f(y')|+2K\le(1+2K)|f(y)-f(y')| ,
\end{align*}
so that $f$ is $(1+2KL_W)$-BL (this is the argument of \cite[Proposition 1.8]{SmilSolo2026}). The restriction of $\Phi_W$ to $\Lambda$ is an $L_W$-BL bijection from $\Lambda$ onto $Y$, so $f\circ\Phi_W|_\Lambda$ is an $L_W(1+2KL_W)$-BL bijection from $\Lambda$ onto $\Z^2$. Finally $L_W\le L_0q_{j_0}^{1/2}$ with $L_0$ absolute, so $c=O(q_{j_0}^2)$, $K=O(q_{j_0}^2)$, and $L_W(1+2KL_W)=O(q_{j_0}^3)$.
\end{proof}

\section{Proofs of the theorems, and remarks}\label{sec.Proofs}

\noindent\emph{Proof of Theorem \ref{thm.Plane}:}
Let $\balph\in\R^2\setminus\Q^2$ and let $W\subseteq\R$ be an interval with $0<|W|<1$. By the first two normalizations of Section \ref{sec.Reduction}, which replace $\Lambda_\balph(W)$ by an isometric copy of itself, we may assume that $\alpha_1\in(0,1)$ is irrational. Proposition \ref{prop.Bottom} then gives a BL bijection from $\Lambda_\balph(W)$ onto $\Z^2$. This completes the proof of Theorem \ref{thm.Plane}, and with it, by the arguments of Section \ref{sec.Reduction}, the proofs of Theorem \ref{thm.Main} and Corollary \ref{cor.CPS}.

\noindent\emph{Finite unions of intervals.}
The companion paper \cite{Hayn2026} shows that the conclusion of Theorem \ref{thm.Main} fails for a compact window which is the closure of its interior and has infinitely many components, so the shape of the window matters. It was left open there whether the conclusion can fail for a window which is a finite union of intervals. When $\alpha_1$ satisfies the Diophantine condition (1.1) of \cite{BuraKlei2002} this is ruled out by the proof of \cite[Theorem 1.2]{BuraKlei2002}, which applies to such windows (see the introduction of \cite{Hayn2026}). The following theorem shows that it cannot happen for any $\balph$.

\begin{theorem}\label{thm.Union}
Let $\balph\in\R^2\setminus\Q^2$, and let $W\subseteq\R$ be a finite union of bounded intervals (each open, closed, or half-open), with $|W|>0$. Then $Y_\balph(W)$ is BL equivalent to $\Z^2$. If $W$ is contained in an interval of length less than $1$, then $\Lambda_\balph(W)$ is also BL equivalent to $\Z^2$.
\end{theorem}

The proof is the same as that of Theorems \ref{thm.Main} and \ref{thm.Plane}, with the changes which we now describe.

\noindent\emph{Proof of Theorem \ref{thm.Union}:}
Suppose first that $W$ is contained in an interval of length less than $1$. For the statement about $\Lambda_\balph(W)$ we may assume, as before, that $\alpha_1\in(0,1)$ is irrational. Let $W_1,\ldots,W_\nu$ be the components of the interior $W^\circ$ of $W$, which are open intervals whose images in $\R/\Z$ are pairwise disjoint. Then the discrepancy $R_M(\gamma)$ of $W^\circ$ and the block counts of Section \ref{sec.Density} are the sums of the corresponding quantities for $W_1,\ldots,W_\nu$, and we define the four parts of the block counts of $W^\circ$ in the same way, with a single choice of $p/q$ in \eqref{eqn.Dirichlet}. Lemma \ref{lem.Block} and Lemmas \ref{lem.SubCell}-\ref{lem.FullCycleDiff} hold for every open interval of length less than $1$ (the condition $j\ge j_0$ is not used in their proofs), so they apply to each of the intervals $W_1,\ldots,W_\nu$. By the triangle inequality, every estimate for these quantities which is used in Sections \ref{sec.Hierarchy} and \ref{sec.CellLemma} therefore holds for $W^\circ$, with its right hand side multiplied by $\nu$. After division by $\nu$ these are the estimates which were used for a single interval, with $|W|/\nu$ in place of $|W|$, and dividing the number of points of $\Lambda_\balph(W^\circ)$ in every cell by $\nu$ does not change the relative densities $v_{j+1}$. It follows that Propositions \ref{prop.CellLemma} and \ref{prop.Telescope} hold for $W^\circ$, with $\varepsilon_j=\nu/(q_j|W|)$ and with $j_0$ the least index for which $q_{j_0}|W|\ge C_0\nu$. In Section \ref{sec.Bottom} the set $\Lambda\setminus\Lambda^\circ$ now has at most $\#\partial W$ points on each horizontal line, instead of at most two, and this changes only the constants in Lemma \ref{lem.Count} and Proposition \ref{prop.Bottom}. Therefore $\Lambda_\balph(W)$ is BL equivalent to $\Z^2$.

For the statement about $Y_\balph(W)$, let $I$ be an interval with $0<|I|<1$ which contains $W$. Then $Y_\balph(W)\subseteq Y_\balph(I)$ is uniformly discrete, by Lemma \ref{lem.NormalForm}, and \cite[Lemma 4]{Hayn2026} gives a bijection $\Psi:Y_\balph(W)\rar\Lambda_\balph(W)$ as in Lemma \ref{lem.NormalForm}. The rest of the argument is the same as in the proof of Theorem \ref{thm.Main}.

Finally, let $W$ be any set as in the theorem, and write $\bn=\bn_\balph$. Then $Y_\balph(W)$ is the set $Y$ of Corollary \ref{cor.CPS} with $x=0$ and with $\sigma$ replaced by the set $\{-s\bn/|\bn|^2:s\in W\}$, no two points of which differ by a nonzero integer vector, since $V_\balph$ is not rational. Step 3 of the proof of Corollary \ref{cor.CPS} applies to this set, with $\xi$ the difference of the endpoints of the smallest segment which contains it: the window $W'=W/(N-\alpha_1)$ which it produces (with $N$ the integer chosen there) is a finite union of intervals, contained in an interval of length less than $1$, and the case of the theorem which we have already proved takes the place of Theorem \ref{thm.Main}. This completes the proof of Theorem \ref{thm.Union}.

We close with a remark on what these results leave open.

\noindent\emph{Higher dimensions.} 
Most of our arguments make sense for cut-and-project sets of dimension $d\ge3$ with one-dimensional internal space and an interval window, that is, for the sets $\Lambda_\balph(W)=\{m\in\Z^d:\langle\balph,m\rangle\in W+\Z\}$ with $\balph\in\R^d$: the reduction of Section \ref{sec.Reduction}, Lemma \ref{lem.Block}, the hierarchy of Section \ref{sec.Hierarchy}, and the final step of Section \ref{sec.Bottom}, which rests on Laczkovich's theorem in $\R^d$, all carry over. The proof of Proposition \ref{prop.CellLemma} does not. At least two of its ingredients are specific to the plane. The first is the estimate (E1) of Lemma \ref{lem.Estimates}, which relies on the fact that every interval of integers contains each residue class modulo $q$ the expected number of times, up to an error less than $1$. For $d\ge3$ the rows are indexed by $\Z^{d-1}$ and the role of the residue classes is played by the cosets of a sublattice of $\Z^{d-1}$, for which the corresponding error is not bounded. The second is the zone lemma (Lemma \ref{lem.Zone}), whose sliding construction produces a map of a planar zone with prescribed trace and Jacobian exactly $1$, and whose $d$-dimensional version we do not have. It seems that new ideas would be needed to overcome these difficulties. Furthermore, the case of cut-and-project sets with internal spaces of dimension greater than $1$ seems to be out of reach of the methods of this paper, since Lemma \ref{lem.Block} and the hierarchy of Section \ref{sec.Hierarchy} depend on the continued fraction expansion of a single real number.

\vspace{.15in}
{\footnotesize
\noindent
Department of Mathematics\\
University of Houston\\
Houston, TX, United States\\
haynes@math.uh.edu
}

\end{document}